\documentclass[twoside,11pt]{article}
\DeclareMathSizes{20}{30}{16}{12}
\usepackage{float}
\usepackage{amsmath}
\usepackage{mathtools}
\usepackage{mathrsfs}
\usepackage{stmaryrd}
\usepackage{bbm}
\usepackage{yfonts}
\usepackage{amsfonts}
\usepackage{tcolorbox}
\usepackage{tikz-cd}
\usepackage[utf8]{inputenc}

\usepackage{caption}
\usepackage{indentfirst}
\usepackage[Symbol]{upgreek}
\usepackage{enumerate}   
\usepackage{upref}

\usepackage{graphicx}
\usepackage[margin=0.8 in]{geometry}
\usepackage{fancyhdr}
\usepackage{hyperref}
\hypersetup{
    colorlinks,
    citecolor=black,
    filecolor=black,
    linkcolor=black,
    urlcolor=black
}
\newcommand{\N}{\mathbb{N}}
\newcommand{\Z}{\mathbb{Z}}

\newcommand{\mbbF}{\mathbb{F}}
\newcommand{\ra}{\rightarrow}

\newcommand{\calA}{\mathcal{A}}
\newcommand{\calI}{\mathcal{I}}

\newcommand{\calE}{\mathcal{E}}

\newcommand{\calD}{\mathcal{D}}

\newcommand{\calB}{\mathcal{B}}

\newcommand{\calM}{\mathcal{M}}
\newcommand{\calN}{\mathcal{N}}
\newcommand{\calP}{\mathcal{P}}

\newcommand{\op}{\operatorname}

\newcommand{\Ox}{\mathcal{O}}
\newcommand{\F}{\mathcal{F}}

\newcommand{\ov}{\overline}

\newcommand{\frakX}{\mathfrak{X}}

\newcommand{\frakU}{\mathfrak{U}}

\newcommand{\frakS}{\mathfrak{S}}
\newcommand{\frakY}{\mathfrak{Y}}

\newcommand{\frakZ}{\mathfrak{Z}}

\newcommand{\isom}{\xrightarrow{\raisebox{-0.4ex}[0ex][0ex]{$\sim$}}}

\long\def\/*#1*/{}

\usepackage{unicode-math}
\usepackage{fontspec}

\usepackage{sectsty}
\sectionfont{\centering}

\usepackage[toc,page]{appendix}
\usepackage{amsthm}

\newtheorem*{proof*}{Proof}

\newtheorem{theorem}[subsubsection]{Theorem}
\newtheorem{theorem*}{Theorem}
\newtheorem{proposition}[subsubsection]{Proposition}
\newtheorem{corollaire}[subsubsection]{Corollary}
\newtheorem{lemma}[subsubsection]{Lemma}

\theoremstyle{definition}

\newtheorem{definition}[subsubsection]{Definition}
\newtheorem{remark}[subsubsection]{Remarks}

\newtheorem{exmp}[subsubsection]{Examples}

\theoremstyle{definition}
\newtheorem{parag}[subsubsection]{}

\usepackage{etoolbox}

\numberwithin{equation}{subsection} 

\makeatletter
\renewcommand{\theequation}{%
  \ifnum\value{subsubsection}>0
    \thesubsubsection.\arabic{equation}%
  \else
    \thesubsection.\arabic{equation}%
  \fi
}
\makeatother

\title{Lifting the local logarithmic Cartier transform modulo $p^n$}
\author{Sami Fersi}
\date{}

\AtEndDocument{\bigskip{\footnotesize%
  \textsc{Sami Fersi, Laboratoire Alexander Grothendieck, Institut des Hautes Études Scientifiques, 35 Route de Chartres, 91440 Bures-sur-Yvette, France} \par
  \textit{E-mail address}: \texttt{fersi@ihes.fr} \par
}}

\begin{document}

\maketitle

\begin{abstract}
Let $p$ be a prime number and $n\ge 1$ an integer. We construct a local logarithmic Cartier transform modulo $p^n$, following an approach by Shiho. More precisely, consider the ring $W$ of Witt vectors of a perfect field of characteristic $p$. Let $f:\mathfrak{X}\ra \mathfrak{S}$ be a log smooth morphism of fs logarithmic $p$-adic formal schemes, whose underlying formal schemes are flat and locally of finite type over $W$, and assume the exact relative Frobenius after reduction modulo $p$ admits a lift $F$. Let $\mathcal{A}$ be Lorenzon's indexed algebra on the reduction modulo $p^n$ of $\mathfrak{X}.$ We define an indexed subalgebra $\mathcal{B}_F\subset\mathcal{A}$ and equip it with an integrable quasi-nilpotent $p$-connection. Our main result is an equivalence between the category of indexed $\mathcal{B}_F$-modules equipped with integrable quasi-nilpotent $p$-connections compatible with the $p$-connection on $\mathcal{B}_F$, and the category of indexed $\mathcal{A}$-modules equipped with integrable quasi-nilpotent connections compatible with the canonical connection on $\mathcal{A}$. This gives a local logarithmic version of the Cartier transform modulo $p^n$.
\end{abstract}

\tableofcontents

\newpage
\section{Introduction}

\begingroup
\renewcommand{\thesubsection}{\thesection}

\begin{parag}
The Cartier transform of Ogus and Vologodsky \cite{OgusVol} is an equivalence, for a smooth scheme $X$ over a perfect field $k$ of positive characteristic, between modules on $X$ equipped with integrable connections and modules on $X'=X\times_{k,F_k}k$ equipped with Higgs fields, both satisfying certain nilpotence conditions. Shiho constructed, in \cite{Shiho}, a local version of this transform under the assumption that both $X$ and the relative Frobenius $X\ra X'$ admit lifts to $p$-adic formal schemes. A key ingredient in proving that Shiho's transform is an equivalence, is faithfully flat descent along the relative Frobenius morphism. Shiho's construction has the advantage of being valid modulo $p^n$ for all positive integers $n,$ and not just in characteristic $p.$ In \cite{SF1}, we generalized Shiho's local transform to the log smooth setting and we proved that the functor we obtain is fully faithful. It is however not necessarily essentially surjective due to the relative Frobenius not being in general flat in the log smooth case. This is the same issue that arises when attempting to generalize Cartier descent to the log smooth setting. In this latter case, one solution was developed by Lorenzon in \cite{Lor2000}. He introduced \emph{indexed algebras} and \emph{indexed modules} and proved an \emph{indexed Cartier descent theorem}. In the present article, we enrich our previous construction of Shiho's local logarithmic transform with an indexed refinement, obtaining a local logarithmic Cartier transform in characteristic $p^n.$ To prove that it is an equivalence of categories, we establish a faithfully flat descent theorem for indexed algebras. 
\end{parag}

\begin{parag}
We now describe our work in more details. Let $p$ be a prime number and $W$ the ring of Witt vectors of a perfect field of characteristic $p.$ For a morphism $X\ra S$ of fs logarithmic schemes of characteristic $p,$ we recall the diagram of the exact relative Frobenius:
$$
\begin{tikzcd}
X \ar[bend right=-30]{rrd}{F_X} \ar[bend right=30]{ddr} \ar{dr}{F_{X/S}} & & \\
 & X' \ar{r} \ar{d} & X \ar{d} \\
 & S \ar{r}{F_S} & S,
\end{tikzcd}
$$
where the square is cartesian in the category of fs logarithmic schemes and $F_S$ and $F_X$ are the absolute Frobenius morphisms of $S$ and $X$ respectively. The morphism $F_{X/S}$ is the \emph{exact relative Frobenius of $X$ with respect to $S$} (\cite{Ogus2018} IV 3.3.9).
\end{parag}

\begin{parag}\label{intro1}
Let $f:\frakX\ra \frakS$ be a log smooth morphism of fs logarithmic $p$-adic formal schemes, flat and locally of finite type over $W.$ Denote by $X\ra S$ its special fiber. Locally, the diagonal immersion $\frakX \ra \frakX\times_{\frakS}\frakX$ factors as
\begin{equation}\label{decomp}
\frakX \xrightarrow{\Delta} \frakY \xrightarrow{g} \frakX \times_{\frakS} \frakX,
\end{equation}
where $\Delta$ is a strict closed immersion and $g$ is log étale and separated (\ref{lem122}). Let $R_{\frakX/\frakS}$ be the dilatation of the ideal of $X\ra \frakX \xrightarrow{\Delta} \frakY$ with respect to $p$ (\cite{Ahmed2010} 3.2.3.4). Although the decomposition \eqref{decomp} is not unique, we prove in \ref{prop224x} and \ref{prop3211} that the formal scheme $R_{\frakX/\frakS}$ is independent of the choice of that decomposition and can hence be defined globally. We obtain a formal scheme $R_{\frakX/\frakS}$ over $\frakX\times_{\frakS}\frakX,$ along with a closed immersion $\frakX \ra R_{\frakX/\frakS}.$ We denote by $P_{\frakX/\frakS,1}$ its PD envelope and, for a positive integer $n,$ we denote by $P_n$ its reduction modulo $p^n$ and by $\calP_n$ the structural ring of $P_n.$ For all $n,$ $\calP_n$ has a Hopf algebra structure (\cite{DXU} 4.2) and we can hence consider stratifications with respect to $P_n:$
\begin{definition}[\cite{DXU} 5.4]
Let $n\ge 1$ be an integer and $\ov{\calI}_n$ the PD ideal of $\calP_n.$ A $\calP_n$-stratification on an $\Ox_{\frakX_n}$-module $\calE$ is an isomorphism $\epsilon:\calP_n \otimes_{\Ox_{\frakX_n}}\calE \isom \calE \otimes_{\Ox_{\frakX_n}} \calP_n,$ where the left hand side (resp. right hand side) $\calP_n$ is considered as an $\Ox_{\frakX_n}$-module via the second projection (resp. the first projection) $P_n \ra \frakX_n$, such that:
\begin{itemize}
\item The base change of $\epsilon$ by the canonical projection $\calP_n \ra \calP_n/\ov{\calI}_n=\Ox_{\frakX_n}$ is the identity.
\item $\epsilon$ satisfies a certain cocycle condition (\cite{DXU} 5.4.1).
\end{itemize}
\end{definition}

\begin{proposition}[\ref{equivstrat}]
For an integer $n\ge 1$ and an $\Ox_{\frakX_n}$-module $\calE,$ the two following data are equivalent:
\begin{enumerate}
\item An integrable quasi-nilpotent $p$-connection on $\calE$ (\ref{parKhaminei} and \ref{def227}).
\item A $\calP_n$-stratification on $\calE.$
\end{enumerate}
\end{proposition}
\end{parag}

\begin{parag}
In section 4, we review the notions of indexed algebras and modules introduced by Lorenzon in \cite{Lor2000} and establish some key results and definitions that we will need in later sections. For a ringed topos $(X,\Ox_X)$ and a sheaf of commutative monoids $\upmu$ of $X,$ we equip the localized topos $X_{/\upmu}$ with the ring $\Ox_{X,\upmu},$ localization of $\Ox_X$ at $\upmu.$ The projections $p_1,p_2:\upmu^2 \ra \upmu$ and the addition maps $\sigma:\upmu^2 \ra \upmu,$ $\tau:\upmu^3 \ra \upmu$ induce morphisms of ringed topoi, abusively denoted
$$
p_1,p_2,\sigma:X_{/\upmu^2} \ra X_{/\upmu},\quad \tau:X_{/\upmu^3} \ra X_{/\upmu}.
$$
A $\upmu$-indexed algebra is a module $\calA$ of $X_{/\upmu}$ equipped with a multiplication map $\calA \boxtimes \calA = p_1^*\calA \otimes_{\Ox_{X,\upmu^2}} p_2^*\calA \ra \sigma^*\calA,$ satisfying certain commutativity, associativity and unitarity conditions. Similarly, for a $\upmu$-indexed algebra $\calA,$ we define the notion of a $\upmu$-indexed $\calA$-module. We recall that the category of $\upmu$-indexed $\calA$-modules has a tensor product defined as follows: the inverse image functor $\sigma^*$ has a left adjoint $\sigma_!.$ If $\calM$ and $\calN$ are two indexed $\calA$-modules, the two actions of $\calA$ on both $\calM$ and $\calN,$ induce two morphisms
\begin{equation}\label{eqintro1}
\begin{tikzcd}
\tau_!\left ( \calA \boxtimes \calM \boxtimes \calN \right ) \ar[shift left=-0.5ex]{r} \ar[shift left=0.5ex]{r} & \sigma_! \left ( \calM \boxtimes \calN \right ).
\end{tikzcd}
\end{equation}
The tensor product $\calM \circledast_{\calA} \calN$ is defined to be the coequalizer of these two morphisms.
\end{parag}

\begin{parag}
The main example of an indexed algebra we will be working with is that coming from a logarithmic structure. Let $X$ be an fs logarithmic scheme and $\upmu=\ov{\calM}_X^{\mathrm{gp}}.$ Take the equivalence relation $\sim$ on $\calM_X^{\mathrm{gp}} \times \Ox_X$ defined by $(\alpha_X^{-1}(a)+m,b)\sim (m,ab)$ for all local sections $m,$ $a$ and $b$ of $\calM_X^{\mathrm{gp}},$ $\Ox_X^*$ and $\Ox_X$ respectively. The composition of the canonical projections
$$
\calM_X^{\mathrm{gp}} \times \Ox_X \ra \calM_X^{\mathrm{gp}} \ra \ov{\calM}_X^{\mathrm{gp}}=\upmu
$$
factors through the quotient $\calA_X=\calM_X^{\mathrm{gp}} \times \Ox_X/\sim,$ turning $\calA_X$ into an object of $X_{\text{ét}/\upmu}.$ For a local section $s$ of $\calM_X^{\mathrm{gp}},$ we have a specific section $e_s=(s,1)$ of $\calA_X.$ We then define a $\upmu$-indexed algebra structure on $\calA_X$ by setting, for all local sections $a$ and $b$ of $\Ox_X$ and $s$ and $t$ of $\calM_X^{\mathrm{gp}},$
$$
(ae_s)\cdot (be_t)=abe_{s+t}.
$$
If $X\ra S$ is a morphism of fs logarithmic schemes, then $\calA_X$ is equipped with a canonical connection $d_{\calA_X}:\calA_X \ra \calA_X \otimes_{\Ox_{X,\upmu}} \omega^1_{X/S,\upmu},$ satisfying, for all local sections $s$ of $\calM_X^{\mathrm{gp}},$
$$
d_{\calA_X}(e_s)=e_s\otimes \op{dlog}s.
$$
Its kernel $\calB_{X/S}$ inherits a $\upmu$-indexed algebra structure from $\calA_X.$ In \ref{thmlor}, we reprove a theorem of Lorenzon stating that $\calA_X$ is locally free of finite rank over $\calB_{X/S}.$ We also prove the following isomorphism:
\end{parag}

\begin{proposition}[\ref{th234}]
Let $X\ra S$ be a log smooth morphism of fs logarithmic schemes in characteristic $p.$ Denote by $F_{X/S}:X\ra X'$ the exact relative Frobenius and $X\times_{X'}X$ the fiber product in the category of fs log schemes. Then, the projections $X\times_{X'}X \ra X$ induce an isomorphism of indexed algebras
\begin{equation}
\calA_{X} \circledast_{\calB_{X/S}} \calA_X \isom \calA_{X\times_{X'}X}.
\end{equation}
\end{proposition}

\begin{parag}
In \cite{SF3}, we constructed an equivalence of categories between $\upmu$-indexed $\calB_{X/S}$-modules, equipped with a $\calB_{X/S}$-linear Higgs field, and $\upmu$-indexed $\calA_X$-modules equipped with integrable connections compatible with $d_{\calA_X},$ both satisfying certain nilpotence conditions. The equivalence is based on an Azumaya algebra property that was proved by Schepler and Ohkawa in \cite{Schepler} and \cite{Ohkawa} respectively. If we want to establish a similar equivalence modulo $p^n,$ we have to find a replacement for the Azumaya algebra approach. First, we need to find a suitable counterpart to $\calB_{X/S}$ modulo $p^n.$

Let $n\ge 1$ be an integer and $f:\widetilde{X} \ra \widetilde{S}$ a log smooth morphism of fs logarithmic schemes, flat over $\Z/p^n\Z.$ Let $X \ra S$ be its reduction modulo $p$ and suppose that the exact relative Frobenius $F_{X/S}:X\ra X'$ lifts to an $\widetilde{S}$-morphism of fs logarithmic schemes $\widetilde{F}:\widetilde{X}\ra \widetilde{X}'$ such that $\widetilde{X}'$ is log smooth over $\widetilde{S}.$ Set $\upmu=\ov{\calM}_X^{\mathrm{gp}}$ and we identify the étale sites of $X,$ $X',$ $\widetilde{X}$ and $\widetilde{X}'$ via the universal homeomorphisms $X\ra \widetilde{X},$ $X' \ra \widetilde{X}'$ and $\widetilde{F}:\widetilde{X} \ra \widetilde{X}'.$ 
\end{parag}

\begin{definition}
We define $\calB_{\widetilde{F}}$ to be the $\upmu$-indexed subalgebra of $\calA_{\widetilde{X}}$ generated locally by $\Ox_{\widetilde{X}'}$ and sections of the form $\alpha_{\widetilde{X}}(m)e_{-m}$ and $e_{\widetilde{F}^{\flat}m'},$ where $m$ and $m'$ are local sections of $\calM_{\widetilde{X}}$ and $\calM_{\widetilde{X}'}^{\mathrm{gp}}$ respectively.
\end{definition}

We prove, in \ref{corBequalizer}, that $\calB_{\widetilde{F}}$ can be defined alternatively as an equalizer:
\begin{proposition}[\ref{corBequalizer}]
Let $\widetilde{Y}=\widetilde{X}\times_{\widetilde{X}'}\widetilde{X}$ be the fiber product of $\widetilde{F}$ with itself in the category of fs logarithmic schemes. Consider the two morphisms $p_1,p_2:\calA_{\widetilde{X}} \ra \calA_{\widetilde{Y}}$ induced by the projections $\widetilde{Y} \ra \widetilde{X}.$ Then, the following sequence is exact:
$$
\begin{tikzcd}
0 \ar{r} & \calB_{\widetilde{F}} \ar{r} & \calA_{\widetilde{X}} \ar[swap, shift right=0.5ex]{r}{p_2} \ar[shift right=-0.5ex]{r}{p_1} & \calA_{\widetilde{Y}}.
\end{tikzcd}
$$
\end{proposition}

Just as in characteristic $p,$ we have the following local freeness result:
\begin{proposition}[\ref{propABlocfree}]\label{propintro2}
$\calA_{\widetilde{X}}$ is a locally free $\upmu$-indexed $\calB_{\widetilde{F}}$-module of rank $p^d,$ where $d$ is the rank of $\omega^1_{\widetilde{X}/\widetilde{S}}.$
\end{proposition}

\begin{theorem}[\ref{propnablaF}]
The $\upmu$-indexed algebra $\calB_{\widetilde{F}}$ has a unique connection $\nabla_{\widetilde{F}}:\calB_{\widetilde{F}} \ra \calB_{\widetilde{F}} \otimes_{\Ox_{\widetilde{X}',\upmu}} \omega^1_{\widetilde{X}'/\widetilde{S},\upmu}$ satisfying:
$$
\nabla_{\widetilde{F}}(a)=d'a,\quad \nabla_{\widetilde{F}}(\alpha_{\widetilde{X}}(m)e_{-m})=0,\quad \nabla_{\widetilde{F}}(e_{\widetilde{F}^{\flat}m'})=e_{\widetilde{F}^{\flat}m'} \otimes \op{dlog}m',
$$
$$
\nabla_{\widetilde{F}}(xy)= x\nabla_{\widetilde{F}}(y)+y\nabla_{\widetilde{F}}(x),
$$
for all local sections $a,$ $m$ and $m'$ of $\Ox_{\widetilde{X}'},$ $\calM_{\widetilde{X}}$ and $\calM_{\widetilde{X}'}$ respectively and for all local sections $x$ and $y$ of $\calB_{\widetilde{F}},$ where $d':\Ox_{\widetilde{X}'} \ra \omega^1_{\widetilde{X}'/\widetilde{S}}$ is the universal derivation. In addition $\nabla_{\widetilde{F}}$ is integrable. We denote by $\theta_{\widetilde{F}}$ the integrable quasi-nilpotent $p$-connection (\ref{def622}) on $\calB_{\widetilde{F}}$ given by $\theta_{\widetilde{F}}=p\nabla_{\widetilde{F}}.$
\end{theorem}

\begin{proposition}[\ref{thmkeyiso}]
Let $\widetilde{X}\times_{\widetilde{X}'}\widetilde{X}$ be the fiber product in the category of fs logarithmic schemes. The projections $\widetilde{X}\times_{\widetilde{X}'}\widetilde{X} \ra \widetilde{X}$ induce an isomorphism of $\upmu$-indexed algebras
\begin{equation}\label{eqintro4}
\calA_{\widetilde{X}} \circledast_{\calB_{\widetilde{F}}} \calA_{\widetilde{X}} \isom \calA_{\widetilde{X}\times_{\widetilde{X}'}\widetilde{X}}.
\end{equation}
\end{proposition}

\begin{parag}
In section 5, we establish a faithfully flat descent theorem for indexed algebras and modules. The definitions, results and proofs are very similar to those of the usual faithfully flat descent. Let $(X,\Ox_X)$ be a ringed topos and $\upmu$ a sheaf of commutative monoids of $X.$ A morphism of $\upmu$-indexed algebras $\calA \ra \calA'$ is said to be \emph{faithfully flat} if a sequence of $\upmu$-indexed $\calA$-modules is exact if and only if so is its base change by $\calA \ra \calA'.$ For a morphism of $\upmu$-indexed algebras $f:\calA\ra \calA'$ and $\upmu$-indexed $\calA$-module $\calE,$ we denote by $f^*\calE$ the base change $\calA'\circledast_{\calA}\calE.$ We have the two following descent theorems:

\begin{theorem}[\ref{inddescentmor}, Descent of indexed morphisms]\label{thmintro1}
Let $\calA \ra \calA'$ be a faithfully flat morphism of $\upmu$-indexed algebras and $\calE_1$ and $\calE_2$ two $\upmu$-indexed $\calA$-modules. Set $\calA''=\calA'\circledast_{\calA}\calA'$ and, for $i=1,2,$ $\calE_i'=\calE_i\circledast_{\calA}\calA'$ and $\calE_i''=\calE'_i\circledast_{\calA}\calA'$ and denote by $\pi_1,\pi_2:\calA' \ra \calA'\circledast_{\calA} \calA'$ the morphisms $x\mapsto x\otimes 1$ and $x\mapsto 1\otimes x$ (\cite{SF3} 3.22.24). Then, the sequence
$$0\ra \op{Hom}_{\calA}(\calE_1,\calE_2) \ra \op{Hom}_{\calA'}(\calE_1',\calE_2') \ra \op{Hom}_{\calA''}(\calE_1'',\calE_2'')$$
is exact, where the first arrow is $u\mapsto u\circledast \op{Id}_{\calA'}$ and the second arrow is $u'\mapsto u'\circledast \pi_2 -u'\circledast \pi_1.$
\end{theorem}

\begin{theorem}[\ref{inddescentmod}, Descent of indexed modules]\label{thmintro2}
Let $f:\calA\ra \calA'$ be a faithfully flat morphism of $\upmu$-indexed algebras. Denote by
$$\pi_1,\pi_2:\calA' \ra \calA'\circledast_{\calA}\calA',\ p_{ij}:\calA'\circledast_{\calA}\calA' \ra \calA'\circledast_{\calA}\calA'\circledast_{\calA}\calA',\ 1\le i<j\le 3$$
the projection morphisms (\cite{SF3} 3.22.24). Let $\calE'$ be a $\upmu$-indexed $\calA'$-module. We suppose given a $\calA'\circledast_{\calA}\calA'$-linear isomorphism
$$\epsilon:\pi_2^*\calE' \isom \pi_1^*\calE'$$
satisfying the cocycle condition
$p_{13}^*\epsilon=p_{12}^*\epsilon \circ p_{23}^*\epsilon.$
Then, there exists a $\upmu$-indexed $\calA$-module $\calE,$ unique up to a unique isomorphism, and an isomorphism
$$\varphi:f^*\calE \isom \calE',$$
making the following diagram commutative
$$
\begin{tikzcd}
\pi_2^*f^*\calE \ar{r}{\pi_2^*\varphi} \ar{d} & \pi_2^*\calE' \ar{d}{\epsilon} \\
\pi_1^*f^*\calE \ar{r}{\pi_1^*\varphi} & \pi_1^*\calE',
\end{tikzcd}
$$
where the left vertical arrow is the canonical isomorphism.
\end{theorem}
\end{parag}

\begin{parag}
In section 6, we construct the local Cartier transform. We fix a positive integer $n$ and a log smooth morphism $\frakX \ra \frakS$ of fs logarithmic $p$-adic formal schemes, whose underlying formal schemes are flat and locally of finite type over $W.$ We denote by $X \ra S$ its special fiber, $F_{X/S}:X\ra X'$ the exact relative Frobenius and we suppose that $F_{X/S}$ lifts to an $\frakS$-morphism $F:\frakX \ra \frakX'$ of fs logarithmic $p$-adic formal schemes, such that $\frakX'$ is log smooth over $\frakS$ and flat and locally of finite type over $W.$ We denote by $\frakX_n \ra \frakS_n$ and $F_n:\frakX_n \ra \frakX_n'$ the reductions modulo $p^n.$ We identify the étale sites on $X,$ $X',$ $\frakX_n$ and $\frakX_n'$ via the nilpotent immersions $X\ra \frakX_n,$ $X'\ra \frakX_n'$ and $F_{X/S},$ which are all universal homeomorphisms. We set $\upmu=\ov{\calM}_X^{\mathrm{gp}}.$
We denote by $p\text{-}\boldsymbol{\op{MIC}}(\calB_{F_n})$ the category of $\upmu$-indexed $\calB_{F_n}$-modules $\calE'$ equipped with an integrable $p$-connection $\theta':\calE' \ra \calE' \otimes_{\Ox_{\frakX_n',\upmu}} \omega^1_{\frakX_n'/\frakS_n,\upmu},$ compatible with $\theta_{F_n}$ i.e. such that, for all local sections $a$ and $x$ of $\calB_{F_n}$ and $\calE$ respectively,
$$
\theta'(ax)=a\theta'(x)+x\otimes \theta_{F_n}(a).
$$
Similarly, we denote by $\boldsymbol{\op{MIC}}(\calA_{\frakX_n})$ the category of $\upmu$-indexed $\calA_{\frakX_n}$-modules equipped with an integrable connection compatible with $d_{\calA_{\frakX_n}}.$ We start by constructing a functor
\begin{equation}\label{eqintro2}
\Phi: p\text{-}\boldsymbol{\op{MIC}}(\calB_{F_n}) \ra \boldsymbol{\op{MIC}}(\calA_{\frakX_n}).
\end{equation}
Since $dF_{X/S}:F_{X/S}^*\omega^1_{X'/S} \ra \omega^1_{X/S}$ vanishes, the image of $dF:F^*\omega^1_{\frakX'/\frakS} \ra \omega^1_{\frakX/\frakS}$ lies in $p\omega^1_{\frakX/\frakS}.$ Then, there exists a unique morphism $p^{-1}dF:F_n^*\omega^1_{\frakX'_n/\frakS_n} \ra \omega^1_{\frakX_n/\frakS_n}$ fitting into the commutative diagram:
$$
\begin{tikzcd}
F^*\omega^1_{\frakX'/\frakS} \ar{r}{dF} \ar{dr} \ar{d} & p\omega^1_{\frakX/\frakS} \\
F_{n}^*\omega^1_{\frakX_n'/\frakS_n} \ar[swap]{dr}{p^{-1}dF} & \omega^1_{\frakX/\frakS} \ar[swap]{u}{\times p} \ar{d} \\
 & \omega^1_{\frakX_n/\frakS_n}.
\end{tikzcd}
$$
Let $(\calE',\theta')$ be an object of $p\text{-}\boldsymbol{\op{MIC}}(\calB_{F_n}).$ Set $\calE=\calA_{\frakX_n} \circledast_{\calB_{F_n}} \calE'$ and consider the composition
$$
\zeta: \calE' \xrightarrow{\theta'} \calE' \otimes_{\Ox_{\frakX_n',\upmu}} \omega^1_{\frakX_n'/\frakS_n,\upmu} \ra \calE' \otimes_{\Ox_{\frakX_n',\upmu}} F_n^*\omega^1_{\frakX_n'/\frakS_n,\upmu} \xrightarrow{\op{Id}_{\calE'} \otimes p^{-1}dF} \calE' \otimes_{\Ox_{\frakX_n',\upmu}} \omega^1_{\frakX_n/\frakS_n,\upmu},
$$
where the second arrow is induced by $\omega^1_{\frakX_n'/\frakS_n} \ra F_n^*\omega^1_{\frakX_n'/\frakS_n}=\Ox_{\frakX_n}\otimes_{\Ox_{\frakX_n'}}\omega^1_{\frakX_n'/\frakS_n},\ \omega\mapsto 1\otimes\omega.$
Let $\sigma$ be the addition map of $\upmu$ and $p_1,p_2:\upmu^2\ra \upmu$ the two projections. Set
$$
\nabla':\begin{array}[t]{clc}
\calA_{\frakX_n} \boxtimes \calE' & \ra & \left (\calA_{\frakX_n} \boxtimes \calE' \right ) \otimes_{\Ox_{\frakX_n,\upmu^2}} \omega^1_{\frakX_n/\frakS_n,\upmu^2} \\
a\otimes x' & \mapsto & a\otimes p_2^*\zeta(x')+1 \otimes x' \otimes p_1^*d_{\calA_{\frakX_n}}(a).
\end{array}
$$
Applying $\sigma_!$ to the integrable connection $\nabla'$ and passing to the coequalizer of \eqref{eqintro1}, we get an integrable connection $\nabla$ on $\calE$ and we define $\Phi(\calE',\theta')=(\calE,\nabla).$

For an integrable connection (resp. integrable $p$-connection) $\nabla$ on an $\Ox_{\frakX_n,\upmu}$-module $\calE,$ we say $\nabla$ is quasi-nilpotent if, for all local sections $s\in \Gamma \left (U,\upmu \right )$ over an étale $\frakX_n$-scheme $U,$ the connection $\nabla_s,$ pullback of $\nabla$ along the morphism of ringed topoi $s:U_{\text{ét}} \ra \frakX_{n/\upmu},$ is quasi-nilpotent. Our goal is to prove that $\Phi$ induces an equivalence between subcategories of quasi-nilpotent objects. To do so, we provide a stratified interpretation. Consider the dilatation $R_{\frakX'/\frakS}$ (\ref{intro1}) and denote by $P_{\frakX'/\frakS,1}$ the PD envelope of the immersion $\frakX' \ra R_{\frakX'/\frakS}.$ We also denote by $P_{\frakX/\frakS,0}$ the log PD envelope of the diagonal immersion $\frakX \ra \frakX\times_{\frakS}\frakX.$ Let $P_n$ and $P_n'$ be the reductions of $P_{\frakX/\frakS,0}$ and $P_{\frakX'/\frakS,1}$ respectively modulo $p^n$ and by $\calP_n$ and $\calP_n'$ their structural rings.
\begin{definition}
The integrable quasi-nilpotent connection $d_{\calA_{\frakX_n}}$ (resp. $p$-connection $\theta_{F_n}$) corresponds to a $\calP_{n,\upmu}$-stratification $\epsilon_{\calA_{\frakX_n}}$ (resp. $\calP'_{n,\upmu}$-stratification $\epsilon'_{\calB_{F_n}}$). We denote by
$$
\boldsymbol{\op{MHS}}(\calA_{\frakX_n}) \left (\text{resp.}\ \boldsymbol{\op{MHS}}(\calB_{F_n}) \right )
$$
the category of $\upmu$-indexed $\calA_{\frakX_n}$-modules (resp. $\upmu$-indexed $\calB_{F_n}$-modules) equipped with an $\epsilon_{\calA_{\frakX_n}}$-linear (resp. $\epsilon'_{\calB_{F_n}}$-linear) $\calP_{n,\upmu}$-stratification (resp. $\calP'_{n,\upmu}$-stratification).
\end{definition}
\begin{proposition}[\ref{indpropstrat2}]\label{propintro1}
We have equivalences of categories
$$\boldsymbol{\op{MHS}}(\calA_{\frakX_n}) \isom \boldsymbol{\op{MIC}}^{\mathrm{qn}}(\calA_{\frakX_n}),\quad \boldsymbol{\op{MHS}}(\calB_{F_n}) \isom p\boldsymbol{\op{\text{-}MIC}}^{\mathrm{qn}}(\calB_{F_n}),$$
where $\op{MIC^{qn}}(\calA_{\frakX_n})$ and $p\op{\text{-}MIC^{qn}}(\calB_{F_n})$ are the subcategories of $\boldsymbol{\op{MIC}}(\calA_{\frakX_n})$ and $p\op{\text{-}MIC}(\calB_{F_n})$ consisting of quasi-nilpotent connections and $p$-connections (\ref{def622}) respectively.
\end{proposition}
We prove, in \ref{prop432}, that the lifting $F$ of the Frobenius induces a morphism $\varphi:P_n \ra P_n'$ fitting into the commutative diagram
\begin{equation}\label{eqintro3}
\begin{tikzcd}
P_n \ar{rr}{\varphi} \ar{d} & & P_n' \ar{d} \\
\frakX_n\times_{\frakS_n} \frakX_n \ar{rr}{F_n\times_{\frakS_n}F_n} & & \frakX_n'\times_{\frakS_n}\frakX_n'.
\end{tikzcd}
\end{equation}
We consider $\calP_n$ as a $\calP'_n$-algebra via $\varphi.$ We then define the functor
\begin{equation}\label{Psiintro}
\Psi:\begin{array}[t]{clc}
\boldsymbol{\op{MHS}}(\calB_{F_n}) & \ra & \boldsymbol{\op{MHS}}(\calA_{\frakX_n}) \\
(\calE',\epsilon') & \mapsto & \left ( \calA_{\frakX_n} \circledast_{\calB_{F_n}} \calE',\epsilon_{\calA_{\frakX_n}} \circledast_{\calP'_{n,\upmu} \otimes_{\Ox_{\frakX_n',\upmu}} \calB_{F_n} } \epsilon' \right ).
\end{array}
\end{equation}
\begin{proposition}[\ref{prop443}]
The diagram
\begin{equation*}
\begin{tikzcd}
\boldsymbol{\op{MHS}(\calB_{F_n})} \ar{r}{\Psi} \ar[swap,sloped]{d}{\sim} & \boldsymbol{\op{MHS}(\calA_{\frakX_n})}\ar[sloped]{d}{\sim} \\
p\text{-}\boldsymbol{\op{MIC}}^{\mathrm{qn}}(\calB_{F_n}) \ar[hook]{d} & \boldsymbol{\op{MIC}}^{\mathrm{qn}}(\calA_{\frakX_n})\ar[hook]{d} \\
\boldsymbol{p\op{\text{-}MIC}(\calB_{F_n})} \ar{r}{\Phi} & \boldsymbol{\op{MIC}(\calA_{\frakX_n})}
\end{tikzcd}
\end{equation*}
where $\Phi$ is defined in \eqref{eqintro2}, the upper vertical arrows are the equivalences of categories in \ref{propintro1} and the lower vertical arrows are the canonical embeddings, is commutative.
\end{proposition}
The previous diagram implies that $\Phi$ preserves quasi-nilpotence. Then, to prove that it induces an equivalence between subcategories of quasi-nilpotent connections and $p$-connections, it is sufficient to prove that $\Psi$ is an equivalence of categories. For that, we complete the diagram \eqref{eqintro3} as follows: the ideal of the diagonal immersion $\frakX \ra \frakX \times_{\frakX'} \frakX$ has a unique PD structure. It follows that there exists a morphism $ \frakX\times_{\frakX'}\frakX \ra P$ whose reduction $\psi$ modulo $p^n$ fits into the commutative diagram
\begin{equation*}
\begin{tikzcd}
 & P_n \ar{rr}{\varphi} \ar{d} & & P_n'\ar{d} \\
\frakX_n\times_{\frakX_n'}\frakX_n \ar{ur}{\psi}  \ar{r} & \frakX_n\times_{\frakS_n}\frakX_n \ar{rr}{F_n\times_{\frakS_n}F_n} & & \frakX_n'\times_{\frakS_n}\frakX_n'.
\end{tikzcd}
\end{equation*}
The key to proving that $\Psi$ is an equivalence of categories is to note that, for an object $(\calE,\epsilon)$ of $\boldsymbol{\op{MHS}}(\calA_{\frakX_n}),$ the pullback $\psi^*\epsilon$ is an isomorphism
$$
\psi^*\epsilon:\calA_{\frakX_n\times_{\frakX_n'}\frakX_n} \circledast_{\calA_{\frakX_n}} \calE \isom \calE \circledast_{\calA_{\frakX_n}} \calA_{\frakX_n\times_{\frakX_n'}\frakX_n}.
$$
Using \ref{eqintro4}, we deduce descent datum of $\calE$ with respect to $\calB_{F_n} \ra \calA_{\frakX_n}.$ Then, using \ref{propintro2}, \ref{thmintro1} and \ref{thmintro2}, we deduce the desired result:

\begin{theorem}[\ref{Psiequiv} and \ref{Phiequiv}]
The functor $\Psi$ \eqref{Psiintro} is an equivalence of categories and the functor $\Phi$ \eqref{eqintro2} induces an equivalence of categories
$$
\Phi:p\boldsymbol{\op{\text{-}MIC}}^{\mathrm{qn}}(\calB_{F_n}) \isom  \boldsymbol{\op{MIC}}^{\mathrm{qn}}(\calA_{\frakX_n}).
$$
\end{theorem}

\end{parag}

\endgroup

\textbf{Acknowledgements.}
This paper was written during my stay at the Morningside Center of Mathematics, Beijing. I thank Daxin Xu for his invitation and the fruitful discussions we had. I also thank Atsushi Shiho for his insightful ideas and remarks. I am likewise grateful to my PhD advisor, Ahmed Abbes, for their guidance throughout my doctoral studies, which laid the foundation for this work.

\section{Notations and conventions}

\begingroup
\renewcommand{\thesubsection}{\thesection}

\begin{parag}\label{not1}
Throughout this article, we fix a prime number $p$ and a ring of Witt vectors $W$ of a perfect field of characteristic $p$. All formal schemes are assumed to be $p$-adic (\cite{Ahmed2010} 2.1.24) and all logarithmic structures are assumed to be fs i.e. fine and saturated. In particular, fiber products of logarithmic schemes (resp. logarithmic formal schemes) are taken in the category of fs logarithmic schemes (resp. fs logarithmic $p$-adic formal schemes).
\end{parag}

\begin{parag}
All considered monoids are assumed to be commutative.
\end{parag}

\begin{parag}\label{not3}
For a monoid $M$ and a ring $R,$ we denote by $R[M]$ the free $R$-module of basis $t^x,\ x\in M,$ equipped with the multiplication $t^xt^y=t^{x+y}$ for all $x,y\in M.$
\end{parag}

\begin{parag}
For an fs monoid $M,$ we denote by $A[M]$ the scheme $\op{Spec}\Z[M]$ equipped with the logarithmic structure induced by the canonical morphism $M\ra \Z[M].$
\end{parag}

\begin{parag}
If a gothic letter $\frakX$ denotes a logarithmic formal scheme and $n$ is a positive integer, we denote by $\frakX_n$ the logarithmic scheme obtained from $\frakX$ by reduction modulo $p^n$ (\cite{SF1} 7.8). The corresponding roman letter $X$ denotes $\frakX_1.$
\end{parag}

\begin{parag}
For a monoid $M,$ we denote by $\ov{M}$ the quotient of $M$ by the equivalence relation $\sim$ given, for all $x,y\in M,$ by $x\sim y$ if there exists an invertible $u\in M$ such that $x=y+u.$ If $\calM$ is a sheaf of monoids, we denote by $\ov{\calM}$ the sheaf associated to the presheaf $U\mapsto \ov{\calM(U)}.$
\end{parag}

\begin{parag}
For integers $1\le i\le d,$ we denote by $\epsilon_{d,i}$ the element of $\N^d$ whose all coefficients are zero except for the $i$th coefficient which is equal to $1.$ If no confusion arises, we denote $\epsilon_{d,i}$ simply by $\epsilon_i.$
\end{parag}

\begin{parag}
For elements $x_1,\hdots,x_d$ of a PD algebra and $I=(I_1,\hdots,I_d)\in \N^d,$ we set $x^{[I]}=\prod_{i=1}^dx_i^{[I_i]}.$
\end{parag}

\begin{parag}
For a ringed topos $(X,\Ox_X)$ and bimodules $\calM$ and $\calN$ over $\Ox_X,$ we denote by $\calM \otimes_{\Ox_X} \calN$ the tensor product of the right $\Ox_X$-module $\calM$ and the left $\Ox_X$-module $\calN.$
\end{parag}

\begin{parag}
For a morphism of logarithmic schemes $X\ra S,$ we denote by $\omega^1_{X/S}$ the $\Ox_X$-module of logarithmic differentials (\cite{Kat89} 1.7).
\end{parag}

\endgroup

\section{Generalities on Frobenii and dilatations}

\subsection{Frobenii and charts}

\begin{parag}\label{parag111}
Recall that all logarithmic schemes are assumed to be fs and their fiber products are taken in the category of fs logarithmic schemes (\ref{not1}). For a logarithmic scheme $T$ of characteristic $p,$ we denote by $F_T$ the absolute Frobenius morphism of $T,$ i.e. $F_T$ is the identity on the underlying topological space, elevation to the $p$th power on $\Ox_T$ and multiplication by $p$ on $\calM_T.$ For a monoid $M,$ let $F_M:M\ra M$ be the morphism defined by $F_M:x\mapsto px.$
Let $f:X\ra S$ be a morphism of logarithmic schemes in characteristic $p.$ We have the following commutative diagram
\begin{equation}\label{diagfrob}
\begin{tikzcd}
X \ar{dr}{F_{X/S}} \ar[bend right=-30]{drr}{F_X} \ar[swap,bend right=30]{ddr}{f} & & \\
 & X' \ar{r}{\pi} \ar{d} & X \ar{d}{f} \\
 & S \ar{r}{F_S} & S,
\end{tikzcd}
\end{equation}
where the square is cartesian. The morphism $F_{X/S}$ is called the exact relative Frobenius (\cite{Ogus2018} IV 3.3.9).

Suppose that $f$ has a chart $\theta:P\ra Q,$ where $P$ and $Q$ are fs monoids. Let $Q'$ be the inverse image of $Q$ by the morphism of groups
$
Q^{\mathrm{gp}}\oplus_{P^{\mathrm{gp}},F_P^{\mathrm{gp}}}P^{\mathrm{gp}} \ra Q^{\mathrm{gp}},\quad (x,t) \mapsto px+\theta^{\mathrm{gp}}(t).
$
Note that $Q'=\left (Q\oplus_{P,F_P}P\right )^{\mathrm{sat}}.$ Indeed, if $(x,t) \in \left (Q\oplus_{P,F_P}P\right )^{\mathrm{sat}} \subset Q^{\mathrm{gp}} \oplus_{P^{\mathrm{gp}},F_P^{\mathrm{gp}}} P^{\mathrm{gp}}$ then there exists a positive integer $n$ such that $n(x,t)\in Q\oplus_{P,F_P}P.$ Then $n(px+\theta^{\mathrm{gp}}(t))\in Q.$ Since $Q$ is saturated, we have $px+\theta^{\mathrm{gp}}(t)\in Q$ and $(x,t)\in Q'.$ Conversely, if $(x,t)\in Q'$ then $px+\theta^{\mathrm{gp}}(t)\in Q.$ It follows, that in $Q^{\mathrm{gp}}\oplus_{P^{\mathrm{gp}},F_P^{\mathrm{gp}}}P^{\mathrm{gp}},$ we have $p(x,t)=(px+\theta^{\mathrm{gp}}(t),0)\in Q\oplus_{P,F_P}P.$ Then, $(x,t)\in \left (Q\oplus_{P,F_P}P\right )^{\mathrm{sat}}.$ 
The relative Frobenius $F_{X/S}$ has a chart
\begin{equation}\label{FQP}
F_{Q/P}:Q' \ra Q,\quad (x,t) \mapsto px+\theta^{\mathrm{gp}}(t).
\end{equation}
\end{parag}

\begin{proposition}\label{prop212}
The morphism $F_{Q/P}$ \eqref{FQP} is exact. In addition, if $\theta^{\mathrm{gp}}$ is injective and the torsion subgroup of its cokernel is finite of order coprime with $p,$ then $F_{Q/P}$ is injective. 
\end{proposition}

\begin{proof}
The exactness follows from the definition of $Q'$ since, for all $(x,t)\in Q'^{\mathrm{gp}} = Q^{\mathrm{gp}}\oplus_{P^{\mathrm{gp}},F_P^{\mathrm{gp}}}P^{\mathrm{gp}},$ $F_{Q/P}^{\mathrm{gp}}(x,t)\in Q$ if and only if $(x,t)\in Q'.$ The injectivity is proved in (\cite{SF1} 3.2).
\end{proof}

\begin{lemma}\label{lem112}
Let $C$ be the cokernel of $F_{Q/P}^{\mathrm{gp}}:Q'^{\mathrm{gp}}=Q^{\mathrm{gp}}\oplus_{P^{\mathrm{gp}},F_P^{\mathrm{gp}}}P^{\mathrm{gp}} \ra Q^{\mathrm{gp}},\quad (x,t) \mapsto px+\theta^{\mathrm{gp}}(t).$ Then, the morphism $\left ( Q\oplus_{Q'} Q\right )^{\mathrm{sat}} \ra Q\oplus C,$ defined by
$$
p_1:\begin{array}[t]{clc}
Q & \ra & Q\oplus C,\\
x & \mapsto & (x,0)
\end{array}
\quad
p_2:\begin{array}[t]{clc}
Q & \ra & Q\oplus C,\\
y & \mapsto & (y,\ov{y})
\end{array}
$$
is an isomorphism.
\end{lemma}

\begin{proof}
First, note that the diagram
$$
\begin{tikzcd}
Q' \ar{r}{F_{Q/P}} \ar[swap]{d}{F_{Q/P}} & Q \ar{d} \ar[bend right=-30]{ddr}{p_1} & \\
Q \ar{r} \ar[bend right=30]{rrd}{p_2} & \left ( Q\oplus_{Q'} Q\right )^{\mathrm{sat}} \ar[dashed]{dr}{\varphi} & \\
 & & Q\oplus C,
\end{tikzcd}
$$
is commutative, hence the existence of $\varphi.$
In addition, $\varphi^{\mathrm{gp}}$ is clearly an isomorphism. It remains to prove that $\varphi^{\mathrm{gp}}\left (\left ( Q\oplus_{Q'} Q\right )^{\mathrm{sat}} \right )=Q\oplus C.$ It is clear that $\varphi^{\mathrm{gp}} \left (\left ( Q\oplus_{Q'} Q\right )^{\mathrm{sat}} \right ) \subset Q\oplus C.$ Conversely, if $(x,y) \in Q^{\mathrm{gp}} \oplus_{Q'^{\mathrm{gp}}}Q^{\mathrm{gp}}$ such that $\varphi^{\mathrm{gp}}(x,y) =(x+y,\ov{y}) \in Q\oplus C,$ then $x+y\in Q.$ But, in $Q^{\mathrm{gp}}\oplus_{Q'^{\mathrm{gp}}}Q^{\mathrm{gp}},$ we have $p(x,y)=(p(x+y),0) \in Q\oplus_{Q'}Q$ so $(x,y)\in \left ( Q\oplus_{Q'} Q\right )^{\mathrm{sat}}.$ This finishes the proof.
\end{proof}

\begin{remark}\label{rem113}
Let $P\ra Q$ be a chart of a morphism of logarithmic schemes $X\ra S$ of characteristic $p,$ such that $P$ and $Q$ are fs, and let $C$ be the cokernel of $F_{Q/P}^{\mathrm{gp}}$ \eqref{FQP}. Then, by \ref{lem112}, $Q\oplus C$ is a chart of $X\times_{X'}X$ and the projections $X\times_{X'}X \ra X$ have charts
\begin{equation}
\begin{array}[t]{clc}
Q & \ra & Q\oplus C,\\
x & \mapsto & (x,0)
\end{array}
\qquad
\begin{array}[t]{clc}
Q & \ra & Q\oplus C\\
y & \mapsto & (y,\ov{y}).
\end{array}
\end{equation}
\end{remark}

\begin{proposition}\label{prop114}
Let $f:X\ra S$ be a morphism of logarithmic schemes in characteristic $p$ and $F_{X/S}:X\ra X'$ the exact relative Frobenius. The diagonal immersion $X\ra X\times_{X'}X$ is strict, closed and a universal homeomorphism.
\end{proposition}

\begin{proof}
The absolute Frobenius is affine and since $F_X=\pi \circ F_{X/S}$ and $\pi$ is the base change of $F_S$ \eqref{diagfrob}, $F_{X/S}$ is also affine. It follows that the diagonal immersion $X\ra Z$ is closed, where $Z$ is the fiber product of $F_{X/S}$ with itself in the category of, not necessarily fs, logarithmic schemes. Since the canonical morphism $X\times_{X'}X\ra Z$ is affine, the immersion $X\ra X\times_{X'}X$ is closed.
Since $F_{X/S}$ is inseparable (\cite{Ogus2018} III 2.4.2), so is the projection $X\times_{X'}X \ra X$ by (\cite{Ogus2018} III 2.4.8). It follows, by (\cite{Ogus2018} III 2.4.8), that the diagonal immersion $X\ra X\times_{X'}X$ is inseparable and in particular a universal homeomorphism and exact. We conclude by the fact that exact immersions are strict.
\end{proof}

\subsection{\texorpdfstring{$p$}{p} %
-connections and dilatations}

\begin{parag}\label{parKhaminei}
Let $f:X\ra S$ be a morphism of logarithmic schemes, $\calE$ an $\Ox_X$-module and $\lambda\in \Gamma \left (X,\Ox_X \right ).$ A \emph{$\lambda$-connection on $\calE$} is a morphism of abelian groups $\nabla:\calE\ra \calE\otimes_{\Ox_X}\omega^1_{X/S}$ satisfying the modified Leibniz rule
$$\nabla(ax)=a\nabla(x)+\lambda x\otimes da,$$
for all local sections $a$ and $x$ of $\Ox_X$ and $\calE$ respectively.
A $\lambda$-connection $\nabla$ on $\calE$ induces, for all integers $i \ge 1,$ a morphism
$$\nabla^i:\calE\otimes_{\Ox_X}\omega^i_{X/S}\ra \calE\otimes_{\Ox_X}\omega^{i+1}_{X/S}$$
defined, for all local sections $x$ and $\omega$ of $\calE$ and $\omega^i_{X/S}$ respectively, by
$$\nabla^i(x\otimes \omega)=\nabla(x)\wedge \omega+\lambda x\otimes d\omega.$$
We say that the $\lambda$-connection $\nabla$ is \emph{integrable} if $\nabla^1\circ \nabla=0.$ We denote by $\boldsymbol{\lambda\text{-}\op{MIC}(X/S)}$ the category of $\Ox_X$-modules equipped with an integrable $\lambda$-connection.
$1$-connections are simply called \emph{connections}.

If $\partial:\omega^1_{X/S} \ra \Ox_X$ is an $\Ox_X$-linear map and $\nabla:\calE \ra \calE \otimes_{\Ox_X}\omega^1_{X/S}$ is a $\lambda$-connection, we denote by $\nabla \left (\partial \right )$ the composition
\begin{equation}\label{nablacomp}
\calE\xrightarrow{\nabla}\calE\otimes_{\Ox_X}\omega^1_{X/S}\xrightarrow{\op{Id}_{\calE}\otimes \partial}\calE.
\end{equation}
Suppose we have local coordinates $m_1,\hdots,m_d \in \Gamma \left (X,\calM_X\right )$ of $X$ with respect to $S,$ i.e. $(\op{dlog}m_1,\hdots,\op{dlog}m_d)$  is a basis of the $\Ox_X$-module $\omega^1_{X/S}.$ Denote by $\left (\partial_{\epsilon_1},\hdots,\partial_{\epsilon_d} \right )$ the dual basis. By (\cite{SF1} 5.2), if $\nabla$ is integrable then $\nabla(\partial_{\epsilon_i})$ and $\nabla(\partial_{\epsilon_j})$ commute for all $1\le i,j\le d.$
\end{parag}

\begin{definition}\label{def227}
Let $X\ra S$ be a log smooth morphism of logarithmic schemes and $n\ge 0$ an integer. For a local system of coordinates $m_1,\hdots,m_d$ of $X/S$ and an integrable $p^n$-connection $\nabla$ on an $\Ox_{X}$-module $\calE,$ let $\left (\partial_{\epsilon_1},\hdots,\partial_{\epsilon_d} \right )$ be the dual basis of $\left (\op{dlog}m_1,\hdots,\op{dlog}m_d \right ).$ We define the following endomorphisms:
$$
\nabla \left ( \partial_{\epsilon_i} \right ) : \calE \xrightarrow{\nabla} \calE \otimes_{\Ox_{X}} \omega^1_{X/S} \xrightarrow{\op{Id}_{\calE} \otimes \partial_{\epsilon_i}} \calE\quad \forall 1\le i\le d,\qquad 
\nabla \left ( \partial_I \right ) = \prod_{i=1}^d \prod_{j=0}^{I_i-1} \left ( \nabla \left (\partial_{\epsilon_i} \right )-p^nj \right )\quad \forall I\in \N^d.
$$
The $p^n$-connection $\nabla$ is said to be \emph{quasi-nilpotent} if, for a local system of coordinates $m_1,\hdots,m_d$ of $X/S,$ for all local sections $x$ of $\calE,$ there exists locally an integer $N$ such that $\nabla\left (\partial_I \right )(x)=0$ for all $I\in \N^d$ such that $|I|\ge N.$ 
\end{definition}

\begin{remark}\label{remquasinil}
\begin{enumerate}
\item The definition of quasi-nilpotence is independent of the choice of local coordinates.
\item If $n$ is a positive integer and $X\ra S$ is a log smooth morphism of logarithmic schemes over $W/p^nW$ and $\nabla$ is an integrable connection on an $\Ox_X$-module $\calE,$ then the $p$-connection $\theta=p\nabla$ is quasi-nilpotent. Indeed, if we have local coordinates, then $\theta \left ( \partial_I \right )\left (\calE \right )\subset p^n\calE=0$ for all $I\in \N^d$ such that $|I|\ge n.$
\end{enumerate}
\end{remark}

\begin{parag}\label{parag324f}
Let $n$ be a positive integer, $X\ra S$ be a morphism of logarithmic schemes over $W/p^nW,$ $P_{X/S}$ the log PD envelope of the diagonal immersion $X\ra X\times_SX,$ $\calP_{X/S}$ its structural ring, $\ov{\calI}$ its PD ideal and $p_1,p_2:P_{X/S} \ra X$ the two projections. The canonical closed immersion $\Delta:X\ra P_{X/S}$ is a nilimmersion, we use it to identify the étale sites of $X$ and $P_{X/S}.$ Consider the exact sequence of monoids (\cite{SF1} (4.3.1))
$$
0 \ra 1+\ov{\calI} \xrightarrow{\lambda} \calM_{P_{X/S}} \xrightarrow{\Delta^{\flat}} \calM_X \ra 0,
$$
where $\lambda=\alpha_{P_{X/S}}^{-1}.$ Suppose that we have local coordinates $m_1,\hdots,m_d$ of $X/S.$ For $1\le i\le d,$ let $\eta_i$ be the unique local section of $\ov{\calI}$ such that
$$
\lambda(1+\eta_i)+p_1^{\flat}m_i=p_2^{\flat}m_i.
$$
Then, by (\cite{Kat89} 6.5), $\calP_{X/S}=\Ox_X \langle \eta_1 ,\hdots,\eta_d \rangle.$ The following result links $\calP_{X/S}$-stratifications to quasi-nilpotent integrable connections:
\end{parag}

\begin{theorem}[\cite{Ogus94} 1.1.8]\label{thmKato}
Let $X\ra S$ be a log smooth morphism of logarithmic schemes and an $\Ox_X$-module $\calE.$ The following two data are equivalent:\begin{enumerate}
\item An integrable quasi-nilpotent connection $\nabla$ on $\calE.$
\item A $\calP_{X/S}$-stratification on $\calE.$
\end{enumerate}
Suppose we have local coordinates $m_1,\hdots,m_d$ of $X/S.$ Let $\eta_1,\hdots,\eta_d$ be as in \ref{parag324f}. If $\nabla$ is an integrable quasi-nilpotent connection on $\calE,$ then the corresponding stratification $\epsilon$ is given by
$$
\epsilon(1\otimes x)=\sum_{I\in \N^d}\nabla \left (\partial_I \right )(x)\otimes \eta^{[I]},
$$
where $x$ is a local section of $\calE$ and $\nabla \left (\partial_I \right )$ is defined in \ref{def227}.
\end{theorem}

\begin{lemma}\label{lem121}
Let $\theta : P\ra Q$ be a morphism of monoids such that $\ov{\theta}:\ov{P} \ra \ov{Q}$ is surjective and let $M$ be the inverse image of $Q$ by $\theta^{\mathrm{gp}}.$ Then, the morphism $M\ra Q$ is strict.
\end{lemma}

\begin{proof}
It is clear that $\ov{M}\ra \ov{Q}$ is surjective. Let $x,x'\in M$ such that $\theta^{\mathrm{gp}}(x)=\theta^{\mathrm{gp}}(x')$ in $\ov{Q}.$ Then, there exists an invertible $u\in Q$ such that $\theta^{\mathrm{gp}}(x)=\theta^{\mathrm{gp}}(x')+u$ in $Q.$ It follows that $x-x'\in M$ and $x'-x\in M$ so $x-x'$ is invertible in $M.$ Since $x=x'+(x-x'),$ we get the injectivity of $\ov{M} \ra \ov{Q}.$
\end{proof}

\begin{lemma}\label{lem122}
Let $\Delta:X\ra Y$ be an immersion of logarithmic schemes or logarithmic locally noetherian formal schemes. Suppose that $\Delta$ admits a chart $\theta:P\ra Q$ such that $\ov{\theta}:\ov{P} \ra \ov{Q}$ is surjective. Then, $\Delta$ factors as a strict closed immersion $i$ followed by a separated and log étale morphism. In addition, $i$ has a chart given by a strict morphism of monoids $M\ra Q.$
\end{lemma}

\begin{proof}
We may suppose that $\Delta$ is closed. Let $M$ be the inverse image of $Q$ by $\theta^{\mathrm{gp}}.$ By \ref{lem121}, the morphism $X\ra Y\times_{A[P]} A[M],$ induced by $X\ra Y$ and $X\ra A[Q] \ra A[M],$ is strict. It is also a closed immersion since so is $X\ra Y.$ Since $Y\times_{A[P]} A[M] \ra Y$ is log étale and affine, the result is proved.
\end{proof}

\begin{lemma}
Let $X\ra S$ be a log smooth morphism of logarithmic schemes or logarithmic locally noetherian formal schemes and suppose that the diagonal immersion $X\ra X\times_SX$ decomposes into a strict closed immersion $\Delta:X\ra Y$ followed by a log étale morphism $g:Y\ra X\times_SX.$ Then, the conormal sheaf $\calN_{\Delta}$ of $\Delta$ is canonically isomorphic to $\omega^1_{X/S}.$
\end{lemma}

\begin{proof}
Let $p_1,p_2:X\times_SX\ra X$ be the projections. By (\cite{Ogus2018} IV 3.2.2), we have the exact sequence of $\Ox_X$-modules:
$$
0\ra \calN_{\Delta} \ra \Delta^*\omega^1_{Y/X,p_1} \ra \omega^1_{X/X} \ra 0.
$$
We get a canonical isomorphism $\calN_{\Delta} \isom \Delta^*\omega^1_{Y/X,p_1}.$ Since $g$ is log étale, we have
$$
\Delta^*\omega^1_{Y/X}=\Delta^*g^*\omega^1_{X\times_SX/X,p_1}=\Delta^*g^*p_2^*\omega^1_{X/S}=\omega^1_{X/S}.
$$
\end{proof}

\begin{proposition}\label{prop224x}
Let $\frakX \ra \frakS$ be a morphism of logarithmic formal schemes, flat and locally of finite type over $W.$ Suppose we have two factorizations $\frakX \xrightarrow{\Delta} \frakY \ra \frakX\times_{\frakS}\frakX$ and $\frakX \xrightarrow{\Delta'} \frakY' \ra \frakX\times_{\frakS}\frakX$ as in \ref{lem122}. Let $R$ (resp. $R'$) be the dilatation of the ideal of $X\ra \frakX \xrightarrow{\Delta} \frakY$ (resp. $X\ra \frakX \xrightarrow{\Delta'} \frakY'$) with respect to $p$ (\cite{Ahmed2010} 3.2.3.4), equipped with the logarithmic structure pullback of $\calM_{\frakY}$ (resp. $\calM_{\frakY'}$). Then, we have a canonical isomorphism of logarithmic formal schemes $R' \ra R.$
\end{proposition}

\begin{proof}
The closed immersions $\Delta$ and $\Delta'$ induce a closed immersion $\frakX \ra \frakY \times_{\frakX \times_{\frakS} \frakX}\frakY'$ which factors, locally by \ref{lem122}, as a strict closed immersion $\frakX \ra \frakZ$ followed by a separated and log étale $\frakZ \ra \frakY \times_{\frakX \times_{\frakS} \frakX}\frakY'.$ Since $\frakY \times_{\frakX\times_{\frakS}\frakX}\frakY' \ra \frakX\times_{\frakS}\frakX$ is also separated and log étale, we can reduce to the case where we have a morphism $g:\frakY' \ra \frakY$ over $\frakX\times_{\frakS} \frakX$ such that the following diagram is commutative:
$$
\begin{tikzcd}
 & \frakY' \ar{d}{g} \\
\frakX \ar{ur}{\Delta'} \ar{r}{\Delta} & \frakY.
\end{tikzcd}
$$
Since $\frakY$ and $\frakY'$ are both log étale on $\frakX\times_{\frakS} \frakX,$ the morphism $g$ is also log étale. Since $\Delta$ and $\Delta'$ are both strict, $g$ is also strict on a neighborhood of $\Delta'(\frakX).$ Replacing $\frakY$ and $\frakY'$ by neighborhoods of $\Delta(\frakX)$ and $\Delta'(\frakX)$ respectively, we may suppose that $g$ is strict and hence étale. The morphism $h:\frakX \ra \frakY'\times_{\frakY}\frakX,$ induced by $\Delta'$ and $\op{Id}_{\frakX},$ is then a section of the étale morphism $\frakY'\times_{\frakY}\frakX \ra \frakX.$ It is hence an open immersion and we have a closed subset $\frakZ=\left ( \frakY'\times_{\frakY}\frakX \right )\setminus h(\frakX).$ The projection $q:\frakY'\times_{\frakY}\frakX \ra \frakY'$ is a closed immersion, as base change of $\Delta.$ Consider the strict open formal subscheme $\frakU=\frakY'\setminus q(\frakZ).$ The open formal subscheme $\frakU$ is an open neighborhood of $\Delta'(\frakX)$ so the canonical morphism $R' \ra \frakY'$ factors through $\frakU.$ Both $\Delta'$ and $R' \ra \frakY'$ factor through $\frakU,$ we can hence replace $\frakY'$ with $\frakU.$ Note that doing so does not change $R'.$ A point $z$ of $\frakY' \times_{\frakY}\frakX$ belongs to the open subset $\frakU \times_{\frakY}\frakX$ if and only if $q(z)\in \frakU$ if and only if $z\notin \frakZ.$ It follows that $\frakU\times_{\frakY}\frakX=h(\frakX),$ we can hence reduce to the case where $h:\frakX \ra \frakY'\times_{\frakY}\frakX$ is an isomorphism. This yields the equality $\calI'=g^*\calI,$ where $\calI$ and $\calI'$ are the ideals of $\Delta$ and $\Delta'$ respectively. We have $\calI\Ox_R\subset p\Ox_R.$ We deduce that the image of $\calI'$ in $R\times_{\frakY}\frakY'$ is a subideal of $(p).$ In addition, $g$ is étale hence flat, thus so is $R\times_{\frakY}\frakY' \ra R.$ It follows that $R\times_{\frakY}\frakY'$ is flat over $W$ and there exists a unique morphism $f:R\times_{\frakY}\frakY' \ra R'$ fitting into the commutative diagram
$$
\begin{tikzcd}
 & R' \ar{d} \\
R\times_{\frakY}\frakY' \ar{r} \ar{ur}{f} & \frakY'.
\end{tikzcd}
$$
By the universal property of $R,$ the composition $R' \ra \frakY' \xrightarrow{g} \frakY$ factors uniquely through the canonical morphism $R \ra \frakY,$ yielding $R' \ra R.$ It induces, along with $R' \ra \frakY',$ a morphism $f':R' \ra R\times_{\frakY}\frakY'.$ The equalities $f\circ f'=\op{Id}_{R'}$ and $f'\circ f=\op{Id}_{R\times_{\frakY}\frakY'}$ are checked using the universal properties of $R$ and $R'.$ We hence have  an isomorphism
$$
f':R' \isom R\times_{\frakY}\frakY'.
$$
Since both $R'$ and $R$ are flat over $W,$ to prove that $R' \xrightarrow{f'} R\times_{\frakY}\frakY' \ra R$ is an isomorphism, it is sufficient to prove that it is so modulo $p.$ Modulo $p,$ the canonical morphism $R_1\ra Y$ factors through $\Delta_1:X\ra Y.$ So we have
$$
R'_1 \isom R_1\times_YY' \isom R_1\times_X(X\times_YY')\xrightarrow{\op{Id}_{R_1}\times h_1^{-1}}R_1.
$$
We have thus constructed locally an isomorphism $R' \isom R.$ By a classical gluing argument, these isomorphisms glue together into the desired isomorphism.
\end{proof}

\begin{parag}\label{prop124}
Let $\frakX \ra \frakS$ be a morphism of logarithmic formal schemes, flat and locally of finite type over $W,$ $\left (\frakS_i \ra \frakS \right )_{i\in I}$ and $\left (\frakX_{ij} \ra \frakX_i=\frakX \times_{\frakS} \frakS_i \right )_{j\in I_i}$ étale covers such that the induced morphism $\frakX_{ij} \ra \frakS_i$ has a chart for all $i$ and $j.$
For $j,k\in I_i,$ set
$$
V_{ijk}=\frakX_{ij}\times_{\frakS_i}\frakX_{ik}.
$$
The morphism
$
\frakX_{ij}\times_{\frakX_i}\frakX_{ik}\ra V_{ijk}
$
is the base change of the diagonal immersion
\begin{equation}\label{eq2251x}
\frakX_i \ra \frakX_i\times_{\frakS_i}\frakX_i
\end{equation}
along the étale morphism $V_{ijk}\ra \frakX_i\times_{\frakS_i}\frakX_i,$ hence is again an immersion and has a chart satisfying the condition of \ref{lem122}.
By \ref{lem122}, the immersion $\frakX_{ij}\times_{\frakX_i}\frakX_{ik}\ra V_{ijk}$ factors as
$
\frakX_{ij}\times_{\frakX_i}\frakX_{ik}\ra \frakY_{ijk}\ra V_{ijk},
$
where $\frakX_{ij}\times_{\frakX_i}\frakX_{ik}\ra \frakY_{ijk}$ is a strict closed immersion and $\frakY_{ijk}\ra V_{ijk}$ is separated and log étale. Let
$
R_{ijk}
$
be the dilatation of the ideal of $X_{ij}\times_{X_i}X_{ik} \ra \frakY_{ijk}$ with respect to $p$ (\cite{Ahmed2010} 3.2.3.4), equipped with the logarithmic structure pullback of $\calM_{\frakY_{ijk}}.$ By the following proposition, the logarithmic formal schemes $R_{ijk}$ glue into a logarithmic formal scheme $R_{\frakX/\frakS}$ over $\frakX\times_{\frakS}\frakX$ and we obtain a strict immersion $\frakX \ra R_{\frakX/\frakS}.$
\end{parag}

\begin{proposition}\label{prop3211}
Keep the hypothesis and notation of \ref{prop124}. Then, the logarithmic formal schemes $R_{ijk}$ glue into a logarithmic formal scheme $R_{\frakX/\frakS}$ over $\frakX\times_{\frakS}\frakX$ and we obtain a strict immersion $\frakX \ra R_{\frakX/\frakS}.$
\end{proposition}

\begin{proof}
For all $i\in I,$ the morphisms $\frakX_{ij}\ra \frakX_i$ for $j\in I_i$ form an étale cover of $\frakX_i.$
Then $$\left (V_{ijk} \ra \frakX_i \times_{\frakS_i} \frakX_i \right )_{j,k\in I_i}$$ is an étale cover.
For $i\in I$ and $j,k\in I_i,$ set
$
D_{ijk}=\frakX_{ij}\times_{\frakX_i}\frakX_{ik}.
$
Fix $i\in I.$
We glue the local objects $R_{ijk}$ over the étale cover $\left (V_{ijk}\ra \frakX_i \times_{\frakS_i}\frakX_i \right )_{j,k\in I_i}.$ On an overlap
$$
V_{ijk}\times_{\frakX_i\times_{\frakS_i}\frakX_i}V_{ij'k'},
$$
both $R_{ijk}$ and $R_{ij'k'}$ are dilatations corresponding to the base change of the same diagonal immersion \eqref{eq2251x}. By \ref{prop224x}, there exists a canonical isomorphism on the overlap
$
R_{ijk} \isom R_{ij'k'}.
$
These isomorphisms satisfy the cocycle condition on triple overlaps by the universal property of the dilatation. We hence have a logarithmic formal scheme
$$
R_{\frakX_i/\frakS_i}\ra \frakX_i\times_{\frakS_i}\frakX_i,
$$
well-defined up to a unique isomorphism, whose restriction to each $V_{ijk}$ is canonically isomorphic to $R_{ijk},$ along with a strict immersion $\frakX_i \ra R_{\frakX_i/\frakS_i}.$
It remains to glue the $R_{\frakX_i/\frakS_i}$ for varying $i.$ This is done similarly to before and so we obtain a logarithmic formal scheme
$$
R_{\frakX/\frakS}\ra \frakX\times_{\frakS} \frakX,
$$
well-defined up to a unique isomorphism, whose restriction to $\frakX_i \times_{\frakS_i}\frakX_i$ is canonically isomorphic to $R_{\frakX_i/\frakS_i},$ along with a strict immersion $\frakX \ra R_{\frakX/\frakS}.$
\end{proof}

\begin{proposition}[\cite{SF1} 11.15]\label{prop125*}
Let $f:\frakX \ra \frakS$ be a log smooth morphism of logarithmic formal schemes, flat and locally of finite type over $W,$ admitting a chart $\theta:P\ra Q$ such that $\theta^{\mathrm{gp}}$ is injective, the torsion subgroup of its cokernel is finite of order coprime with $p,$ and the morphism $\frakX \ra \frakS \times_{\op{Spf}\Z_p\langle P\rangle} \op{Spf}\Z_p\langle Q\rangle,$ induced by $f$ and the chart $\frakX \ra \op{Spf}\Z_p\langle Q\rangle,$ is étale and strict. Consider a local decomposition of the diagonal immersion into a strict closed immersion $\Delta:\frakX \ra \frakY$ followed by a log étale separated morphism, as in \ref{lem122}. Suppose we have local coordinates $m_1,\hdots,m_d\in \calM_{\frakX}.$ Let $\calI$ be the ideal of $\Delta,$ $p_1,p_2:\frakY \ra \frakX$ the projections and consider the exact sequence of monoids (\cite{SF1} (4.3.1))
$$
0\ra \Delta^{-1}(1+\calI) \xrightarrow{\lambda} \Delta^{-1}\calM_{\frakY} \xrightarrow{\Delta^{\flat}} \calM_{\frakX} \ra 0.
$$
For $1\le i\le d,$ let $\widetilde{\eta}_i$ be the unique local section of $\calI$ such that $\lambda(1+\widetilde{\eta}_i)+\Delta^{-1}p_1^{\flat}m_i = \Delta^{-1}p_2^{\flat}m_i.$ Then, we have an isomorphism over $\frakS:$
$$
\Ox_{\frakX} \{ x_1,\hdots,x_d \} \isom p_{1*}\Ox_{R_{\frakX/\frakS}},\quad x_i \mapsto \frac{\widetilde{\eta}_i}{p},
$$
where $\Ox_{\frakX} \{ x_1,\hdots,x_d \}$ denotes the $p$-adic completion of the polynomial algebra $\Ox_{\frakX} [ x_1,\hdots,x_d ].$ Under these conditions, we will write
\begin{equation}
p_{1*}\Ox_{R_{\frakX/\frakS}} = \Ox_{\frakX} \left \{ \frac{\widetilde{\eta}_1}{p},\hdots, \frac{\widetilde{\eta}_d}{p} \right \}.
\end{equation}
If $P_{\frakX/\frakS,1}$ is the PD envelope of the canonical immersion $\frakX \ra R_{\frakX/\frakS},$ then
\begin{equation}
\Ox_{P_{\frakX/\frakS,1}} = \Ox_{\frakX} \left \langle \left \langle \frac{\widetilde{\eta}_1}{p},\hdots, \frac{\widetilde{\eta}_d}{p} \right \rangle \right \rangle,
\end{equation}
where $\Ox_{\frakX} \left \langle \left \langle \frac{\widetilde{\eta}_1}{p},\hdots, \frac{\widetilde{\eta}_d}{p} \right \rangle \right \rangle$ is the $p$-adic completion of the PD algebra $\Ox_{\frakX} \left \langle \frac{\widetilde{\eta}_1}{p},\hdots, \frac{\widetilde{\eta}_d}{p} \right \rangle.$
\end{proposition}

\begin{parag}\label{prop125}
Let $n$ be a positive integer, $\frakX \ra \frakS$ a log smooth morphism of logarithmic formal schemes, flat and locally of finite type over $W.$ Let $\frakX \xrightarrow{\Delta} \frakY \ra \frakX \times_{\frakS} \frakX$ be a local decomposition of the diagonal immersion into a strict closed immersion followed by a log étale and separated morphism (\ref{lem122}). Consider the formal scheme $R_{\frakX/\frakS}$ defined in \ref{prop124} and let $P_{\frakX/\frakS,1}$ be the PD envelope of the canonical immersion $\frakX \ra R_{\frakX/\frakS},$ $P_n$ its reduction modulo $p^n$ and $\calP_n$ its structural ring. Locally, we have local coordinates $m_1,\hdots,m_d$ of $\frakX_n/\frakS_n,$ that lift to local sections $\widetilde{m}_1,\hdots,\widetilde{m}_d$ of $\calM_{\frakX}.$ Let $p_1,p_2:\frakY \ra \frakX$ be the projections and $\calI$ the ideal of $\Delta.$ Consider the exact sequence of monoids
$$
0\ra \Delta^{-1}(1+\calI) \xrightarrow{\lambda} \Delta^{-1}\calM_{\frakY} \xrightarrow{\Delta^{\flat}} \calM_{\frakX} \ra 0,
$$
where $\lambda=\Delta^{-1}\alpha_{\frakY}^{-1}.$
For $1\le i\le d,$ let $\widetilde{\eta}_i$ be the unique local section of $\Delta^{-1}\calI$ such that $\lambda(1+\widetilde{\eta}_i)+\Delta^{-1}p_1^{\flat}m_i=\Delta^{-1}p_2^{\flat}m_i.$ For $1\le i\le d,$ let $\xi_i$ be the image of $\frac{\widetilde{\eta}_i}{p}$ in $\calP_n.$ Then, by \ref{prop125*}, $\calP_n=\Ox_{\frakX_n}\langle \xi_1,\hdots, \xi_d\rangle.$
\end{parag}

\begin{proposition}[\cite{SF1} 11.27]\label{equivstrat}
Let $\frakX \ra \frakS$ be a log smooth morphism of logarithmic formal schemes, flat and locally of finite type over $W,$ and $P_{\frakX/\frakS,1}$ the PD envelope of the canonical immersion $\frakX \ra R_{\frakX/\frakS}.$ Then, $P_{\frakX/\frakS,1}$ has a formal groupoid structure, given locally under the hypothesis and notation of \ref{prop125}, by the morphisms
\begin{alignat*}{2}
\delta: & \frac{\eta_i}{p} \mapsto 1\otimes \frac{\eta_i}{p} + \frac{\eta_i}{p} \otimes 1 + p\frac{\eta_i}{p} \otimes \frac{\eta_i}{p}\quad \text{(comultiplication)},\\
\sigma: & \frac{\eta_i}{p} \mapsto -\frac{\eta_i}{p}+\frac{\eta_i^2}{p}-\frac{\eta_i^3}{p}+\hdots\quad \text{(antipode)} ,\qquad \pi: \frac{\eta_i}{p} \mapsto 0\quad \text{(counit)}.
\end{alignat*}
In addition, stratifications relative to the reduction $P_n$ of $P_{\frakX/\frakS,1}$ modulo $p^n,$ are equivalent to integrable and quasi-nilpotent $p$-connections. This equivalence is given locally, under the hypothesis and notation of \ref{prop125}, as follows: if $\epsilon$ is a $P_n$-stratification on an $\Ox_{\frakX_n}$-module $\calE,$ then, for all local sections $x$ of $\calE,$ there exist unique local sections $x_I$ of $\calE,$ such that $\epsilon(1\otimes x)=\sum_{I\in \N^d}x_I \otimes \xi^{[I]}.$ The corresponding $p$-connection is then given by
$$
\nabla:\calE \ra \calE \otimes_{\Ox_{\frakX_n}} \omega^1_{\frakX_n/\frakS_n},\quad x\mapsto \sum_{i=1}^dx_{\epsilon_i}\otimes \op{dlog}m_i.
$$
Conversely, given an integrable quasi-nilpotent $p$-connection $\nabla:\calE \ra \calE\otimes_{\Ox_{\frakX_n}}\omega^1_{\frakX_n/\frakS_n},$ the corresponding stratification is given by
$$
\epsilon:1\otimes x\mapsto \sum_{I\in \N^d} \nabla \left (\partial_I \right )(x)\otimes \xi^{[I]},
$$
where $\nabla\left (\partial_I \right )$ is defined in \ref{def227}.
\end{proposition}

\section{Indexed algebras and modules}

\subsection{Generalities}

\begin{parag}
We refer the reader to (\cite{SF3} 3) for a more detailed exposition of indexed objects. Let $(X,\Ox_X)$ be a ringed topos on a site $C$ and $\upmu$ a sheaf of monoids of $X.$ We have the localized topos $X_{/\upmu},$ whose objects are sheaves of $X$ equipped with a morphism to $\upmu,$ and the localization morphism of topoi $j_{\upmu}:X_{/\upmu} \ra X.$ For a sheaf $\F$ of $X,$ we denote by $\F_{\upmu}$ the localization $j_{\upmu}^{-1}\F=(\F\times \upmu \ra \upmu)$ of $\F.$ We equip $X_{/\upmu}$ with the ring $\Ox_{X,\upmu}.$ The addition and projection morphisms $\sigma,p_1,p_2:\upmu^2 \ra \upmu$ induce morphisms of ringed topoi
$$
\sigma,p_1,p_2:X_{/\upmu^2} \ra X_{/\upmu}.
$$
For modules $\calM$ and $\calN$ of $X_{/\upmu},$ we denote by $\calM \boxtimes \calN=p_1^*\calM \otimes_{\Ox_{X,\upmu^2}} p_2^*\calN$ their exterior tensor product. Finally, for a module $\calM$ of $X_{/\upmu}$ and a local section $s \in \Gamma(U,\upmu)$ over an object $U$ of $C,$ we denote by $\calM_s$ the pullback of $\calM$ along the morphism of ringed topoi $X_{/U} \ra X_{/\upmu}$ induced by $s:U\ra \upmu.$ We call it the fiber of $\calM$ over $s.$ Note that, if $X$ is the étale topos of a scheme $T,$ then $U$ is an étale $T$-scheme and $X_{/U}=U_{\text{ét}}.$
\end{parag}

\begin{definition}[\cite{SF3} 3.10]
\begin{enumerate}
\item A $\upmu$-indexed algebra is a module $\calA$ of $X_{/\upmu}$ equipped with a morphism $\calA \boxtimes \calA \ra \sigma^*\calA$ satisfying associativity, commutativity and unitarity conditions. Similarly, for a $\upmu$-indexed algebra $\calA,$ a $\upmu$-indexed $\calA$-module is a module $\calE$ of $X_{/\upmu}$ equipped with a morphism $\calA \boxtimes \calE \ra \sigma^*\calE$ satisfying associativity and unitarity conditions. 
\item A morphism of $\upmu$-indexed algebras $\calA \ra \calA'$ is an $\Ox_{X,\upmu}$-linear morphism $f:\calA \ra \calA'$ such that the diagram
$$
\begin{tikzcd}
\calA \boxtimes \calA \ar{r}{f\boxtimes f} \ar{d} & \calA'\boxtimes \calA' \ar{d} \\
\sigma^*\calA \ar{r}{\sigma^*f} & \sigma^*\calA'
\end{tikzcd}
$$
is commutative, where the vertical arrows are the $\upmu$-indexed algebra structures. We similarly define morphisms of $\upmu$-indexed modules\footnote{Note that, contrary to (\cite{SF3} 3.15), we do not consider morphisms that shift indices.}.
\end{enumerate}
\end{definition}

\begin{remark}
A $\upmu$-indexed algebra structure on an $\Ox_{X,\upmu}$-module $\calA$ is equivalent to the data of a family of morphisms
$$
\calA_s \otimes_{\Ox_U} \calA_t \ra \calA_{s+t},
$$
for every object $U$ of $C$ and $s,t\in \Gamma (U,\upmu),$ which are functorial and satisfy associativity, commutativity and unitarity conditions. We have a similar equivalence for indexed modules.
\end{remark}

\begin{exmp}\label{example314}
The main examples of an indexed algebra we will work with are those of the indexed algebra associated to a log structure and the subalgebra of horizontal sections relative to its canonical connection. Our ringed topos will be the étale topos of a logarithmic scheme $X$ and $\upmu=\ov{\calM}_X^{\mathrm{gp}}.$
\begin{enumerate}
\item We denote by $\calA_X$ the sheaf $\calM_X^{\mathrm{gp}} \times \Ox_X/\sim$ where $(u+m,a)\sim (m,\alpha_X(u)a)$ for all local sections $u,\ m$ and $a$ of $\calM_X^{\times},$ $\calM_X$ and $\Ox_X$ respectively. The canonical morphism $\calM_X^{\mathrm{gp}} \times \Ox_X \ra \calM_X^{\mathrm{gp}} \ra \upmu$ factors through $\calA_X$ so that $\calA_X$ can be considered as an object of $X_{\text{ét}/\upmu}.$ If $s$ is a local section of $\calM_X^{\mathrm{gp}}$ over an étale $X$-scheme $U,$ $\ov{s}$ its image in $\upmu$ and $\calM_{X,\ov{s}}^{\mathrm{gp}}$ the subsheaf of $\calM_U^{\mathrm{gp}}$ consisting of local sections over $\ov{s},$ then $\calA_{X,\ov{s}}=\calM_{X,\ov{s}}^{\mathrm{gp}}\times \Ox_U/\sim.$ The $\Ox_U$-module $\calA_{X,\ov{s}}$ is free of rank one and has $e_s=(s,1)$ as a basis. The morphisms
$$
\calA_{X,\ov{s}} \otimes_{\Ox_U} \calA_{X,\ov{t}} \ra \calA_{X,\ov{s+t}},\quad ae_s \otimes be_t \mapsto abe_{s+t},
$$
for local sections $s$ and $t$ of $\calM_X^{\mathrm{gp}}$ and $a$ and $b$ of $\Ox_X,$ define a $\upmu$-indexed algebra structure on $\calA_X.$
\item Suppose we have a morphism of logarithmic schemes $X\ra S.$ By (\cite{Lor2000} II 1.6), there exists a unique connection $d_{\calA_X}:\calA_X \ra \calA_X \otimes_{\Ox_{X,\upmu}} \omega^1_{X/S,\upmu}$ satisfying
\begin{equation}\label{dAX}
d_{\calA_X,\ov{s}}(e_s)=e_s\otimes \op{dlog}s,
\end{equation}
for all local sections $s$ of $\calM_X^{\mathrm{gp}}.$ In addition, by (\cite{Lor2000} II 1.7), $d_{\calA_X}$ is integrable and, if $S$ is of characteristic $p,$ then $d_{\calA_X}$ has vanishing $p$-curvature. The subsheaf $\calB_{X/S}=\op{ker}d_{\calA_X}$ inherits a $\upmu$-indexed algebra structure from $\calA_X.$ 
\end{enumerate}
\end{exmp}

\begin{parag}\label{strictA}
A strict morphism of logarithmic schemes $f:X\ra Y$ induces an isomorphism $f^{-1}\ov{\calM}_Y^{\mathrm{gp}} \isom \ov{\calM}_X^{\mathrm{gp}}$ and then an isomorphism of $\upmu$-indexed algebras $f^*\calA_Y \isom \calA_X.$ 
\end{parag}

\begin{proposition}\label{prop01}
Let $g:X\ra Y$ be an étale and strict morphism of logarithmic $S$-schemes in characteristic $p$ and $g':X' \ra Y'$ its base change along the Frobenius $F_S$ in the category of fs logarithmic schemes. Then, $g$ and $g'$ induce isomorphisms of indexed algebras
$$
g^*\calA_Y \isom \calA_X,\quad g'^*\calB_{Y/S} \isom \calB_{X/S}.
$$
\end{proposition}

\begin{proof}
We have an exact sequence of $\Ox_{Y',\upmu}$-modules
$$
0\ra \calB_{Y/S} \ra F_{Y/S*} \calA_Y \xrightarrow{F_{Y/S*}d_{\calA_Y}} F_{Y/S*}\left ( \calA_Y \otimes_{\Ox_{Y}}\omega^1_{Y/S} \right ).
$$
The morphism $g$ is étale hence so is $g'.$ It follows that $g'$ is flat and we obtain an exact sequence 
$$
0\ra g'^*\calB_{Y/S} \ra g'^*F_{Y/S*} \calA_Y \xrightarrow{g'^*F_{Y/S*}d_{\calA_Y}} g'^*F_{Y/S*}\left ( \calA_Y \otimes_{\Ox_{Y}}\omega^1_{Y/S} \right ).
$$
By the cartesian square
$$
\begin{tikzcd}
X \ar{r}{F_{X/S}} \ar[swap]{d}{g} & X' \ar{d}{g'}  \\
Y \ar{r}{F_{Y/S}} & Y',
\end{tikzcd}
$$
we have a canonical isomorphism $g'^*F_{Y/S*}=F_{X/S*}g^*.$ Also, since $g$ is log étale, we have a canonical isomorphism $g^*\omega^1_{Y/S}=\omega^1_{X/S}.$ It follows that we have the following exact sequence of $\Ox_{X',\upmu}$-modules:
$$
0\ra g'^*\calB_{Y/S} \ra F_{X/S*} \calA_X \xrightarrow{F_{X/S*}d_{\calA_X}} F_{X/S*}\left ( \calA_X \otimes_{\Ox_{X}}\omega^1_{X/S} \right ).
$$
We conclude that $g'^*\calB_{Y/S}$ and $\calB_{X/S}$ are isomorphic.
\end{proof}

\begin{definition}
Let $(X,\Ox_X)$ be a ringed topos on a site $C$ and $\upmu$ a sheaf of monoids of $X.$ For a $\upmu$-indexed algebra $\calA$ and $\upmu$-indexed $\calA$-modules $\calM$ and $\calN$, we define their tensor product $\calM\circledast_{\calA}\calN$ to be the sheaf of $X$ associated to the presheaf
$$
U \mapsto \coprod_{s\in \Gamma(U,\upmu)} \bigoplus_{\substack{t,r\in \Gamma(U,\upmu) \\ t+r=s}} \calM_t(U) \otimes_{\Ox_X(U)} \calN_r(U)/\sim,
$$
equipped with the natural map to $\upmu,$
where $(ax) \otimes y\sim x\otimes (ay)$ for all $t,r,u\in \Gamma(U,\upmu),$ $a\in \Gamma(U,\calA_u),$ $x\in \Gamma(U,\calM_t)$ and $y\in \Gamma(U,\calN_r)$ such that $t+r+u=s.$
\end{definition}

\begin{remark}\label{rem318}
Let $(X,\Ox_X)$ be a ringed topos on a site $C,$ $\upmu$ a sheaf of monoids of $X,$ $\calA$ a $\upmu$-indexed algebra and $\calM$ and $\calN$ two $\upmu$-indexed $\calA$-modules. Let $\sigma:\upmu^2\ra \upmu$ be the addition map and $\sigma_!$ the left adjoint of the pullback functor $\sigma^*:X_{/\upmu} \ra X_{/\upmu^2}.$ The sheaf $\sigma_! \left (\calM \boxtimes \calN \right )$ is the sheaf on $C$ associated to the presheaf
$$
U\mapsto \coprod_{s\in \Gamma(U,\upmu)} \bigoplus_{\substack{u,v\in \Gamma(U,\upmu)\\ u+v=s}} \calM_{u}(U) \otimes_{\Ox_X(U)} \calN_v(U),
$$
equipped with the canonical morphism to $\upmu.$ If $\tau:\upmu^3\ra \upmu$ is the addition map, then the $\upmu$-indexed $\calA$-algebra structures on $\calM$ and $\calN$ induce two natural morphisms
$\tau_! \left ( \calA \boxtimes \calM \boxtimes \calN \right ) \ra \sigma_! \left ( \calM \boxtimes \calN \right ).$
The tensor product $\calM \circledast_{\calA} \calN$ is the coequalizer of these two morphisms i.e. we have an exact sequence
\begin{equation}\label{exseq4182}
\tau_! \left ( \calA \boxtimes \calM \boxtimes \calN \right ) \ra \sigma_! \left ( \calM \boxtimes \calN \right ) \ra \calM\circledast_{\calA} \calN \ra 0.
\end{equation}
\end{remark}

\subsection{The local model}

\begin{parag}\label{parag221}
Let $f:X\ra S$ be a log smooth morphism of logarithmic schemes of characteristic $p.$ Set $\upmu=\ov{\calM}_X^{\mathrm{gp}}$ and recall the $\upmu$-indexed algebras $\calA_X$ and $\calB_{X/S}$ (\ref{example314}). By (\cite{Ogus2018} IV 3.3.1), étale locally around a geometric point $\ov{x} \ra X$ and its image $f(\ov{x}) \ra S,$ there exists an fs chart $\theta:P\ra Q$ of $X\ra S$ satisfying the three following conditions:
\begin{itemize}
\item $\theta^{\mathrm{gp}}$ is injective and the torsion subgroup of its cokernel is finite of order invertible in $\Ox_X.$ 
\item The morphism $X\ra S\times_{A[P]}A[Q],$ induced by $X\ra S$ and the chart $X\ra A[Q],$ is étale and strict.
\item The chart $X\ra A[Q]$ is exact at $\ov{x}$ i.e. $\ov{Q} \ra \ov{\calM}_{X,\ov{x}}$ is an isomorphism.
\end{itemize}
Set $T=S\times_{A[P]}A[Q]$ and $g:X\ra T$ the morphism given above. In this case, we have
\begin{equation}
T'=S\times_{S,F_S}T=S\times_{S,F_S}S\times_{A[P]}A[Q]=S\times_{A[P],F_{A[P]}}A[P]\times_{A[P]}A[Q]=S\times_{A[P]}A[Q'],
\end{equation}
where $Q'$ is defined in \ref{parag111}.
Let $q\in Q^{\mathrm{gp}}$ and consider it as a section of $\calM_X^{\mathrm{gp}}$ via the chart $Q\ra \calM_X.$ Denote by $\ov{q}$ the image of $q$ in $\upmu_T=\ov{\calM}_T^{\mathrm{gp}}.$ Then
\begin{equation}\label{eq2211}
\calA_{T,\ov{q}}=\Ox_Te_q=\left ( \Ox_S \otimes_{\mathbb{F}_p[P]} \mathbb{F}_p[Q] \right )e_q,
\end{equation}
where $e_q$ is given in \ref{example314}. We also have
$$
\omega^1_{T/S}=\Ox_T \otimes_{\Z[Q]} \omega^1_{A[Q] / A[P]}=\Ox_T \otimes_{\Z} \left ( Q^{\mathrm{gp}}/P^{\mathrm{gp}} \right )=\Ox_T\otimes_{\mathbb{F}_p} C,
$$
where $C$ is the cokernel of $F_{Q/P}^{\mathrm{gp}}:Q^{\mathrm{gp}} \oplus_{P^{\mathrm{gp}},F_P^{\mathrm{gp}}}P^{\mathrm{gp}} \ra Q^{\mathrm{gp}},\ (x,t) \mapsto px+\theta^{\mathrm{gp}}(t).$
Recalling the notation \ref{not3}, the connection $d_{\calA_T,\ov{q}}$ \eqref{dAX} is given by
$$
d_{\calA_T,\ov{q}}: \begin{array}[t]{clc}
\left ( \Ox_S \otimes_{\mathbb{F}_p[P]}\mathbb{F}_p[Q] \right )e_q & \ra & \left ( \Ox_S \otimes_{\mathbb{F}_p[P]} \mathbb{F}_p[Q] \right ) \otimes_{\mathbb{F}_p}C \\
(a\otimes t^x )e_q & \mapsto & a\otimes t^x \otimes (x+q),
\end{array}
$$
for $x\in Q$ and $a$ a local section of $\Ox_S.$ For $c\in C$ or $c\in Q^{\mathrm{gp}},$ define the $\mathbb{F}_p[P]$-subspace (or even $\mbbF_p[Q']$-subspace) of $\mathbb{F}_p[Q]:$
$$
\mathbb{F}_p[Q]_c=\bigoplus_{\substack{a\in Q\\ a=c\ \text{in}\ C}} \mathbb{F}_pt^a.
$$
We then have direct sums
$$
\mathbb{F}_p[Q]=\bigoplus_{c\in C}\mathbb{F}_p[Q]_c,\quad \calA_{T,\ov{q}}=\bigoplus_{c\in C} \left ( \Ox_S \otimes_{\mathbb{F}_p[P]} \mathbb{F}_p[Q]_c\right ) e_q.
$$
For $c\in C$ or $c\in Q^{\mathrm{gp}},$ let $Q_c\subset Q$ be the submonoid consisting of elements $x\in Q$ that are equal to $c$ in $C.$ Let
$$
b=\sum_{x\in Q}(a_x\otimes t^x)e_q=\sum_{c\in C} \sum_{x\in Q_c}(a_x\otimes t^x)e_q\in \calB_{T/S,\ov{q}},
$$
where the $a_x$ are local sections of $\Ox_S.$ Then
$$
0=d_{\calA_T,\ov{q}}(b)=\sum_{c\in C}\sum_{x\in Q_c} a_x\otimes t^x\otimes (q+c) \in \bigoplus_{c\in C} \left ( \Ox_S\otimes_{\mathbb{F}_p[P]}\mathbb{F}_p[Q]_c  \otimes (q+c)\right ).
$$
It follows that for all $c\in C$ such that $q+c \neq 0$ in $C,$ we have
$
\sum_{x\in Q_c} a_x \otimes t^x=0,
$
and so $b=\sum_{x\in Q_{-q}}a_x\otimes t^x e_q.$
We conclude that
\begin{equation}\label{Blocdesc}
\calB_{T/S,\ov{q}}=\left (\Ox_S\otimes_{\mathbb{F}_p[P]}\mathbb{F}_p[Q]_{-q}\right )e_q=\left ( \Ox_S \otimes_{\mbbF_p[P]} \mbbF_p[Q'] \otimes_{\mbbF_p[Q']}\mbbF_p[Q]_{-q} \right ) e_q=\left ( \Ox_{T'}\otimes_{\mbbF_p[Q']}\mbbF_p[Q]_{-q} \right ) e_q.
\end{equation}
By \ref{prop01}, we get
\begin{equation}\label{eq2214}
\calB_{X/S,\ov{q}}=\Ox_{X'}\otimes_{\Ox_{T'}}\calB_{T/S,\ov{q}}=\left ( \Ox_{X'} \otimes_{\mbbF_p[Q']}\mbbF_p[Q]_{-q} \right )e_q.
\end{equation}
\end{parag}

\begin{theorem}[Lorenzon]\label{thmlor}
If $X\ra S$ is a log smooth morphism of logarithmic schemes of characteristic $p,$ then the $\upmu$-indexed $\calB_{X/S}$-module $\calA_X$ is locally free of finite rank equal to $p^d,$ where $d$ is the rank of $\omega^1_{X/S}.$
\end{theorem}

\begin{proof}
The question being local, we may assume that $X\ra S$ has a chart as in \ref{parag221}. Keep the notation of \ref{parag221}. For every $c\in C,$ choose a lifting $q_c\in Q^{\mathrm{gp}}$ of $c.$
Then, by \eqref{eq2211} and \eqref{Blocdesc}, we have
$$
\calA_{T,\ov{q}}=\left ( \Ox_S \otimes_{\mathbb{F}_p[P]}\mathbb{F}_p[Q] \right )e_q=\bigoplus_{c\in C} \left ( \Ox_S\otimes_{\mathbb{F}_p[P]}\mathbb{F}_p[Q]_c \right )e_q=\bigoplus_{c\in C}\calB_{T/S,-\ov{q}_c}e_{q+q_c}.
$$
Denote by $\calB_{T/S}(\ov{q}_c)$ the sheaf $\calB_{T/S}$ but with indices shifted by $\ov{q}_c.$ We obtain an isomorphism
$$
F_{T/S*}\calA_T \isom \bigoplus_{c\in C}\calB_{T/S}(\ov{q}_c).
$$
Pulling back along the morphism $g':X'\ra T'$ induced by $g$ and using \ref{prop01} and the cartesian square
$$
\begin{tikzcd}
X \ar{r}{g} \ar[swap]{d}{F_{X/S}} & T \ar{d}{F_{T/S}} \\
X' \ar{r}{g'} & T',
\end{tikzcd}
$$
we get the isomorphism
$
F_{X/S*}\calA_X \isom \bigoplus_{c\in C}\calB_{X/S}(\ov{q}_c).
$
Since $\omega^1_{X/S}=\Ox_X\otimes_{\mathbb{F}_p}C,$ we get that $\op{dim}_{\mathbb{F}_p}C=\op{rk}_{\Ox_X}\omega^1_{X/S}=d,$ so $|C|=p^d$ and the dimension follows.
\end{proof}

\begin{proposition}\label{prop02}
The indexed algebra $\calB_{X/S}$ is generated by $\Ox_{X'} \subset \calB_{X/S,0}$ and local sections of the form $\alpha_X(m)e_{-m}$ and $e_{F_{X/S}^{\flat}m'}$ where $m$ and $m'$ are local sections of $\calM_X$ and $\calM_{X'}^{\mathrm{gp}}$ respectively.
\end{proposition}

\begin{proof}
By definition of $d_{\calA_X}$ \eqref{dAX}, it is clear that local sections of one of these forms is a local section of $\calB_{X/S}.$ For the converse, the question being local, we may assume that we have a chart $P\ra Q$ as in \ref{parag221}. We keep the notation we introduced in \ref{parag221}. Proposition \ref{prop01} allows us then to reduce to the case $X=T.$ By the local description \eqref{Blocdesc}, it is sufficient to prove that, for $q\in Q^{\mathrm{gp}}$ and $x\in Q_{-q}$, the local section $t^xe_q \in \mathbb{F}_p[Q]_{-q}e_q$ can be expressed in terms of the sections in the proposition. Recall that we have a chart $F_{Q/P}$ \eqref{FQP} of $F_{X/S}.$ Since $x+q=0$ in $C=\op{coker} F_{Q/P}^{\mathrm{gp}},$ there exist $(a,b)\in Q^{\mathrm{gp}}\oplus_{P^{\mathrm{gp}},F_P^{\mathrm{gp}}}P^{\mathrm{gp}}$ such that $x+q=F_{Q/P}^{\mathrm{gp}}(a,b)$ in $Q.$ Then
$$
t^xe_q=t^xe_{-x}e_{F_{Q/P}^{\mathrm{gp}}(a,b)}.
$$
This concludes the proof as $t^xe_{-x}$ is of the form $\alpha_X(m)e_{-m}.$
\end{proof}

\begin{lemma}\label{lem224}
Let $\beta:Q\ra \calM_X$ be a chart on a logarithmic scheme $X$ and $\ov{x}\ra X$ a geometric point. Then $\calM_{X,\ov{x}}=\alpha_{X,\ov{x}}^{-1} \left ( \Ox_{X,\ov{x}}^* \right ) + \beta(Q).$
\end{lemma}

\begin{proof}
This is a direct consequence of the definition of a chart: $\calM_X$ is the logarithmic structure associated to the prelogarithmic structure $\gamma:Q\xrightarrow{\beta} \calM_X \xrightarrow{\alpha_X} \Ox_X$ so $\calM_X=Q\oplus_{\gamma^{-1}\Ox_X^*}\Ox_X^*$ and $\beta(a)=(a,1)$ and $\alpha_X(a,b)=\beta(a)+\alpha_X^{-1}(b)$ for $a\in Q$ and $b$ a local section of $\Ox_X^*.$
\end{proof}

\begin{lemma}\label{lem225}
Let $\beta:Q\ra \calM_X$ be a chart on a logarithmic scheme $X.$ If $\beta$ is exact at a geometric point $\ov{x}$ i.e. $\ov{Q} \ra \ov{\calM}_{X,\ov{x}}$ is an isomorphism, then the stalk of the canonical morphism
\begin{equation}\label{eq2251}
\Ox_X \otimes_{\Z[\beta(Q)]} \Z[\beta(Q)^{\mathrm{gp}}] \ra \Ox_X \otimes_{\Z[\calM_X]} \Z[\calM_X^{\mathrm{gp}}]
\end{equation}
at $\ov{x}$ is an isomorphism of $\Ox_{X,\ov{x}}$-algebras.
\end{lemma}

\begin{proof}
By \ref{lem224}, $\calM_{X,\ov{x}}=\alpha_{X,\ov{x}}^{-1}(\Ox_{X,\ov{x}}^*)+\beta(Q).$ For simplicity, we denote this by $\calM_{X,\ov{x}}=\Ox_{X,\ov{x}}^*\beta(Q).$ Set
$$
g:\Ox_{X,\ov{x}} \otimes_{\Z[\Ox_{X,\ov{x}}^*\beta(Q)]}\Z[\Ox_{X,\ov{x}}^*\beta(Q)^{\mathrm{gp}}] \ra \Ox_X \otimes_{\Z[\beta(Q)]} \Z[\beta(Q)^{\mathrm{gp}}],\quad 1\otimes t^{\alpha_{X,\ov{x}}^{-1}(a)+\beta(m)} \mapsto a\otimes t^{\beta(m)},
$$
for $m\in Q^{\mathrm{gp}}$ and $a\in \Ox_{X,\ov{x}}^*.$ Let us check that $g$ is well-defined: suppose $m,m'\in Q^{\mathrm{gp}}$ and $a,a'\in \Ox_{X,\ov{x}}^*$ such that $\alpha_{X,\ov{x}}^{-1}(a)+\beta(m)=\alpha_{X,\ov{x}}^{-1}(a')+\beta(m')$ in $\calM_{X,\ov{x}}^{\mathrm{gp}}.$ Since $\beta(m-m')\in \calM_{X,\ov{x}}$ and the chart $\beta$ is exact at $\ov{x},$ we have $m-m'\in Q$ and $\alpha_{X,\ov{x}}^{-1}(a'a^{-1})=\beta(m-m')\in \beta(Q).$ Then, in $\Ox_{X,\ov{x}} \otimes_{\Z[\beta(Q)]} \Z[\beta(Q)^{\mathrm{gp}}],$ we have
$$
a\otimes t^{\beta(m)} = a\otimes t^{\beta(m')}t^{\alpha_{X,\ov{x}}^{-1}(a'a^{-1})}=a'\otimes t^{\beta(m')}.
$$
It follows that $g$ is well-defined. It is clearly an inverse to the stalk of \eqref{eq2251} at $\ov{x}.$
\end{proof}

\begin{theorem}\label{thm326}
Let $f:X \ra S$ be a log smooth morphism of logarithmic schemes in characteristic $p$ and $s$ a section of $\ov{\calM}_X^{\mathrm{gp}}.$ Set $\calP=\calM_{X'}^{\mathrm{gp}} \oplus \ov{\calM}_X$ and denote by $\calP_s$ the subsheaf of $\calP$ consisting of local sections $(m',m)$ of $\calP$ such that $F_{X/S}^{\flat}(m')-m=s.$ Equip $\calP_s$ with an action of $\calM_{X'}$ by setting
$$
a \cdot (m',\ov{m}) = (a+m',\ov{m}+F_{X/S}^{\flat}(a)),
$$
for local sections $a,$ $m'$ and $\ov{m}$ of $\calM_{X'},$ $\calM_{X'}^{\mathrm{gp}}$ and $\ov{\calM}_X$ respectively sarisfying $F_{X/S}^{\flat}(m')-m=s.$
Then, the $\Ox_{X'}$-linear morphism
$$
\varphi:\Ox_{X'} \otimes_{\mbbF_p[\calM_{X'}]}\mbbF_p[\calP_s] \ra \calB_{X/S,s},\quad 1\otimes t^{(m',m)} \mapsto \alpha_X(m)e_{F_{X/S}^{\flat}m'-m},
$$
is an isomorphism.
\end{theorem}

\begin{proof}
Take a geometric point $\ov{x} \ra X.$ The question being local, we may suppose that we have a chart $\theta:P\ra Q$ of $X\ra S$ satisfying the conditions in \ref{parag221}. Denote by $\beta:Q\ra \calM_X$ the chart on $X$ and $\beta':Q'\ra \calM_{X'}$ the deduced chart on $X'.$ Recall the chart $F_{Q/P}:Q'\ra Q$ of the exact relative Frobenius $F_{X/S}$ (\ref{parag111}). Since $F_{Q/P}$ is exact by \ref{prop212}, $\beta$ is exact at $\ov{x}$ and $\calM_{X',\ov{x}} \ra \calM_{X,\ov{x}}$ is exact, the chart $\beta'$ is also exact at $\ov{x}.$ The charts $\beta$ and $\beta'$ then induce isomorphisms $\ov{Q} \isom \ov{\calM}_{X,\ov{x}}$ and $\ov{Q}' \isom \ov{\calM}_{X',\ov{x}}.$ By \ref{lem225}, we have a canonical isomorphism
$$
\Ox_{X',\ov{x}} \otimes_{\mbbF_p[\calM_{X',\ov{x}}]} \mbbF_p[\calP_{\ov{x}}] \isom \Ox_{X',\ov{x}} \otimes_{\mbbF_p[\beta'(Q')]} \mbbF_p[\beta'(Q'^{\mathrm{gp}})\oplus \ov{Q}].
$$
Under this isomorphism, we have the following identification
$$
\Ox_{X',\ov{x}} \otimes_{\mbbF_p[\calM_{X',\ov{x}}]} \mbbF_p[\calP_{s,\ov{x}}] \isom \Ox_{X',\ov{x}} \otimes_{\mbbF_p[\beta'(Q')]} \mbbF_p[\left ( \beta'(Q'^{\mathrm{gp}})\oplus \ov{Q}\right )_s],
$$
where $\left ( \beta'(Q'^{\mathrm{gp}})\oplus \ov{Q}\right )_s$ is the submonoid of $\beta'(Q'^{\mathrm{gp}})\oplus \ov{Q}$ consisting of elements $(\beta'(m'),\ov{a})$ such that $F_{Q/P}(m')-a=s$ in $\ov{Q}^{\mathrm{gp}}.$ Take a lift $q\in Q^{\mathrm{gp}}$ of $s.$ By \eqref{eq2214}, we have
$$
\calB_{X/S,s}=\left ( \Ox_{X'} \otimes_{\mbbF_p[Q']}\mbbF_p[Q]_{-q} \right )e_q.
$$
Denote by $C$ the cokernel of $F_{Q/P}^{\mathrm{gp}}.$ For $b\in Q$ such that $b+q=0$ in $C,$ there exists $t_b\in Q'^{\mathrm{gp}}$ such that $b+q=F_{Q/P}^{\mathrm{gp}}(t),$ which is unique by injectivity of $F_{Q/P}^{\mathrm{gp}}$ (\ref{prop212}). It follows that, in $\ov{Q}^{\mathrm{gp}},$ we have
$$
F_{Q/P}^{\mathrm{gp}}(t_b)-b=q=s,
$$
so $(\beta'(t_b),\ov{b})\in \left ( \beta'(Q'^{\mathrm{gp}})\oplus \ov{Q}\right )_s$ and we have a morphism
$$
\psi:\begin{array}[t]{clc}
\left ( \Ox_{X',\ov{x}} \otimes_{\mbbF_p[Q']}\mbbF_p[Q]_{-q} \right )e_q=\calB_{X/S,s,\ov{x}} & \ra & \Ox_{X',\ov{x}} \otimes_{\mbbF_p[\beta'(Q')]} \mbbF_p[\left ( \beta'(Q'^{\mathrm{gp}})\oplus \ov{Q}\right )_s] \\
(a\otimes t^b) e_q & \mapsto & a \otimes t^{(\beta'(t_b),\ov{b})}.
\end{array}
$$
The stalk of $\varphi$ at $\ov{x}$ and $\psi$ are inverse to each other.
\end{proof}

\subsection{A key isomorphism}

\begin{parag}
In this subsection, we consider a log smooth morphism $X\ra S$ of logarithmic schemes of characteristic $p.$ We denote by $F_{X/S}:X\ra X'$ the exact relative Frobenius and set $Y=X\times_{X'}X.$ Let $p_1,p_2:Y\ra X$ be the projections. Recall that the diagonal immersion $X\ra Y$ is a universal homeomorphism (\ref{prop114}). We use it to identify the étale sites of $X$ and $Y.$ We also identify the étale sites of $X$ and $X'$ via the $F_{X/S}.$ Set $\upmu=\ov{\calM}_X^{\mathrm{gp}},$ $\sigma:\upmu^2\ra \upmu$ the addition morphism and recall the $\upmu$-indexed algebras $\calA_X,$ $\calB_{X/S}$ and $\calA_Y$ (\ref{example314}). Recall also the sections $e_s$ given in \ref{example314}.
Etale locally on $X$ and $S,$ the following hypothesis is satisfied:

\textbf{Hypothesis $(H)$}: there exists an fs chart $\theta:P \ra Q$ of $X\ra S$ satisfying the following conditions:
\begin{itemize}
\item $\theta^{\mathrm{gp}}$ is injective and the torsion subgroup of $\op{coker} \theta^{\mathrm{gp}}$ is finite of order coprime with $p.$
\item The morphism $X\ra S\times _{A[P]}A[Q],$ induced by the chart $X\ra A[Q]$ and $X\ra S,$ is étale and strict.
\end{itemize}
Set $T=S\times_{A[P]}A[Q],$ $Y=X\times_{X'}X$ and $Z=T\times_{T'}T$ and let $g:X\ra T$ be the morphism given above. Since $\theta^{\mathrm{gp}}$ is injective and the monoids $P$ and $Q$ are integral, $\theta$ is also injective and we identify $P$ with its image $\theta(P).$ Recall that $F_{Q/P}:Q' \ra Q$ \eqref{FQP} is a chart of the Frobenius $F_{X/S}.$
\end{parag}

\begin{lemma}\label{lem01}
Suppose that $(H)$ is satisfied and $X=T.$ Let $p_i:\calA_X \ra \calA_Y$ be the morphism of indexed algebras induced by the projection $p_i.$ Then $p_1$ and $p_2$ are equal on $\calB_{X/S}.$
\end{lemma}

\begin{proof}
Let $C=Q^{\mathrm{gp}} / (pQ^{\mathrm{gp}}+P^{\mathrm{gp}})$ be the cokernel of $F_{Q/P}^{\mathrm{gp}}$ and $q\in Q^{\mathrm{gp}}$ and denote by $\ov{q}$ its image in $\ov{Q}^{\mathrm{gp}}.$ Since $X=T,$ we have $Y=T\times_{T'}T=S\times_{A[P]}A[\left ( Q\oplus_{Q'}Q \right )^{\mathrm{sat}}].$ By \ref{lem112}, $Y=S\times_{A[P]}A[ Q\oplus C].$ Then, by \eqref{eq2211} and \eqref{Blocdesc}, we have
$$\calA_{X,\ov{q}}=\left ( \Ox_S \otimes_{\mathbb{F}_p[P]}\mathbb{F}_p[Q] \right )e_q,\quad \calB_{X/S,\ov{q}}=\left (\Ox_S \otimes_{\mbbF_p[P]}\mbbF_p[Q]_{-q} \right )e_q,\quad 
\calA_{Y,\ov{q}}=\left ( \Ox_S\otimes_{\mathbb{F}_p[P]}\mathbb{F}_p[Q\oplus C] \right ) e_{(q,0)},
$$
where $\mbbF_p[Q]_{-q}=\bigoplus_{\substack{x\in Q\\ x+q=0\ \text{in}\ C}}\mbbF_pt^x.$ By \ref{rem113}, the projections $p_1,p_2$ are given, for all $a\in Q,$ by
$$
p_1((1\otimes t^a)e_q)=(1\otimes t^{(a,0)})e_{(q,0)},\quad p_2((1\otimes t^a)e_q)=(1\otimes t^{(a,a)})e_{(q,q)}.
$$
Since $e_{(q,q)}=t^{(0,q)}e_{(q,0)}$ in $\calA_{Y,\ov{q}},$ we get $p_2(t^ae_q)=t^{(a,a+q)}e_{(q,0)}.$ We deduce that, if $a+q=0$ in $C,$ then
$$
p_1((1\otimes t^a)e_q)=p_2((1\otimes t^a)e_q).
$$
The morphisms $p_1$ and $p_2$ are clearly equal on $\Ox_S,$ hence the result.
\end{proof}

\begin{lemma}\label{lem03}
Suppose that $(H)$ is satisfied. The morphism $\calA_T \boxtimes \calA_T \ra \sigma^*\calA_Z,$ defined by the projections $p_1,p_2:Z\ra T,$ induces an isomorphism of indexed algebras
\begin{equation}
\varphi:\calA_{T} \circledast_{\calB_{T/S}} \calA_T \isom \calA_Z.
\end{equation}
\end{lemma}

\begin{proof}
By \ref{lem01}, $\varphi$ is well-defined as a morphism of indexed algebras. Let $q\in Q^{\mathrm{gp}}$ and $\ov{q}$ its image in $\ov{Q}^{\mathrm{gp}}.$ By \ref{lem112}, we have
$$Z=T\times_{T'}T=S\times_{A[P]}A[ \left (Q\oplus_{Q'}Q\right )^{\mathrm{sat}}]=S\times_{A[P]}A[Q\oplus C],$$
where $C$ is the cokernel of $F_{Q/P}^{\mathrm{gp}}$ \eqref{FQP}. Then
$$
\calA_{Z,\ov{q}}=\left (\Ox_S \otimes_{\mbbF_p[P]}\mathbb{F}_p[Q\oplus C] \right )e_{(q,0)}=\Ox_S \otimes_{\mbbF_p[P]} \mathbb{F}_p[Q] \otimes_{\mathbb{F}_p} \mathbb{F}_p[C] e_{(q,0)},\quad \omega^1_{T/S}=\Ox_T \otimes_{\mathbb{F}_p}C.
$$
Set $d=\op{rk}_{\Ox_T}\omega^1_{T/S}.$ Then $\op{dim}_{\mathbb{F}_p}C=d$ and the cardinal of $C$ is $|C|=p^d.$ It follows that the rank of $\calA_{Z,\ov{q}}$ as an $\left ( \Ox_S \otimes_{\mbbF_p[P]} \mathbb{F}_p[Q] \right )$-module is finite and equal to $p^d.$ On the other hand, by \ref{thmlor}, $\calA_T$ is a locally free indexed $\calB_{T/S}$-module of rank $p^d.$ It follows that
$\calA_T\circledast_{\calB_{T/S}}\calA_T$ is a locally free indexed $\calA_T$-module of rank $p^d.$
The $\left (\Ox_S \otimes_{\mbbF_p[P]} \mathbb{F}_p[Q] \right )$-module $\left ( \calA_T\circledast_{\calB_{T/S}} \calA_T \right )_{\ov{q}}$ is hence free of rank $p^d.$ It is thus sufficient to prove that $\varphi$ is surjective. Let $(a,b) \in Q\oplus C.$ Let $x=a\in Q$ and $y\in Q^{\mathrm{gp}}$ that lifts $b\in C.$ We then have $(1\otimes t^x)e_{q-y}\in \left (\Ox_S \otimes_{\mbbF_p[P]} \mathbb{F}_p[Q] \right )e_{q-y}=\calA_{T,\ov{q}-\ov{y}},$ $e_y\in \left (\Ox_S \otimes_{\mbbF_p[P]} \mathbb{F}_p[Q] \right )e_y=\calA_{T,\ov{y}}$ and
$$
p_1((1\otimes t^x)e_{q-y})p_2(e_y)=(1\otimes t^{(a,0)})e_{(q,b)} =(1\otimes t^{(a,b)})e_{(q,0)} \in \left (\Ox_S \otimes_{\mbbF_p[P]}\mathbb{F}_p[Q\oplus C] \right )e_{(q,0)}=\calA_{Z,\ov{q}}.
$$
\end{proof}

\begin{theorem}\label{th234}
The morphism $\calA_X \boxtimes \calA_X \ra \sigma^*\calA_Y,$ defined by the projections $p_1,p_2:Y=X\times_{X'}X\ra X,$ induces an isomorphism of indexed algebras
$$
\calA_{X} \circledast_{\calB_{X/S}} \calA_X \isom \calA_Y.
$$
\end{theorem}

\begin{proof}
When forming the indexed tensor product $\calA_{X} \circledast_{\calB_{X/S}} \calA_X,$ both $\calA_X$ are considered as $\Ox_{X',\upmu}$-modules via the relative Frobenius $F_{X/S}.$ Similarly, $\calA_Y$ is considered as an $\Ox_{X',\upmu}$-module via the canonical morphism $\varphi:Y\ra X'.$ So the isomorphism we have to prove is actually
$$
\left (F_{X/S*}\calA_{X} \right )\circledast_{\calB_{X/S}} \left (F_{X/S*}\calA_X \right ) \isom \varphi_*\calA_Y.
$$
We may work étale locally on $X$ and $S$ and suppose that $(H)$ is satisfied. Denote by $\varphi':Z\ra T'$ the canonical morphism, by $g':X'\ra T'$ the base change of $g$ along the Frobenius $F_S$ and by $h:Y\ra Z=T\times_{T'}T$ the morphism induced by $g$ and $g'.$ By \ref{lem03}, the projections $Z\ra T$ induce an isomorphism
$$
\left (F_{T/S*}\calA_{T} \right ) \circledast_{\calB_{T/S}} \left ( F_{T/S*} \calA_T \right ) \isom \varphi'_*\calA_Z.
$$
Applying $g'^*$ and (\cite{SF3} 3.30), we get
$$
\left (g'^*F_{T/S*}\calA_{T} \right ) \circledast_{g'^*\calB_{T/S}} \left (g'^* F_{T/S*} \calA_T \right ) \isom g'^*\varphi'_*\calA_Z.
$$
By \ref{prop01} and the cartesian square
$$
\begin{tikzcd}
Y \ar{r}{\varphi}  \ar[swap]{d}{h} & X' \ar{d}{g'} \\
Z \ar{r}{\varphi'} & T',
\end{tikzcd}
$$
we get a canonical isomorphism
$$
\left (F_{X/S*}\calA_X \right ) \circledast_{\calB_{X/S}} \left (F_{X/S*}\calA_X \right ) \isom \varphi_* h^*\calA_Z.
$$
Given that $h$ is strict, we deduce that
$
\left (F_{X/S*}\calA_{X} \right )\circledast_{\calB_{X/S}} \left (F_{X/S*}\calA_X \right ) \isom \varphi_*\calA_Y.
$
\end{proof}

\subsection{Modulo $p^n$}

\begin{parag}
In this subsection, we fix a positive integer $n$ and consider a log smooth morphism $f:\widetilde{X} \ra \widetilde{S}$ of logarithmic schemes, flat over $W/p^nW.$ We denote by $f_1:X\ra S$ its reduction modulo $p$ and suppose that the exact relative Frobenius $F_{X/S}:X\ra X'$ lifts to an $\widetilde{S}$-morphism $\widetilde{F}:\widetilde{X} \ra \widetilde{X}',$ such that $\widetilde{X}'$ is log smooth over $\widetilde{S}.$ We identify the étale sites of $X$ and $\widetilde{X}$ (resp. $X'$ and $\widetilde{X}'$) via the nilpotent immersion $X\ra \widetilde{X}$ (resp. $X'\ra \widetilde{X}'$) and set $\upmu=\ov{\calM}_X^{\mathrm{gp}},$ which identifies with $\ov{\calM}_{\widetilde{X}}^{\mathrm{gp}}.$ We also identify the étale sites of $X$ and $X'$ via the universal homeomorphism $F_{X/S}.$ Recall the $\upmu$-indexed algebras $\calA_{\widetilde{X}}$ and $\calB_{X/S}$ and the sections $e_s$ of $\calA_{\widetilde{X}}$ (\ref{example314}).
\end{parag}

\begin{proposition}\label{logflatlift}
The lifting $\widetilde{F}:\widetilde{X}\ra \widetilde{X}'$ is log flat (\cite{Ogus2018} IV 4.1.1).
\end{proposition}

\begin{proof}
Recall that the exact relative Frobenius $F_{X/S}:X\ra X'$ is log flat by (\cite{SF1} 3.5). By (\cite{Ogus2018} 4.1.2), $\widetilde{X}$ and $\widetilde{X}'$ are log smooth, hence log flat, over $\widetilde{S}$. Let $\ov{x}$ be a geometric point of $\widetilde{X},$ $\ov{s}$ its image in $\widetilde{S}$ and $\kappa(\ov{s})$ the residue field at $\ov{s}.$ Since $p$ is nilpotent on $\widetilde{S},$ the fiber $\op{Spec}\kappa(\ov{s}) \ra \widetilde{S}$ factors through $S,$ where we equip the fiber with the logarithmic structure pullback of that of $\widetilde{S}.$ The base change $\widetilde{F}_{\ov{s}}$ of $\widetilde{F}$ along $\op{Spec} \kappa(\ov{s}) \ra \widetilde{S}$ is then a base change of $F_{X/S}.$ We deduce that $\widetilde{F}_{\ov{s}}$ is log flat. We conclude the log flatness of $\widetilde{F}$ by (\cite{Ogus2018} IV 4.2.2). 
\end{proof}

\begin{definition}\label{defBF}
We denote by $\calB_{\widetilde{F}},$ or $\calB_{\widetilde{X}/\widetilde{S}}$ if no confusion arises, the $\upmu$-indexed $\Ox_{\widetilde{X}',\upmu}$-subalgebra of $\calA_{\widetilde{X}}$ (\ref{example314}) generated locally by sections of the form $\alpha_{\widetilde{X}}(m)e_{-m}$ and $e_{\widetilde{F}^{\flat}m'}$ for local sections $m$ and $m'$ of $\calM_{\widetilde{X}}$ and $\calM_{\widetilde{X}'}^{\mathrm{gp}}$ respectively.
\end{definition}

\begin{lemma}\label{lem243}
Let $n$ be a positive integer and $g:M\ra N$ a morphism of $\Z/p^n\Z$-modules such that $N$ is flat over $\Z/p^n\Z.$ If the reduction $g_1$ of $g$ modulo $p$ is an isomorphism, then so is $g.$
\end{lemma}

\begin{proof}
If $y\in N$ then there exist $x\in M$ and $a\in N$ such that $y=g(x)+pa.$ It follows that $\op{coker}g= p\op{coker}g= \hdots = p^n\op{coker}g=0.$ This proves the surjectivity. Now let $K=\op{ker}g.$ The exact sequence of $\Z/p^n\Z$-modules
$$
0 \ra K\ra M\ra N \ra 0
$$
yields the exact sequence
$$
\hdots \ra \op{Tor}^1_{\Z/p^n\Z}\left (N,\mathbb{F}_p \right ) \ra K/pK \ra M/pM \ra N/pN \ra 0.
$$
The flatness of $N$ over $\Z/p^n\Z$ implies that $\op{Tor}^1_{\Z/p^n\Z} \left ( N,\mathbb{F}_p \right )=0.$ The injectivity of $g_1$ then implies that $K/pK=0$ so
$
K\subset pK \subset \hdots \subset p^nK=0.
$
\end{proof}

\begin{proposition}\label{propflatBF}
The $\Ox_{\widetilde{X}',\upmu}$-module $\calB_{\widetilde{F}}$ is flat over $\Z/p^n\Z$ and $\calB_{\widetilde{F}} \cap p\calA_{\widetilde{X}}=p\calB_{\widetilde{F}}.$
\end{proposition}

\begin{proof}
For a local section $s$ of $\upmu=\ov{\calM}_{\widetilde{X}}^{\mathrm{gp}},$ set $\calP=\calM_{\widetilde{X}'}^{\mathrm{gp}}\oplus \ov{\calM}_{\widetilde{X}}.$ Consider the subsheaf $\calP_s$ of $\calP$ consisting of local sections $(m',\ov{m})$ such that $\widetilde{F}^{\flat}m'-m=s$ in $\ov{\calM}_{\widetilde{X}}^{\mathrm{gp}}$ and
$$
f_s:\calD_s=\Ox_{\widetilde{X}'} \otimes_{\Z [\calM_{\widetilde{X}'}]} \Z[\calP_s] \ra \calA_{\widetilde{X},s}, \quad
a\otimes t^{(m',\ov{m})}  \mapsto  a\alpha_{\widetilde{X}}(m)e_{\widetilde{F}^{\flat}m'-m}.
$$
The image of $f_s$ is $\calB_{\widetilde{F},s}.$ By \ref{thm326}, $f_s$ induces an isomorphism $\calD_s/p\calD_s \isom \calB_{X/S,s}.$
Suppose we have local coordinates $m_1,\hdots,m_d\in \Gamma \left ( \widetilde{X},\calM_{\widetilde{X}} \right )$ of $\widetilde{X}$ with respect to $\widetilde{S}.$ For $I\in \N^d,$ set $m_I=\sum_{i=1}^dI_im_i$ and $e_I=e_{m_I}$ the local section of $\calA_{\widetilde{X}}$ given in \ref{example314}.
For all local sections $s$ of $\upmu,$ the morphism of $\Ox_{\widetilde{X}'}$-modules
$$
\psi:\bigoplus_{I\in \llbracket 0,p-1\rrbracket ^d} \calD_{s-m_I} \ra \calA_{\widetilde{X},s},\quad (x_I) \mapsto \sum_I f(x_I)e_{m_I}
$$
is, modulo $p,$ an isomorphism by (\cite{Lor2000} II 2.5). It follows, by \ref{lem243}, that $\psi$ itself is an isomorphism. Looking at the component corresponding to $I=0,$ we get that $\calB_{\widetilde{F},s}$ is a direct summand of $\calA_{\widetilde{X},s}.$ The result follows.
\end{proof}

\begin{remark}
By \ref{prop02}, $\calB_{\widetilde{F}}$ modulo $p$ is simply $\calB_{X/S}.$ Indeed, the canonical immersion $X\ra \widetilde{X}$ is strict so, by \ref{strictA}, we have $\calA_{\widetilde{X}}/p\calA_{\widetilde{X}}=\calA_X.$ The composition of the embedding $\calB_{\widetilde{F}} \hookrightarrow \calA_{\widetilde{X}}$ with the projection $\calA_{\widetilde{X}} \ra \calA_X$ induces, by \ref{prop02}, a surjection
$
\calB_{\widetilde{F}} \ra \calB_{X/S}.
$
Its kernel is $\calB_{\widetilde{F}} \cap p\calA_{\widetilde{X}},$ which is equal to $p\calB_{\widetilde{F}}$ by \ref{propflatBF}.
\end{remark}

\begin{proposition}\label{propABlocfree}
$\calA_{\widetilde{X}}$ is a locally free $\upmu$-indexed $\calB_{\widetilde{F}}$-module of rank $p^d,$ where $d$ is the rank of $\omega^1_{\widetilde{X}/\widetilde{S}}.$
\end{proposition}

\begin{proof}
The assertion being local, we may assume that we have local coordinates $m_1,\hdots,m_d\in \Gamma \left ( \widetilde{X},\calM_{\widetilde{X}} \right ).$ For $I\in \N^d,$ set $m_I=\sum_{i=1}^dI_im_i$ and $e_I=e_{m_I}$ the local section of $\calA_{\widetilde{X}}$ given in \ref{example314}. Consider then the morphism
$$
\varphi:\bigoplus_{I\in \llbracket 0,p-1 \rrbracket^d} \calB_{\widetilde{F}}\left ( m_I \right ) \ra \calA_{\widetilde{X}},\quad (b_I) \mapsto \sum_I b_Ie_I,
$$
where $\calB_{\widetilde{F}}\left ( m_I \right )$ is the $\upmu$-indexed algebra $\calB_{\widetilde{F}}$ whose indices are shifted by $-m_I.$
Modulo $p,$ $\varphi$ is an isomorphism by (\cite{Lor2000} II 2.5). Since $\calA_{\widetilde{X}}$ is locally free over $\Ox_{\widetilde{X},\upmu}$ which is, by hypothesis, flat over $\Z/p^n\Z,$ we conclude by \ref{lem243}.
\end{proof}

\begin{proposition}\label{prop347}
Let $s$ be a local section of $\ov{\calM}_{\widetilde{X}}^{\mathrm{gp}},$ $\calP=\calM_{\widetilde{X}'}^{\mathrm{gp}}\oplus \ov{\calM}_{\widetilde{X}}$ and $\calP_s$ the subsheaf of $\calP$ consisting of local sections $(m',\ov{m})$ such that $\widetilde{F}^{\flat}m'-m=s$ in $\ov{\calM}_{\widetilde{X}}^{\mathrm{gp}}.$ Then, the morphism of $\Ox_{\widetilde{X}'}$-modules
$$
\Ox_{\widetilde{X}'} \otimes_{\Z [\calM_{\widetilde{X}'}]} \Z[\calP_s] \ra \calB_{\widetilde{F},s}, \quad
a\otimes t^{(m',\ov{m})}  \mapsto  a\alpha_{\widetilde{X}}(m)e_{\widetilde{F}^{\flat}m'-m}
$$
is an isomorphism and is compatible with multiplication as $s$ varies.
\end{proposition}

\begin{proof}
This follows from \ref{thm326}, \ref{lem243} and \ref{propflatBF}.
\end{proof}

\begin{theorem}\label{thmkeyiso}
Let $\widetilde{Y}=\widetilde{X}\times_{\widetilde{X}'}\widetilde{X},$ $Y=X\times_{X'}X$ and identify the étale sites of $X,$ $Y$ and $\widetilde{Y}$ via the diagonal immersion $X\ra Y$ (which is a universal homeomorphism by \ref{prop114}) and the canonical nilpotent closed immersion $Y\ra \widetilde{Y}.$ The projections $\widetilde{p}_1,\widetilde{p}_2:\widetilde{Y} \ra \widetilde{X}$ induce an isomorphism of $\upmu$-indexed algebras
$$
\calA_{\widetilde{X}} \circledast_{\calB_{\widetilde{F}}} \calA_{\widetilde{X}} \isom \calA_{\widetilde{Y}}.
$$
\end{theorem}

\begin{proof}
We abusively denote by $\widetilde{p}_1,\widetilde{p}_2:\calA_{\widetilde{X}} \ra \calA_{\widetilde{Y}}$ the morphisms of $\upmu$-indexed algebras induces by the projections $\widetilde{p}_1$ and $\widetilde{p}_2.$ The morphisms $\widetilde{p}_1$ and $\widetilde{p}_2$ clearly agree on $\Ox_{\widetilde{X}'}$ and local sections of the form $e_{\widetilde{F}^{\flat}m'}$ where $m'$ is a local section of $\calM_{\widetilde{X}'}^{\mathrm{gp}}.$ Denote by $\calI$ the ideal of the closed diagonal immersion $\widetilde{\Delta}:\widetilde{X}\ra \widetilde{Y}.$ By (\cite{SF1} 4.3), we have an exact sequence of monoids
$$
0 \ra 1+\calI \xrightarrow{\alpha_{\widetilde{Y}}^{-1}} \calM_{\widetilde{Y}} \xrightarrow{\widetilde{\Delta}^{\flat}} \calM_{\widetilde{X}} \ra 0.
$$
Let $m$ be a local section of $\calM_{\widetilde{X}}.$ There exists a local section $a$ of $1+\calI \subset \Ox_{\widetilde{Y}}^{\times}$ such that
$$
\widetilde{p}_2^{\flat}(m)=\widetilde{p}_1^{\flat}(m)+\alpha_{\widetilde{Y}}^{-1}(a).
$$
Then, we have
\begin{alignat*}{2}
\widetilde{p}_2(\alpha_{\widetilde{X}}(m)e_{-m})&=\widetilde{p}_2^{\#}(\alpha_{\widetilde{X}}(m))e_{-\widetilde{p}_2^{\flat}(m)}=\alpha_{\widetilde{Y}}(\widetilde{p}_2^{\flat}(m))e_{-\widetilde{p}_2^{\flat}(m)} \\
&= a\alpha_{\widetilde{Y}}(\widetilde{p}_1^{\flat}(m))e_{-\alpha_{\widetilde{Y}}^{-1}(a)}e_{-\widetilde{p}_1^{\flat}(m)}=ae_{-\alpha_{\widetilde{Y}}^{-1}(a)}\widetilde{p}_1(\alpha_{\widetilde{X}}(m)e_{-m})=\widetilde{p}_1(\alpha_{\widetilde{X}}(m)e_{-m}).
\end{alignat*}
This proves that the morphism is well-defined. To conclude that it is an isomorphism, it is sufficient to prove that $\calA_{\widetilde{Y}}$ is flat over $\Z/p^n\Z$ and to use \ref{lem243} and \ref{th234}. By \ref{prop114}, the projections $\widetilde{Y} \ra \widetilde{X}$ are strict. By \ref{logflatlift}, these projections are also log flat. They are hence flat. Since $\widetilde{X}$ is flat over $\Z/p^n\Z,$ we deduce that $\widetilde{Y}$ is flat over $\Z/p^n\Z$ hence so is $\calA_{\widetilde{Y}}.$
\end{proof}

\begin{corollaire}\label{corBequalizer}
The $\upmu$-indexed algebra $\calB_{\widetilde{F}}$ is the equalizer of the two morphisms $\calA_{\widetilde{X}} \ra \calA_{\widetilde{Y}}$ induced by the two projections $\widetilde{Y} \ra \widetilde{X}.$
\end{corollaire}

\begin{proof}
By \ref{propABlocfree}, $\calA_{\widetilde{X}}$ is locally free of finite rank over $\calB_{\widetilde{F}}.$ By \ref{example305}, $\calA_{\widetilde{X}}$ is faithfully flat over $\calB_{\widetilde{F}}.$ By \ref{indfflem2}, we have the exact sequence
$$
0 \ra \calB_{\widetilde{F}} \ra \calA_{\widetilde{X}} \ra \calA_{\widetilde{X}} \circledast_{\calB_{\widetilde{F}}} \calA_{\widetilde{X}}.
$$
We conclude by \ref{thmkeyiso}.
\end{proof}

\begin{proposition}\label{propnablaF}
There exists a unique connection $\nabla_{\widetilde{F}}:\calB_{\widetilde{F}} \ra \calB_{\widetilde{F}} \otimes_{\Ox_{\widetilde{X}',\upmu}} \omega^1_{\widetilde{X}'/\widetilde{S},\upmu}$ satisfying:
\begin{enumerate}
\item $\nabla_{\widetilde{F}}(\alpha_{\widetilde{X}}(m)e_{-m})=0$ for all local sections $m$ of $\calM_{\widetilde{X}}.$
\item $\nabla_{\widetilde{F}}(e_{\widetilde{F}^{\flat}m'})=e_{\widetilde{F}^{\flat}m'}\otimes \op{dlog}m'$ for all local sections $m'$ of $\calM_{\widetilde{X}'}^{\mathrm{gp}}.$
\item $\nabla_{\widetilde{F}}(a)=1\otimes da$ for all local sections $a$ of $\Ox_{\widetilde{X}'}\subset \calB_{\widetilde{F},0}.$
\item $\nabla_{\widetilde{F}}(ab)= a\nabla_{\widetilde{F}}(b)+b\nabla_{\widetilde{F}}(a)$ for all local sections $a$ and $b$ of $\calB_{\widetilde{F}}.$
\end{enumerate}
In addition, $\nabla_{\widetilde{F}}$ is integrable.
\end{proposition}

\begin{proof}
Since $\calB_{\widetilde{F}}$ is, by definition, locally generated by $\Ox_{\widetilde{X}'}$ and local sections of the form $\alpha_{\widetilde{X}}(m)e_{-m},$ $e_{\widetilde{F}^{\flat}m'}$ for local sections $m$ and $m'$ of $\calM_{\widetilde{X}}$ and $\calM_{\widetilde{X}'}^{\mathrm{gp}}$ respectively, a connection satisfying the fourth property is entirely determined by the image of these local sections. It is hence sufficient to prove the existence of $\nabla_{\widetilde{F}}$ locally. Let $s$ be a local section of $\upmu.$ By \ref{prop347}, we may identify $\calB_{\widetilde{F},s}$ with $\Ox_{\widetilde{X}'} \otimes_{\Z [\calM_{\widetilde{X}'}]} \Z[\calP_s],$ where $\calP_s$ is given in \ref{prop347}. We set
$$
\nabla_{\widetilde{F},s}:\begin{array}[t]{clc}
\Ox_{\widetilde{X}'} \otimes_{\Z [\calM_{\widetilde{X}'}]} \Z [\calP_s] & \ra & \left ( \Ox_{\widetilde{X}'} \otimes_{\Z [\calM_{\widetilde{X}'}]} \Z [\calP_s] \right ) \otimes_{\Ox_{\widetilde{X}'}} \omega^1_{\widetilde{X}'/\widetilde{S}} \\
a \otimes t^{(m',\ov{m})} & \mapsto & 1\otimes t^{(m',\ov{m})} \otimes da + a\otimes t^{(m',\ov{m})} \otimes \op{dlog}m',
\end{array}
$$
where $m'$ and $\ov{m}$ are local sections of $\calM_{\widetilde{X}'}^{\mathrm{gp}}$ and $\ov{\calM}_{\widetilde{X}}$ respectively.
One easily checks that this is well-defined and that it is a connection satisfying the required properties.
\end{proof}

\begin{definition}\label{defthetaBF}
We denote by $\theta_{\widetilde{F}}:\calB_{\widetilde{F}} \ra \calB_{\widetilde{F}} \otimes_{\Ox_{\widetilde{X}',\upmu}} \omega^1_{\widetilde{X}'/\widetilde{S},\upmu}$ the integrable $p$-connection defined by $\theta_{\widetilde{F}}=p\nabla_{\widetilde{F}},$ where $\nabla_{\widetilde{F}}$ is defined in \ref{propnablaF}.
\end{definition}

\begin{remark}
The $p$-connection $\theta_{\widetilde{F}},$ defined in \ref{defthetaBF}, is quasi-nilpotent by \ref{remquasinil}.
\end{remark}

\section{Indexed faithfully flat descent}

\begin{parag}
In this section, we consider a ringed topos $(X,\Ox)$ on a site $C$ and a sheaf of monoids $\upmu$ in $X.$ We suppose that $X$ has enough points. We equip the localized topos $X_{/\upmu}$ with the ring $\Ox_{\upmu},$ localization of $\Ox$ at $\upmu.$ For a morphism $f:\calA \ra \calA'$ of $\upmu$-indexed algebras and a $\upmu$-indexed $\calA$-module $\calE,$ we denote by $f^*\calE$ the $\upmu$-indexed $\calA'$-module $\calA'\circledast_{\calA}\calE.$ For a $\upmu$-indexed algebra $\calA,$ we denote by $\boldsymbol{\op{Mod}}_{\calA}$ the category of $\upmu$-indexed $\calA$-modules and by $\boldsymbol{\op{Mod}}_{\Ox_{\upmu}}$ the category of $\Ox_{\upmu}$-modules.
\end{parag}

\subsection{Faithfully flat indexed algebras}

\begin{proposition}
Let $\calA\ra \calA'$ be a morphism of $\upmu$-indexed algebras. Then the functor
$$
\boldsymbol{\op{Mod}}_{\calA} \ra \boldsymbol{\op{Mod}}_{\calA'}, \quad
\calE  \mapsto  \calE\circledast_{\calA}\calA'
$$
is right exact.
\end{proposition}

\begin{proof}
Let $\calE' \ra \calE \ra \calE'' \ra 0$
be an exact sequence of $\upmu$-indexed $\calA$-modules. Let $p_1,p_2:\upmu^2\ra \upmu$ be the canonical projections and $\sigma:\upmu^2\ra \upmu$ and $\tau:\upmu^3\ra \upmu$ the addition maps. We abusively denote by
$$
p_1,p_2,\sigma:X_{/\upmu^2} \ra X_{/\upmu},\quad \tau : X_{/\upmu^3} \ra X_{/\upmu}
$$
the corresponding morphisms of ringed topoi. Then the sequence
$p_1^*\calE' \ra p_1^*\calE \ra p_1^*\calE'' \ra 0$
is exact. The functor
$
\boldsymbol{\op{Mod}}_{\Ox_{\upmu^2}}  \ra  \boldsymbol{\op{Mod}}_{\Ox_{\upmu^2}},\quad
\calM  \mapsto  \calM\otimes_{\Ox_{\upmu^2}}p_2^*\calA'
$
is right exact so the sequence
$$
\calE'\boxtimes \calA' \ra \calE\boxtimes \calA' \ra \calE''\boxtimes \calA' \ra 0
$$
is exact. Since $\sigma_!$ is a left adjoint by definition, it is right exact and so the sequence
\begin{equation}\label{exseq2}
\sigma_!(\calE'\boxtimes \calA') \ra \sigma_!(\calE\boxtimes \calA') \ra \sigma_!(\calE''\boxtimes \calA') \ra 0
\end{equation}
is exact. Similarly, we have an exact sequence
\begin{equation}\label{exseq1}
\tau_!(\calA\boxtimes\calE'\boxtimes \calA') \ra \tau_!(\calA\boxtimes\calE\boxtimes \calA') \ra \tau_!(\calA\boxtimes\calE''\boxtimes \calA') \ra 0.
\end{equation}
The sequences \eqref{exseq2} and \eqref{exseq1} fit into the following commutative diagram
$$
\begin{tikzcd}
\tau_!(\calA\boxtimes\calE'\boxtimes \calA') \ar{r}{u'} \ar{d}{\alpha'} & \tau_!(\calA\boxtimes\calE\boxtimes \calA') \ar{r}{v'} \ar{d}{\alpha} & \tau_!(\calA\boxtimes\calE''\boxtimes \calA') \ar{r} \ar{d}{\alpha''} & 0 \\
\sigma_!(\calE'\boxtimes \calA') \ar{r}{u} \ar{d}{\beta'} & \sigma_!(\calE\boxtimes \calA') \ar{r}{v} \ar{d}{\beta} & \sigma_!(\calE''\boxtimes \calA') \ar{r} \ar{d}{\beta''} & 0 \\
\calE'\circledast_{\calA}\calA' \ar{r}{u''} \ar{d} & \calE\circledast_{\calA}\calA' \ar{r}{v''} \ar{d} & \calE''\circledast_{\calA}\calA' \ar{r} \ar{d} & 0 \\
0 & 0 & 0 &
\end{tikzcd}
$$
where the vertical exact sequences are given in \eqref{exseq4182}. Since the upper and middle horizontal sequences and the vertical sequences are exact, the lower horizontal sequence is exact. Indeed, the morphism $v''$ is clearly surjective. Let us prove that $\op{Im}u''=\op{Ker}v''.$ Since $X$ has enough points, by (\cite{SGA43} Exposé III 6.7.4), the topos $X_{/\upmu}$ has also enough points. Let $P$ be a point of $X_{/\upmu},$ $x\in (\op{Ker}v'')_P$ and $y\in \sigma_!(\calE \boxtimes \calA')_P$ such that
$\beta(y)=x.$
Since
$\beta''(v(y))=v''(\beta(y))=v''(x)=0,$
there exists $x''\in \tau_!(\calA \boxtimes \calE'' \boxtimes \calA')_P$ such that
$\alpha''(x'')=v(y).$
Let $t\in \tau_!(\calA \boxtimes \calE \boxtimes \calA')_P$ such that
$v'(t)=x''.$
Then $$v\circ \alpha(t)=\alpha''\circ v'(t)=\alpha''(x'')=v(y)$$
so $y-\alpha(t) \in (\op{Ker}v)_P$ and there exists $y'\in \sigma_!(\calE' \boxtimes \calA')_P$ such that
$u(y')=y-\alpha(t).$
Set $x'=\beta'(y').$ Then
$$u''(x')=u''\circ \beta'(y')=\beta\circ u(y')=\beta(y)-\beta\circ \alpha(t)=\beta(y)=x.$$ 
\end{proof}

\begin{proposition}\label{indprop5}
Let $\calA \ra \calA'$ be a morphism of $\upmu$-indexed algebras. Then the following assertions are equivalent:
\begin{enumerate}
\item $\calA \ra \calA'$ is flat and for all $\upmu$-indexed $\calA$-modules $\calE,$ 
$\calE \circledast_{\calA}\calA'=0$ implies $\calE=0.$
\item Any sequence of $\upmu$-indexed $\calA$-modules
$0 \ra \calE' \ra \calE \ra \calE'' \ra 0$
is exact if and only if
$$0 \ra \calE'\circledast_{\calA}\calA' \ra \calE\circledast_{\calA}\calA' \ra \calE''\circledast_{\calA}\calA'\ra 0$$
is exact.
\end{enumerate}
\end{proposition}

\begin{proof}
Suppose the first assertion is satisfied. Consider a sequence of $\upmu$-indexed $\calA$-modules $(E)$ and its base change $(E')$:
$$
(E): 0\ra \calE' \xrightarrow{u} \calE \xrightarrow{v} \calE'' \ra 0,\quad (E'):0\ra \calE'\circledast_{\calA}\calA' \xrightarrow{u'} \calE \circledast_{\calA}\calA' \xrightarrow{v'} \calE'' \circledast_{\calA}\calA'\ra 0.
$$
If $(E)$ is exact then $(E')$ is exact by the flatness of $\calA \ra \calA'.$
Suppose $(E')$ is exact. Since the functor $-\circledast_{\calA}\calA'$ is exact, we have $\calA'\circledast_{\calA}\op{Im}(v\circ u)=\op{Im}(v'\circ u')=0.$ By the first assertion, we get $v\circ u=0.$ We also have
$$\op{Ker}(\calE' \ra \calE)\circledast_{\calA}\calA'=\op{Ker}(\calE'\circledast_{\calA}\calA'\ra \calE\circledast_{\calA}\calA'),\quad \left (\op{Ker}v/\op{Im}u\right )\circledast_{\calA}\calA'=\op{Ker}(v')/\op{Im}(u'),$$
$$\op{coker}(\calE\circledast_{\calA}\calA'\ra \calE''\circledast_{\calA}\calA')=\op{coker}(\calE\ra \calE'')\circledast_{\calA}\calA'.$$
Since $(E')$ is exact, we have
$$\op{Ker}(\calE'\circledast_{\calA}\calA' \ra \calE\circledast_{\calA}\calA') = \op{Ker}(v')/\op{Im}(u') = \op{coker}(\calE\circledast_{\calA}\calA'\ra \calE''\circledast_{\calA}\calA') = 0.$$
So
$\op{Ker}(\calE' \ra \calE)\circledast_{\calA}\calA'=\left (\op{Ker}v/\op{Im}u\right )\circledast_{\calA}\calA'=\op{coker}(\calE\ra \calE'')\circledast_{\calA}\calA'=0.$
It follows, by the first assertion, that
$$\op{Ker}(\calE'\ra \calE) = \op{Ker}v/\op{Im}u = \op{coker}(\calE \ra \calE'') = 0.$$
Hence the exactness of $(E).$

Now suppose that the second assertion is satisfied. The flatness of $\calA \ra \calA'$ follows immediately. Let $\calE$ be a $\upmu$-indexed $\calA$-module such that
$\calE\circledast_{\calA}\calA'=0.$
The zero morphism
$\calE\circledast_{\calA}\calA'\xrightarrow{0} \calE\circledast_{\calA}\calA'$
is injective so by the second assertion,
$\calE \xrightarrow{0} \calE$
is injective. It follows that $\calE=0.$
\end{proof}

\begin{definition}\label{indff}
\begin{enumerate}
\item A morphism $\calA \ra \calA'$ of $\upmu$-indexed algebras is said to be \emph{flat} if the functor
$$
\boldsymbol{\op{Mod}}_{\calA}  \ra  \boldsymbol{\op{Mod}}_{\calA'}, \quad
\calE  \mapsto  \calE \circledast_{\calA}\calA'
$$
is exact. In this case, we say that $\calA'$ is a \emph{flat $\upmu$-indexed $\calA$-algebra.}
\item A morphism $\calA \ra \calA'$ of $\upmu$-indexed algebras is said to be \emph{faithfully flat} if it satisfies one of the two equivalent properties of \ref{indprop5}.
In this case, we say that $\calA'$ is a \emph{faithfully flat $\upmu$-indexed $\calA$-algebra.}
\end{enumerate}
\end{definition}

\begin{exmp}\label{example305}
\begin{enumerate}
\item If $\upmu$ is a sheaf of abelian groups, then locally free indexed algebras are faithfully flat. Indeed, let $\calA \ra \calA'$ be a morphism of $\upmu$-indexed algebras such that $\calA'$ is locally free over $\calA.$ Since exactness of sequences of $\upmu$-indexed modules is checked locally, we may suppose that there exist sections $e_i,\ i\in I$ of $\upmu$ such that $\calA' = \bigoplus_{i\in I}\calA(e_i).$ The result then follows from the fact that, for $i\in I$ and for a $\upmu$-indexed $\calA$-module $\calM,$ we have a canonical isomorphism $\calM \circledast_{\calA}\calA(e_i) \xrightarrow{\raisebox{-0.4ex}[0ex][0ex]{$\sim$}}  \calM(e_i).$

\item Let $\calA$ be a $\upmu$-indexed algebra and $R \ra S$ a faithfully flat morphism of $\Ox$-algebras. Then $\calA \otimes_{\Ox_{\upmu}} R_{\upmu} \ra \calA \otimes_{\Ox_{\upmu}} S_{\upmu}$ is a faithfully flat morphism of $\upmu$-indexed algebras. Indeed, by (\cite{SF3} 3.26), both $\calA \otimes_{\Ox_{\upmu}} R_{\upmu}$ and $\calA \otimes_{\Ox_{\upmu}} S_{\upmu}$ have natural structures of $\upmu$-indexed algebras and, by (\cite{SF3} 3.27), for a $\upmu$-indexed $\calA \otimes_{\Ox_{\upmu}} R_{\upmu}$-module $\calM,$ $\calM$ has a natural $R_{\upmu}$-module structure and we have a canonical isomorphism
$$
\calM \circledast_{\calA \otimes_{\Ox_{\upmu}} R_{\upmu}} \left (\calA \otimes_{\Ox_{\upmu}} S_{\upmu} \right ) \isom \calM \otimes_{R_{\upmu}} S_{\upmu}.
$$
\end{enumerate}
\end{exmp}

\subsection{Descent}

\begin{parag}
Let $\calA\ra \calA'$ be a morphism of $\upmu$-indexed $\Ox$-algebras and $\calE$ a $\upmu$-indexed $\calA$-module. For all positive integers $n,$ we denote by $\calA'^{\circledast n}$ the product $\calA' \circledast_{\calA}\hdots \circledast_{\calA}\calA'$ of $n$ copies of $\calA'$ and we consider the morphism
$$d^n: \calE\circledast_{\calA}\calA'^{\circledast n} \ra \calE\circledast_{\calA}\calA'^{\circledast (n+1)}$$
induced by the morphisms of presheaves defined, for all sections $k$ of $\upmu$ over an object $U$ of $C,$ by 
$$\begin{array}[t]{clc}
\bigoplus\limits_{\substack{i_0,\hdots,i_n\in \upmu(U)\\ i_0+\hdots+i_n=k}} \left (\calE_{i_0}\otimes \calA'_{i_1}\otimes \hdots \otimes \calA'_{i_n} \right ) (U) & \ra  & \bigoplus\limits_{\substack{i'_0,\hdots,i'_{n+1}\in \upmu(U) \\ i'_0+\hdots+i'_{n+1}=k}} \left (\calE_{i'_0}\otimes \calA'_{i'_1}\otimes \hdots \otimes \calA'_{i'_{n+1}} \right )(U) \\
x_0 \otimes a_1 \otimes \hdots \otimes a_n & \mapsto & \sum\limits_{j=0}^n(-1)^jx_0 \otimes a_1 \otimes \hdots a_j\otimes 1\otimes a_{j+1} \otimes \hdots \otimes a_n
 \end{array}$$
The morphisms $d^n$ fit into a complex of $\upmu$-indexed $\calA$-modules
\begin{equation}\label{indcomplex}
0\ra \calE \xrightarrow{d^0} \calE \circledast_{\calA} \calA' \xrightarrow{d^1} \calE \circledast_{\calA} \calA' \circledast_{\calA} \calA' \xrightarrow{d^2} \hdots
\end{equation}
The following lemma provides a sufficient condition for this complex to be exact.
\end{parag}

\begin{lemma}\label{indfflem1}
Let $\calA\ra \calA'$ be a morphism of $\upmu$-indexed $\Ox$-algebras and $\calE$ a $\upmu$-indexed $\calA$-module. If the canonical morphism
$$d^0:\calE \ra \calE \circledast_{\calA}\calA'$$
admits a section $s:\calE \circledast_{\calA} \calA' \ra \calE$ as $\upmu$-indexed $\calA$-modules, then the complex \eqref{indcomplex}
is exact.
\end{lemma}

\begin{proof}
For all positive integers $n,$ let $$s^n=s\circledast \op{Id}_{\calA'^{\circledast (n-1)}}:\calE \circledast_{\calA}\calA'^{\circledast n} \ra \calE \circledast_{\calA}\calA'^{\circledast (n-1)}.$$
Then the diagram
$$
\begin{tikzcd}
\hdots \ar{r} & \calE\circledast_{\calA}\calA'^{\circledast (n-1)} \ar{r}{d^{n-1}} & \calE\circledast_{\calA}\calA'^{\circledast n}\ar{r}{d^n}\ar{dl}{s^n} & \calE\circledast_{\calA}\calA'^{\circledast (n+1)} \ar{r}\ar{dl}{s^{n+1}} & \hdots \\
\hdots \ar{r} & \calE\circledast_{\calA}\calA'^{\circledast (n-1)}\ar[swap]{r}{d^{n-1}} & \calE\circledast_{\calA}\calA'^{\circledast n} \ar{r} & \calE\circledast_{\calA}\calA'^{\circledast (n+1)}\ar{r} & \hdots 
\end{tikzcd}
$$
satisfies
$$s^{n+1} \circ d^n+d^{n-1}\circ s^n=\op{Id}_{\calE\circledast_{\calA}\calA'^{\circledast n}}.$$
Indeed, let $k,i_0,\hdots,i_n\in \upmu(U)$ be sections over an object $U$ of $C$ such that $i_0+\hdots+i_n=k$ and $x\in \Gamma(U,\calE_{i_0}),$ $a_1\in \Gamma(U,\calA'_{i_1}),\hdots,a_n\in \Gamma(U,\calA'_{i_n}).$
Denote
\begin{alignat*}{2}
y=x\otimes a_1 \otimes \hdots \otimes a_n &\in \calE_{i_0}(U) \otimes \calA'_{i_1}(U) \otimes \hdots \otimes \calA'_{i_n}(U)\\
&\subset \bigoplus\limits_{j_0+\hdots+j_n=k}\calE_{j_0}(U) \otimes \calA'_{j_1}(U) \otimes \hdots \otimes \calA'_{j_n}(U).
\end{alignat*}
Then
\begin{alignat*}{2}
s^{n+1} \circ d^n(y) &= \sum_{j=0}^n(-1)^js^{n+1}(x\otimes a_1 \otimes \hdots \otimes  a_j\otimes 1\otimes a_{j+1} \otimes \hdots \otimes a_n) \\
&= \sum_{j=0}^n(-1)^js(x\otimes a_1) \otimes \hdots \otimes  a_j\otimes 1\otimes a_{j+1} \otimes \hdots \otimes a_n \\
&= y+\sum_{j=1}^n(-1)^js(x\otimes a_1) \otimes \hdots \otimes  a_j\otimes 1\otimes a_{j+1} \otimes \hdots \otimes a_n = y-d^{n-1}\circ s^n(y).
\end{alignat*}
It follows that the identity map of the complex \eqref{indcomplex} is homotopic to the zero map. The proof is thus complete.
\end{proof}

\begin{proposition}\label{indfflem2}
Let $\calA \ra \calA'$ a faithfully flat morphism of $\upmu$-indexed algebras and $\calE$ a $\upmu$-indexed $\calA$-module. Then the complex \eqref{indcomplex} is exact.
\end{proposition}

\begin{proof}
By the definition of faithful flatness (\ref{indff}), it is sufficient to prove that the complex
$$0\ra \calE\circledast_{\calA} \calA' \rightarrow \calE \circledast_{\calA} \calA'\circledast_{\calA} \calA' \rightarrow \calE \circledast_{\calA} \calA' \circledast_{\calA} \calA'\circledast_{\calA} \calA' \rightarrow \hdots$$
is exact. By \ref{indfflem1}, it is sufficient to prove that the canonical morphism
$$\calE\circledast_{\calA} \calA' \ra \calE\circledast_{\calA} \calA'\circledast_{\calA} \calA'$$
admits a section $s$ as a map of $\upmu$-indexed $\calA$-modules. Indeed, we can take $s$ to be equal to $\op{Id}_{\calE}\circledast v,$ where $v$ is the multiplication $\calA' \circledast_{\calA} \calA' \ra \calA'$ (\cite{SF3} 3.22.22).
\end{proof}

\begin{theorem}[Descent of indexed morphisms]\label{inddescentmor}
Let $\calA \ra \calA'$ be a faithfully flat morphism of $\upmu$-indexed algebras and $\calE_1$ and $\calE_2$ two $\upmu$-indexed $\calA$-modules. Set $\calA''=\calA'\circledast_{\calA}\calA'$ and, for $i=1,2,$ $\calE_i'=\calE_i\circledast_{\calA}\calA'$ and $\calE_i''=\calE'_i\circledast_{\calA}\calA'$ and denote by $\pi_1,\pi_2:\calA' \ra \calA'\circledast_{\calA} \calA'$ the morphisms $x\mapsto x\otimes 1$ and $x\mapsto 1\otimes x$ (\cite{SF3} 3.22.24). Then, the sequence
$$0\ra \op{Hom}_{\calA}(\calE_1,\calE_2) \ra \op{Hom}_{\calA'}(\calE_1',\calE_2') \ra \op{Hom}_{\calA''}(\calE_1'',\calE_2'')$$
is exact, where the first arrow is $u\mapsto u\circledast \op{Id}_{\calA'}$ and the second arrow is $u'\mapsto u'\circledast \pi_2 -u'\circledast \pi_1.$
\end{theorem}

\begin{proof}
The injectivity of the first arrow $\op{Hom}_{\calA}(\calE_1,\calE_2) \ra \op{Hom}_{\calA'}(\calE_1',\calE_2')$ results from the injectivity of $\calE \ra \calE \circledast_{\calA} \calA'$ (\ref{indfflem2}). Now let $u':\calE'_1\ra \calE'_2$ be a morphism of $\upmu$-indexed $\calA'$-modules such that $u'\circledast \pi_2 -u'\circledast \pi_1=0.$ Since the canonical morphism $\calE_i \ra \calE'_i$ is injective, we can view $\calE_i$ as a $\upmu$-indexed $\calA$-submodule of $\calE'_i$ and it is thus sufficient to prove that $$u'(\calE_1)\subset \calE_2.$$
Consider the morphism
$v_1:\calE_1 \ra \calE''_1=\calE_1 \circledast_{\calA} \calA'\circledast_{\calA}\calA'$
induced by the morphisms
$$\begin{array}[t]{clclc}
\calE_{1,l}(U) & \ra & \calE_{1,l}(U)\otimes \calA'_0(U)\otimes \calA'_0(U) & \hookrightarrow & \bigoplus_{i+j+k=l}\calE_{1,i}(U)\otimes\calA'_j(U)\otimes\calA'_k(U) \\
x & \mapsto & x\otimes 1\otimes 1 & & 
\end{array}$$
for every section $l\in\upmu(U)$ over an object $U$ of $C.$
The diagram
$$
\begin{tikzcd}
\calE_1 \ar{d}{d^0} \ar{rr}{v_1} & & \calE''_1\ar{d}{u'\circledast\pi_2 -u'\circledast \pi_1} \\
\calE'_1 \ar{r}{u'} & \calE'_2 \ar{r}{d^1} & \calE''_2
\end{tikzcd}
$$
is commutative and so
$u'(\calE_1)\subset \op{ker}d^1.$ But $\op{ker}d^1=\calE_2$ by \ref{indfflem2}. This concludes the proof.
\end{proof}

\begin{theorem}[Descent of indexed modules]\label{inddescentmod}
Let $f:\calA\ra \calA'$ be a faithfully flat morphism of $\upmu$-indexed algebras. Denote by
$$\pi_1,\pi_2:\calA' \ra \calA' \circledast_{\calA} \calA',\ p_{ij}:\calA'\circledast_{\calA} \calA' \ra \calA'\circledast_{\calA} \calA'\circledast_{\calA} \calA',\ 1\le i<j\le 3$$
the projection morphisms (\cite{SF3} 3.22.24). Let $\calE'$ be a $\upmu$-indexed $\calA'$-module. We suppose given a $\calA'\circledast_{\calA} \calA'$-linear isomorphism
$$\epsilon:\pi_2^*\calE' \isom \pi_1^*\calE'$$
satisfying the cocycle condition
$p_{13}^*\epsilon=p_{12}^*\epsilon \circ p_{23}^*\epsilon.$
Then, there exists a $\upmu$-indexed $\calA$-module $\calE,$ unique up to a unique isomorphism, and an isomorphism
$$\varphi:f^*\calE \isom \calE',$$
making the following diagram commutative
\begin{equation}\label{diag525}
\begin{tikzcd}
\pi_2^*f^*\calE \ar{r}{\pi_2^*\varphi} \ar[swap]{d}{c} & \pi_2^*\calE' \ar{d}{\epsilon} \\
\pi_1^*f^*\calE \ar{r}{\pi_1^*\varphi} & \pi_1^*\calE',
\end{tikzcd}
\end{equation}
where the left vertical arrow $c$ is the canonical isomorphism.
\end{theorem}

\begin{proof}
We start by proving the existence. Let
\begin{equation}\label{eq5251}
\calE=\op{Ker} \left ( \op{Id}_{\calE'}\circledast \pi_1-\epsilon\circ (\op{Id}_{\calE'}\circledast \pi_2) \right ).
\end{equation}
Then $\calE$ is a $\upmu$-indexed $\calA$-module as the kernel of an $\calA$-linear morphism.
Let $g=\op{Id}_{\calE'}\circledast \pi_1-\epsilon\circ (\op{Id}_{\calE'}\circledast \pi_2)$ and $p_i:\calA'\ra \calA' \circledast_{\calA} \calA'\circledast_{\calA} \calA', 1\le i\le 3$ the projection morphisms defined in (\cite{SF3} 3.22.24). The exact sequence
$$0\ra \calE \hookrightarrow \calE' \xrightarrow{g} \pi_1^*\calE'$$
yields the exact sequence
$$0\ra \calA'\circledast_{\calA}\calE \hookrightarrow \pi_2^*\calE' \xrightarrow{g} p_2^*\calE'.$$
Consider the diagram
$$
\begin{tikzcd}
0\ar{r} & \calA' \circledast_{\calA} \calE\ar[dashed]{d}{\varphi}  \ar{r} & \pi_2^*\calE' \ar{r}{g}\ar{d}{\epsilon} & p_2^*\calE' \ar{d}{-p_{12}^*\epsilon} \\
0 \ar{r} & \calE' \ar{r}{d^0} & \pi_1^*\calE' \ar{r}{d^1} & p_1^*\calE'.
\end{tikzcd}
$$
The lower sequence is exact by \ref{indfflem2}. The right square is commutative by the cocycle condition of $\epsilon.$ The existence of the dashed arrow $\varphi$ follows. It is an isomorphism since the two other vertical arrows are isomorphisms. For a local section $x$ of $\calE,$ we have $\varphi(1\otimes x)=x.$ Then, by definition \eqref{eq5251} of $\calE,$ we have
$$
\epsilon \circ \pi_2^*\varphi(1\otimes 1 \otimes x)=\epsilon(1\otimes x)=x\otimes 1=\pi_1^*\varphi(1\otimes 1\otimes x)=\pi_1^*\varphi \circ c(1\otimes 1 \otimes x).
$$
We deduce the commutativity of \eqref{diag525}. Suppose we have two such descent modules $(\calE_1,\varphi_1)$ and $(\calE_2,\varphi_2).$ Consider the isomorphism $u'=\varphi_2^{-1}\circ \varphi_1.$ By \eqref{diag525}, we have $\pi_1^*u'=c_2\circ \pi_2^*u' \circ c_1^{-1},$ where $c_i:\pi_2^*f^*\calE_i \isom \pi_1^*f^*\calE_i$ are the canonical isomorphisms. By \ref{inddescentmor}, there exists a unique $u:\calE_1 \ra \calE_2$ such that $f^*u=u'.$ Furthermore, $u$ is an isomorphism since $u'$ is. 
\end{proof}

\section{Local Cartier transform}

\begin{parag}\label{indparag91}
In this section, we will provide an indexed version of the log Shiho functor (\cite{SF1} 9). Fix a log smooth morphism $f:\frakX\ra \frakS$ of logarithmic formal schemes, flat and locally of finite type over $W,$ and denote by $f_1:X\ra S$ its special fiber. We suppose that the exact relative Frobenius $F_{X/S}:X\ra X'$ lifts to an $\frakS$-morphism $F:\frakX\ra \frakX'$ such that $\frakX'$ is log smooth over $\frakS$ and flat and locally of finite type over $W.$ We identify the étale sites of $X,$ $X',$ $\frakX$ and $\frakX'$ via $F_{X/S}$ and the canonical immersions $X\ra \frakX$ and $X' \ra \frakX',$ which are universal homeomorphisms. We denote by $\upmu$ the sheaf $\ov{\calM}^{\mathrm{gp}}_X.$ Fix a positive integer $n$ and consider the localization morphism of topoi
\begin{equation}
j_{\upmu}:\frakX_{n,\text{ét}/\upmu}\ra \frakX_{n,\text{ét}}.
\end{equation}
For a sheaf $\F$ of $\frakX_{n,\text{ét}},$ we denote by $\F_{\upmu}$ its pullback along $j_{\upmu}.$
Recall the definition of the $\upmu$-indexed algebra $\calA_{\frakX_n}$ (\ref{example314}), the $\upmu$-indexed subalgebra $\calB_{F_n}$ (\ref{defBF}), the canonical connection $d_{\calA_{\frakX_n}}$ (\ref{dAX}) and the $p$-connection $\theta_{F_n}$ (\ref{defthetaBF}).
We denote by $p\text{-}\boldsymbol{\op{MIC}}(\calB_{F_n})$ the category of $\upmu$-indexed $\calB_{F_n}$-modules $\calE'$ equipped with an integrable $p$-connection $\theta:\calE' \ra \calE'\otimes_{\Ox_{\frakX_n',\upmu}}\omega^1_{\frakX_n'/\frakS_n,\upmu}$ compatible with $\theta_{F_n}$ i.e. such that
$$
\theta(bx')=b\theta(x')+x'\otimes \theta_{F_n}(b),
$$ 
for all local sections $b$ and $x'$ of $\calB_{F_n}$ and $\calE'$ respectively. We also denote by $\boldsymbol{\op{MIC}}(\calA_{\frakX_n})$ the category of $\upmu$-indexed $\calA_{\frakX_n}$-modules $\calE$ equipped with an integrable connection $\nabla:\calE \ra \calE \otimes_{\Ox_{\frakX_n,\upmu}}\omega^1_{\frakX_n/\frakS_n,\upmu}$ which is admissible i.e. compatible with $d_{\calA_{\frakX_n}}$ in the sense that 
$$
\nabla(ax)=a\nabla(x)+x\otimes d_{\calA_{\frakX_n}}(a),
$$
for all local sections $a$ and $x$ of $\calA_{\frakX_n}$ and $\calE$ respectively.
\end{parag}

\subsection{Divided Frobenius map and the Shiho functor}

\begin{parag}
The Frobenius lift $F:\frakX \ra \frakX'$ induces a morphism $dF:F^*\omega^1_{\frakX'/\frakS} \ra \omega^1_{\frakX/\frakS}.$ Its reduction modulo $p$ is equal to $dF_{X/S}$ which vanishes, we deduce that $\op{Im}dF\subset p\omega^1_{\frakX/\frakS}.$ Since $\frakX \ra \frakS$ is log smooth and $\frakX$ is flat over $W,$ we deduce that $\omega^1_{\frakX/\frakS}$ is flat over $W$ and multiplication by $p$ induces an isomorphism $\omega^1_{\frakX/\frakS} \isom p\omega^1_{\frakX/\frakS}.$ We conclude the existence of a unique morphism $p^{-1}dF:F^*\omega^1_{\frakX'/\frakS} \ra \omega^1_{\frakX/\frakS}$ fitting into the commutative diagram
$$
\begin{tikzcd}
F^*\omega^1_{\frakX'/\frakS} \ar{r}{dF} \ar[bend right=30]{rr}{p^{-1}dF} & p\omega^1_{\frakX/\frakS} & \omega^1_{\frakX/\frakS}. \ar[swap]{l}{\times p}
\end{tikzcd}
$$ 
\end{parag}

\begin{lemma}[\cite{SF1} 9.3]\label{lem422}
Consider the composition $X' \xrightarrow{\pi} X \ra \frakX$ where $\pi$ is given in \eqref{diagfrob}. Let $\frakU$ be an étale formal $\frakX$-scheme and $\frakU'$ the unique étale formal $\frakX'$-scheme such that $U'=\frakU' \times_{\frakX'}X'.$ Let $m$ be a local section of $\calM_{\frakX}$ (resp. $\calM_{\frakX_n}$) over $\frakU$ and $m_1$ its image in $\Gamma \left (U,\calM_X \right ).$ After eventually reducing $\frakU,$ there exists a lifting $m'$ of $\pi^{\flat}m_1\in \Gamma \left (U',\calM_{X'} \right )$ to $\Gamma \left (\frakU',\calM_{\frakX'} \right )$ (resp. $\Gamma \left (\frakU',\calM_{\frakX'_n} \right )$). Then there exists an invertible local section $u$ of $\calM_{\frakX}$ (resp. $\calM_{\frakX_n}$) such that $u+pm=F^{\flat}(m')$ (resp. $u+pm=F_n^{\flat}(m')$) and $\alpha_{\frakX}(u)=1+pb$ (resp. $\alpha_{\frakX_n}(u)=1+pb$) for a local section $b$ of $\Ox_{\frakX}$ (resp. $\Ox_{\frakX_n}$).
\end{lemma}

\begin{lemma}
Let $m,$ $m',$ $u$ and $b$ as in \ref{lem422}. Then
\begin{equation}\label{eq4131}
(p^{-1}dF)(\op{dlog}m')=\op{dlog}m+\frac{db}{\alpha_{\frakX}(u)}.
\end{equation}
\end{lemma}

\begin{proof}
This follows from the following computation and the injectivity of multiplication by $p$ on $\omega^1_{\frakX/\frakS}:$
\begin{alignat*}{2}
p(p^{-1}dF)(\op{dlog}m') &= dF(\op{dlog}m') = \op{dlog} F^{\flat}m' = p\op{dlog}m+\op{dlog}u \\
&= p\op{dlog}m +\frac{\alpha_{\frakX}(u)\op{dlog}u}{\alpha_{\frakX}(u)} = p\op{dlog}m + \frac{d\alpha_{\frakX}(u)}{\alpha_{\frakX}(u)} = p\op{dlog}m + p\frac{db}{\alpha_{\frakX}(u)}.
\end{alignat*}
\end{proof}

\begin{parag}
Let $(\calE',\theta)$ be an object of $p\text{-}\boldsymbol{\op{MIC}}(\calB_{\frakX_n/\frakS_n}).$ We consider the projections $\upmu^2\ra \upmu$ and denote by 
$$p_i:\frakX_{n,\text{ét}/\upmu^2} \ra \frakX_{n,\text{ét}/\upmu},\quad p_i':\frakX'_{n,\text{ét}/\upmu^2} \ra \frakX'_{n,\text{ét}/\upmu}$$
the morphisms of topoi they induce.
Let $\zeta$ be the composition
\begin{alignat*}{2}
 \zeta : p_2'^*(\calE' \otimes_{\Ox_{\frakX_n',\upmu}}\omega^1_{\frakX_n'/\frakS_n,\upmu}) \isom & \ (p_2'^*\calE') \otimes_{\Ox_{\frakX_n',\upmu^2}}\omega^1_{\frakX_n'/\frakS_n,\upmu^2} \\
\ra & \ (p_2'^*\calE') \otimes_{\Ox_{\frakX_n',\upmu^2}}(F_n^*\omega^1_{\frakX_n'/\frakS_n})_{\upmu^2}
\ra \ (p_2'^*\calE') \otimes_{\Ox_{\frakX_n',\upmu^2}}\omega^1_{\frakX_n/\frakS_n,\upmu^2},
\end{alignat*}
where the first two morphisms are the canonical ones and the last one is $\op{Id}_{p_2'^*\calE'}\otimes (p^{-1}dF)_{\upmu^2}.$
For sections $s$ and $t$ of $\upmu$ over an étale $\frakX_n$-scheme $U,$ we have
\begin{equation}
\zeta_{(\ov{s},\ov{t})}:\begin{array}[t]{clc}
\calE'_{\ov{s}}\otimes_{\Ox_{U'}}\omega^1_{U'/\frakS_n} & \ra & \calE'_{\ov{s}}\otimes_{\Ox_{U'}}\omega^1_{U/\frakS_n} \\
x\otimes \omega & \mapsto & x\otimes (p^{-1}dF)(\omega).
\end{array}
\end{equation}
For local sections $x'$ and $b$ of $\calE'$ and $\calB_{\frakX_n/\frakS_n}$ respectively, we have
$
\theta(bx')=b\theta(x')+x'\otimes \theta_{\widetilde{F}}(b).
$
Then,
\begin{equation}\label{eq5143}
\zeta \circ p_2'^*\theta(bx') = \begin{cases}
b p_2'^*\zeta \left (p_2'^*\theta(x') \right )+(bx')\otimes \op{dlog}F_n^{\flat}m' & \text{if}\quad b=e_{F_n^{\flat}m'} \\
b\zeta \left (p_2'^*\theta(x') \right ) & \text{if}\quad b=\alpha_{\frakX_n}(m)e_{-m} \\
b\zeta \left (p_2'^*\theta(x') \right ) + x' \otimes d(F_n^{\#}b) & \text{if}\quad b\in \Ox_{\frakX_n'},
\end{cases}
\end{equation}
where $m'$ and $m$ are local sections of $\calM_{\frakX_n'}$ and $\calM_{\frakX_n}$ respectively.
\end{parag}

\begin{lemma}\label{indlem92}
Keep the same hypothesis and notations of the previous paragraph. Then the morphism
\begin{equation}
\nabla_1:\begin{array}[t]{clc} 
 \calA_{\frakX_n} \boxtimes \calE'  & \ra & \left ( \calA_{\frakX_n} \boxtimes \calE' \right ) \otimes_{\Ox_{\frakX_n,\upmu^2}} \omega^1_{\frakX_n/\frakS_n,\upmu^2} \\ 
a\otimes x & \mapsto & a\otimes \zeta \left ( (p_2'^*\theta)(x)\right ) + p_1'^*d_{\calA_{\frakX_n}}(a) \otimes x
\end{array}
\end{equation}
is well-defined and is an integrable connection on the $\Ox_{\frakX_n,\upmu^2}$-module $\calA_{\frakX_n} \boxtimes \calE'.$ In addition, it induces an integrable admissible connection on the $\upmu$-indexed $\calA_{\frakX_n}$-module $\calA_{\frakX_n} \circledast_{\calB_{\frakX_n/\frakS_n}} \calE'$
\begin{equation}\label{eq922}
\nabla: \calA_{\frakX_n} \circledast_{\calB_{\frakX_n/\frakS_n}} \calE' \ra \left (\calA_{\frakX_n} \circledast_{\calB_{\frakX_n/\frakS_n}} \calE' \right ) \otimes_{\Ox_{\frakX_n,\upmu}} \omega^1_{\frakX_n/\frakS_n,\upmu}.
\end{equation}
\end{lemma}

\begin{proof}
To prove the first part, it is sufficient to prove that for any sections $s$ and $t$ of $\calM_{\frakX_n}^{\mathrm{gp}}$ over an étale $\frakX_n$-scheme $U,$ the morphism
$$
\nabla_{1,(\ov{s},\ov{t})}:\begin{array}[t]{clc}
\calA_{\frakX_n,\ov{t}} \otimes_{\Ox_{U'}} \calE'_{\ov{s}} & \ra & \left ( \calA_{\frakX_n,\ov{t}} \otimes_{\Ox_{U'}} \calE'_{\ov{s}} \right ) \otimes_{\Ox_U}\omega^1_{U/\frakS_n} \\
a \otimes x & \mapsto & a\otimes \zeta_{(\ov{s},\ov{t})}(\theta_{\ov{s}}(x))+ d_{\calA_{\frakX_n},\ov{t}}(a)\otimes x
\end{array}
$$
is well-defined and is integrable, which is similar to that of (\cite{SF1} 9.6). For the second part, denote by $\sigma:\upmu^2\ra \upmu$ the addition map. The adjunction morphism
$ \calA_{\frakX_n} \boxtimes \calE' \ra \sigma^* \sigma_!(\calA_{\frakX_n} \boxtimes \calE')$
induces a morphism
$$\left ( \calA_{\frakX_n} \boxtimes \calE' \right ) \otimes_{\Ox_{\frakX_n,\upmu^2}} \omega^1_{\frakX_n/\frakS_n,\upmu^2} \ra \sigma^* \sigma_! \left (\calA_{\frakX_n} \boxtimes \calE' \right )\otimes_{\Ox_{\frakX_n,\upmu^2}} \omega^1_{\frakX_n/\frakS_n,\upmu^2}.$$
Composing it with the connection $\nabla_1,$ we get a connection
$$\calA_{\frakX_n} \boxtimes \calE' \ra \sigma^* \sigma_! \left ( \calA_{\frakX_n} \boxtimes \calE' \right )\otimes_{\Ox_{\frakX_n,\upmu^2}} \omega^1_{\frakX_n/\frakS_n,\upmu^2}.$$
We have a canonical isomorphism
$$\sigma^* \sigma_! \left (\calA_{\frakX_n} \boxtimes \calE' \right ) \otimes_{\Ox_{\frakX_n,\upmu^2}} \omega^1_{\frakX_n/\frakS_n,\upmu^2} \isom \sigma^* \left ( \sigma_! \left (\calA_{\frakX_n} \boxtimes \calE' \right )\otimes_{\Ox_{\frakX_n,\upmu}} \omega^1_{\frakX_n/\frakS_n,\upmu} \right )$$
and so we get a connection
$$\nabla_2:\sigma_! \left (\calA_{\frakX_n} \boxtimes \calE' \right ) \ra \sigma_! \left (\calA_{\frakX_n} \boxtimes \calE' \right )\otimes_{\Ox_{\frakX_n,\upmu}} \omega^1_{\frakX_n/\frakS_n,\upmu}.$$
By \eqref{eq5143} and \ref{rem318}, the composition of $\nabla_2$ with the projection
$$\sigma_! \left (\calA_{\frakX_n} \boxtimes \calE' \right )\otimes_{\Ox_{\frakX_n,\upmu}} \omega^1_{\frakX_n/\frakS_n,\upmu} \ra \left (\calA_{\frakX_n} \circledast_{\calB_{\frakX_n/\frakS_n}} \calE' \right )\otimes_{\Ox_{\frakX_n,\upmu}} \omega^1_{\frakX_n/\frakS_n,\upmu}$$
factors through $\calA_{\frakX_n} \circledast_{\calB_{\frakX_n/\frakS_n}} \calE'$ and we get a connection
$$\nabla:\calA_{\frakX_n} \circledast_{\calB_{\frakX_n/\frakS_n}} \calE' \ra \left ( \calA_{\frakX_n} \circledast_{\calB_{\frakX_n/\frakS_n}} \calE' \right ) \otimes_{\Ox_{\frakX_n,\upmu}} \omega^1_{\frakX_n/\frakS_n,\upmu}.$$
More precisely, for a section $r$ of $\upmu$ over an étale $\frakX_n$-scheme $U,$ the morphism $\nabla_r$ is induced by
\begin{equation}\label{explicitconnimage}
\begin{array}[t]{clc}
\bigoplus\limits_{s+t=r}\calA_{\frakX_n,t} \otimes_{\Ox_{U'}} \calE'_s & \ra & \bigoplus\limits_{s+t=r} \left (\calA_{\frakX_n,t} \otimes_{\Ox_{U'}} \calE'_s \right ) \otimes_{\Ox_{U}}\omega^1_{U/\frakS_n} \\
a\otimes x & \mapsto & \nabla_{1,(s,t)}(a\otimes x).
\end{array}
\end{equation}
The properties of integrability and admissibility of $\nabla$ are inherited from $\nabla_1.$
\end{proof}

\begin{parag}
Keep the assumptions and notations of \ref{indparag91}. By \ref{indlem92}, we have a functor
\begin{equation}\label{Phiind}
\Phi :\begin{array}[t]{clc} p\text{-}\boldsymbol{\op{MIC}}(\calB_{\frakX_n/\frakS_n}) & \ra & \boldsymbol{\op{MIC}}(\calA_{\frakX_n}) \\
(\calE',\theta') & \mapsto & \left ( \calA_{\frakX_n} \circledast_{\calB_{\frakX_n/\frakS_n}} \calE', \nabla \right ), \end{array}
\end{equation}
where $\nabla$ is the connection defined in (\ref{eq922}).
\end{parag}

\begin{parag}
Recall the logarithmic formal scheme $R_{\frakX'/\frakS}$ defined in \ref{prop124}. Let $P_{\frakX/\frakS,0}$ be the log PD envelope of $\frakX \ra \frakX\times_{\frakS}\frakX$ (\cite{Kat89} 5.4) and $P_{\frakX'/\frakS,1}$ the PD envelope of $\frakX' \ra R_{\frakX'/\frakS},$ equipped with the logarithmic structure pullback of that of $R_{\frakX'/\frakS}.$ The canonical immersion $\frakX' \ra P_{\frakX'/\frakS,1}$ is then strict and hence so are the projections $P_{\frakX'/\frakS,1} \ra \frakX'.$
\end{parag}

\begin{proposition}\label{prop432}
There exists a unique morphism of formal schemes $\varphi : P_{\frakX/\frakS,0} \ra R_{\frakX'/\frakS}$ fitting into the commutative diagram
\begin{equation}\label{diag4321}
\begin{tikzcd}
\frakX \ar{r} \ar[swap]{d}{F} & P_{\frakX/\frakS,0} \ar{r} \ar{d}{\varphi} & \frakX \times_{\frakS} \frakX \ar{d}{F\times_{\frakS}F} \\
\frakX' \ar{r} & R_{\frakX'/\frakS} \ar{r} & \frakX'\times_{\frakS}\frakX'.
\end{tikzcd}
\end{equation}
\end{proposition}

\begin{proof}
We prove the existence of $\varphi$ locally.
Let $\frakX' \xrightarrow{\Delta'} \frakY' \ra \frakX' \times_{\frakS} \frakX'$ be a local factorization of the diagonal immersion $\frakX' \ra \frakX'\times_{\frakS}\frakX'$ into a strict closed immersion followed by a log étale separated morphism (\ref{lem122}). Let $\frakY$ be the fiber product of $\frakY' \ra \frakX'\times_{\frakS}\frakX'$ and $F\times_{\frakS}F:\frakX\times_{\frakS}\frakX \ra \frakX'\times_{\frakS}\frakX'.$ Let $G:\frakY \ra \frakY'$ be the projection. Let $\Delta:\frakX \ra \frakY$ be the morphism induced by the diagonal immersion $\frakX \ra \frakX\times_{\frakS}\frakX$ and the composition $\frakX \xrightarrow{F} \frakX' \xrightarrow{\Delta'} \frakY'.$ The morphism $\Delta$ is an exact immersion. Indeed, $F$ is exact and the immersion $\Delta'$ is strict hence exact. Since $\Delta' \circ F=G\circ \Delta,$ we deduce, by (\cite{Ogus2018} III 2.2.1), that $\Delta$ is exact hence strict. We can suppose that $P_{\frakX/\frakS,0} \ra \frakX\times_{\frakS}\frakX$ and $R_{\frakX'/\frakS} \ra \frakX'\times_{\frakS} \frakX'$ factor through $\frakY$ and $\frakY'$ respectively. To prove that $P_{\frakX/\frakS,0} \ra \frakY \xrightarrow{G} \frakY'$ factors uniquely through $R_{\frakX'/\frakS},$ it is sufficient to check the flatness of $P_{\frakX/\frakS,0}$ over $W$ and that $\calI' \Ox_{P_{\frakX/\frakS,0}}\subset p\Ox_{P_{\frakX/\frakS,0}},$ where $\calI'$ is the ideal of $\Delta'.$ The flatness follows from the log smoothness of $\frakX \ra \frakS$ and the flatness of $\frakX$ over $W.$ The inclusion $\calI' \Ox_{P_{\frakX/\frakS,0}}\subset p\Ox_{P_{\frakX/\frakS,0}}$ is a consequence of the computation \eqref{eqGlocal} below.

We now prove this morphism is independent of the choice of the factorization $\frakX' \ra \frakY' \ra \frakX' \times_{\frakS}\frakX'.$ Suppose that we have two such factorizations $\frakX' \xrightarrow{\Delta_i'} \frakY_i' \ra \frakX_i' \times_{\frakS} \frakX_i'$ and let $\varphi_i:P_{\frakX/\frakS,0} \ra R_{\frakX'/\frakS}$ be the morphism constructed as above. Factoring the morphism $\frakX' \ra \frakY'_1\times_{\frakX'\times_{\frakS}\frakX'}\frakY'_2,$ induced by $\Delta'_1$ and $\Delta'_2,$ locally into a strict closed immersion followed by a log étale separated morphism (\ref{lem122}), we obtain a third factorization $\frakX' \xrightarrow{\Delta'_{12}} \frakZ \ra \frakX'\times_{\frakS} \frakX'$ and a third morphism $\varphi_{12}:P_{\frakX/\frakS,0} \ra R_{\frakX'/\frakS}.$ The universal property of the dilatation $R_{\frakX'/\frakS}$ implies that $\varphi_1=\varphi_2=\varphi_{12}.$ We hence obtain the morphism of formal schemes $\varphi$ by a gluing argument.
\end{proof}

\begin{proposition}\label{prop433}
The morphism of formal schemes $P_{\frakX/\frakS,0} \ra R_{\frakX'/\frakS},$ given in \ref{prop432}, factors uniquely through $P_{\frakX'/\frakS,1},$ yielding a morphism of formal schemes $\varphi:P_{\frakX/\frakS,0} \ra P_{\frakX'/\frakS,1}.$
\end{proposition}

\begin{proof}
This is a consequence of the universal property of the PD envelope.
\end{proof}

\begin{remark}\label{loccoord}
We give a local description of the morphism $\varphi:P_{\frakX/\frakS,0} \ra P_{\frakX'/\frakS,1},$ given in \ref{prop433}. Suppose that we have local coordinates $m_1,\hdots,m_d$ of $X/S$ and that $\frakX' \ra \frakX'\times_{\frakS}\frakX'$ factors as a strict closed immersion $\Delta':\frakX' \ra \frakY'$ followed by a separated log étale morphism $\frakY'\ra \frakX'\times_{\frakS}\frakX'.$ We deduce a factorization $\frakX \xrightarrow{\Delta} \frakY \ra \frakX\times_{\frakS}\frakX$ and a morphism $G:\frakY \ra \frakY'$ as in the proof of \ref{prop432}. For $1\le i\le d,$ set $m_i'=\pi^{\flat}m_i.$ Suppose that $m_i$ and $m_i'$ lift to local sections $\widetilde{m}_i$ and $\widetilde{m}_i'$ of $\calM_{\frakX}$ and $\calM_{\frakX'}$ respectively. Consider the commutative diagram with exact rows
$$
\begin{tikzcd}
0 \ar{r} & F^{-1} \Delta'^{-1}(1+\calI') \ar{r}{\lambda'} \ar{d}{\Delta^{-1}G} & F^{-1}\Delta'^{-1}\calM_{\frakY'} \ar{r}{F^{-1}\Delta'^{\flat}} \ar{d}{\Delta^{-1}G^{\flat}} & F^{-1}\calM_{\frakX'} \ar{d}{F^{\flat}} \ar{r} & 0 \\
0 \ar{r} & \Delta^{-1}(1+\calI) \ar{r}{\lambda} & \Delta^{-1}\calM_{\frakY} \ar{r}{\Delta^{\flat}} & \calM_{\frakX} \ar{r} & 0,
\end{tikzcd}
$$
where $\calI$ and $\calI'$ are the ideals of $\Delta$ and $\Delta'$ respectively. Let $p_1,p_2:\frakY \ra \frakX$ and $p_1',p_2':\frakY' \ra \frakX'$ be the projections. For $1\le i\le d,$ let $\widetilde{\eta}_i$ and $\widetilde{\eta}_i'$ be the unique local sections of $\Delta^{-1}\calI$ and $\Delta'^{-1}\calI'$ respectively such that
\begin{equation}\label{4341}
\lambda(1+\widetilde{\eta}_i)+\Delta^{-1}p_1^{\flat}\widetilde{m}_i=\Delta^{-1}p_2^{\flat}\widetilde{m}_i,\quad \lambda'(1+\widetilde{\eta}_i')+\Delta'^{-1}p_1'^{\flat}\widetilde{m}'_i=\Delta'^{-1}p_2'^{\flat}\widetilde{m}_i'.
\end{equation}
Then, $\Ox_{P_{\frakX'/\frakS,1}}=\Ox_{\frakX'} \langle \langle \frac{\widetilde{\eta}'_1}{p},\hdots,\frac{\widetilde{\eta}'_d}{p} \rangle \rangle$ and $\Ox_{P_{\frakX/\frakS,0}}=\Ox_{\frakX} \langle \langle \widetilde{\eta}_1,\hdots,\widetilde{\eta}_d \rangle \rangle.$
We drop the notation $\Delta^{-1}$ and $\Delta'^{-1}$ for simplicity. By \ref{lem422}, there exist local sections $u_i$ and $b_i$ of $\calM_{\frakX}$ and $\Ox_{\frakX}$ respectively such that $u_i$ is invertible, $F^{\flat}(\widetilde{m}_i')=p\widetilde{m}_i+u_i$ and $\alpha_{\frakX}(u_i)=1+pb_i.$ Then
$$
G^{\flat}p_i'^{\flat}\widetilde{m}_j'=p_i^{\flat}F^{\flat}\widetilde{m}_j'=p_i^{\flat}(p\widetilde{m}_j+u_j)=pp_i^{\flat}(\widetilde{m}_j)+p_i^{\flat}u_j.
$$
Applying $G^{\flat}$ to the second equality of \eqref{4341} and using the first, we get
$$
\lambda(G(1+\widetilde{\eta}_i'))+pp_1^{\flat}\widetilde{m}_i+p_1^{\flat}u_i=pp_2^{\flat}\widetilde{m}_i+p_2^{\flat}u_i=pp_1^{\flat}\widetilde{m}_i+p\lambda(1+\widetilde{\eta}_i)+p_2^{\flat}u_i.
$$
Since $\calM_{\frakY}$ is integral, and since $\alpha_{\frakY}\circ p_i^{\flat} = p_i \circ \alpha_{\frakX},$ we get
$$
\lambda(G(1+\widetilde{\eta}_i'))=\lambda\left ( (1+\widetilde{\eta}_i)^p\right )+p_2^{\flat}u_i-p_1^{\flat}u_i=\lambda\left ( (1+\widetilde{\eta}_i)^p\right )+\lambda \left ( \frac{1+pp_2(b_i)}{1+pp_1(b_i)}\right )=\lambda \left ( (1+\widetilde{\eta}_i)^p\frac{1+pp_2(b_i)}{1+pp_1(b_i)} \right ).
$$
We deduce that $G(1+\widetilde{\eta}_i')=(1+\widetilde{\eta}_i)^p\frac{1+pp_2(b_i)}{1+pp_1(b_i)}.$ Then, we have
\begin{alignat}{2}
G(\widetilde{\eta}_i') &= \left ( \widetilde{\eta}_i^p +\sum_{k=1}^{p-1}\begin{pmatrix}p\\k\end{pmatrix}\widetilde{\eta}_i^k\right )\frac{1+pp_2(b_i)}{1+pp_1(b_i)}+\frac{1+pp_2(b_i)}{1+pp_1(b_i)}-1 \nonumber\\
&= \left ( \widetilde{\eta}_i^p +\sum_{k=1}^{p-1}\begin{pmatrix}p\\k\end{pmatrix}\widetilde{\eta}_i^k\right )\frac{1+pp_2(b_i)}{1+pp_1(b_i)}+p\frac{p_2^{\#}b_i-p_1^{\#}b_i}{1+pp_1^{\#}b_i}. \label{eqGlocal}
\end{alignat}
We conclude that $\varphi:\Ox_{P_{\frakX'/\frakS,1}}=\Ox_{\frakX'} \langle \langle \frac{\widetilde{\eta}'_1}{p},\hdots,\frac{\widetilde{\eta}'_d}{p} \rangle \rangle\ra \Ox_{\frakX} \langle \langle \widetilde{\eta}_1,\hdots,\widetilde{\eta}_d \rangle \rangle=\Ox_{P_{\frakX/\frakS,0}}$ satisfies
\begin{equation}\label{eq4243}
\varphi \left ( \frac{\widetilde{\eta}_i'}{p} \right ) = \left ( (p-1)!\widetilde{\eta}_i^{[p]} +\sum_{k=1}^{p-1}\frac{(p-1)!}{k!(p-k)!} \widetilde{\eta}_i^k\right )\frac{1+pp_2^{\#}b_i}{1+pp_1^{\#}b_i}+\frac{p_2^{\#}b_i-p_1^{\#}b_i}{1+pp_1^{\#}b_i}.
\end{equation}
Let $n$ be a positive integer. Let $\varphi_n,$ $\xi_i$ and $\eta_i$ be the reductions of $\varphi,$ $\frac{\widetilde{\eta}_i'}{p}$ and $\widetilde{\eta}_i$ modulo $p^n$ respectively. Abusively denote by $m_i$ and $m_i'$ the images of $\widetilde{m}_i$ and $\widetilde{m}'_i$ in $\calM_{\frakX_n}$ and $\calM_{\frakX_n'}$ respectively. Let $P_n$ and $P_n'$ be the schemes obtained from $P_{\frakX/\frakS,0}$ and $P_{\frakX'/\frakS,1}$ by reduction modulo $p^n$ and $\ov{\calI}$ and $\ov{\calI}'$ the PD ideals of $P_n$ and $P_n'$ respectively. Then, by \eqref{eq4243}, we have $\varphi_n \left (\ov{\calI}' \right ) \subset \ov{\calI}.$ Modulo $\ov{\calI}^{[2]},$ since $\frac{1+pp_2^{\#}b_i}{1+pp_1^{\#}b_i}-1$ is a local section of $\ov{\calI},$ we have
\begin{equation}\label{eq4243*}
\varphi_n(\xi_i)=\eta_i+\frac{p_2^{\#}b_i-p_1^{\#}b_i}{1+pp_1^{\#}b_i}.
\end{equation}
Via the isomorphism $\ov{\calI}/\ov{\calI}^{[2]}\isom \omega^1_{\frakX_n/\frakS_n},$ the local sections $p_2^{\#}b_i-p_1^{\#}b_i$ and $\eta_i$ identify with $db_i$ and $\op{dlog}m_i$ respectively. The section $\varphi_n(\xi_i)$ then identifies with
\begin{equation}\label{eq4243**}
\varphi_n(\xi_i)=\op{dlog}m_i+\frac{db_i}{\alpha_{\frakX_n}(u_i)}.
\end{equation}
In addition, via the isomorphism $\ov{\calI}'/\ov{\calI}'^{[2]}\isom \omega^1_{\frakX'_n/\frakS_n},$ the local section $\xi_i$ identifies with $\op{dlog}m_i'.$ The equalities \eqref{eq4243**} and \eqref{eq4131} yield the commutativity of the diagram
\begin{equation}\label{phidF}
\begin{tikzcd}
\ov{\calI}' \ar{d}{\varphi_n} \ar{r} & \ov{\calI}'/\ov{\calI}'^{[2]} \ar{d} \ar{r}{\sim} & \omega^1_{\frakX'_n/\frakS_n} \ar{d}{p^{-1}dF} \\
\ov{\calI} \ar{r} & \ov{\calI}/\ov{\calI}^{[2]} \ar{r}{\sim} & \omega^1_{\frakX_n/\frakS_n}.
\end{tikzcd}
\end{equation}
\end{remark}

\begin{remark}\label{philog}
Consider the morphism of formal schemes $\varphi:P_{\frakX/\frakS,0} \ra P_{\frakX'/\frakS,1}$ given in \ref{prop433}. Denote by $p_1',p_2':P_{\frakX'/\frakS,1} \ra \frakX'$ and $p_1,p_2:P_{\frakX/\frakS,0} \ra \frakX$ the projections. Since $p_i$ and $p_i'$ are strict, they induce isomorphisms of logarithmic structures $p_i^*\calM_{\frakX} \isom \calM_{P_{\frakX/\frakS,0}}$ and $p_i'^*\calM_{\frakX'} \isom \calM_{P_{\frakX'/\frakS,1}}.$ Consider the morphism
$$
\varphi_i^{\flat}=p_i^*F^{\flat}:\varphi^*\calM_{P_{\frakX'/\frakS,1}} = \varphi^*p_i'^*\calM_{\frakX'} = p_i^*F^*\calM_{\frakX'} \ra p_i^*\calM_{\frakX}=\calM_{P_{\frakX/\frakS,0}}.
$$
We then have a morphism of logarithmic formal schemes $(\varphi,\varphi_i^{\flat}).$ We denote it by $\varphi_i.$ For a positive integer $n,$ we denote by $\varphi_{i,n}$ the reduction of $\varphi_i$ modulo $p^n.$
\end{remark}

\begin{proposition}\label{faithfulflatphi}
Let $n$ be a positive integer, $\varphi:P_n \ra P_n'$ the morphism of schemes induced by \ref{prop433} modulo $p^n$ and $\calP_n' \ra \calP_n$ the corresponding morphism of structural rings. Then, the morphisms of $\upmu$-indexed algebras
\begin{equation}\label{eq5260}
\calP_{n,\upmu}' \otimes_{\Ox_{\frakX_n',\upmu}} \calA_{\frakX_n} \ra \calP_{n,\upmu} \otimes_{\Ox_{\frakX_n,\upmu}} \calA_{\frakX_n},
\end{equation}
\begin{equation}\label{eq5261}
\calP_{n,\upmu}' \otimes_{\Ox_{\frakX_n',\upmu}} \calB_{\frakX_n / \frakS_n} \ra \calP_{n,\upmu} \otimes_{\Ox_{\frakX_n,\upmu}} \calA_{\frakX_n} 
\end{equation}
are faithfully flat. In addition, the morphisms $\varphi_{i,n},$ given in \ref{philog}, are log flat.
\end{proposition}

\begin{proof}
Being faithfully flat is a statement about the exactness of sequences of indexed modules, it can hence be checked locally. We can thus assume that the hypothesis of \ref{loccoord} is satisfied. We keep the same notations of \ref{loccoord}. In this case, we have
$$
\calP_{n,\upmu}' \otimes_{\Ox_{\frakX_n',\upmu}} \calB_{\frakX_n / \frakS_n}=\calB_{\frakX_n/\frakS_n} \langle \xi_1,\hdots ,\xi_d\rangle,\quad 
\calP_{n,\upmu} \otimes_{\Ox_{\frakX_n,\upmu}} \calA_{\frakX_n} = \calA_{\frakX_n}\langle \eta_1,\hdots,\eta_d \rangle.
$$
The morphism \eqref{eq5261} then factors as follows
$$
\calB_{\frakX_n/\frakS_n} \langle \xi_1,\hdots ,\xi_d\rangle \xhookrightarrow{\tau} \calA_{\frakX_n} \langle \xi_1,\hdots ,\xi_d\rangle \xrightarrow{\nu} \calA_{\frakX_n}\langle \eta_1,\hdots,\eta_d \rangle,
$$
where $\tau$ is the canonical embedding and $\nu$ is the morphism of $\upmu$-indexed $\calA_{\frakX_n}$-algebras sending $\xi_i$ to
$$\left ( (p-1)!\widetilde{\eta}_i^{[p]} +\sum_{k=1}^{p-1}\frac{(p-1)!}{k!(p-k)!} \widetilde{\eta}_i^k\right )\frac{1+pp_2^{\#}b_i}{1+pp_1^{\#}b_i}+\frac{p_2^{\#}b_i-p_1^{\#}b_i}{1+pp_1^{\#}b_i}.$$
Also, in this case, \eqref{eq5260} identifies with $\nu.$ By \ref{propABlocfree}, $\calA_{\frakX_n}$ is locally free of finite rank over $\calB_{\frakX_n/\frakS_n}.$ Then, by \ref{example305}, $\tau$ is faithfully flat. For $\nu,$ it is sufficient, by \ref{example305}, to prove that $\Ox_{\frakX_n} \langle \xi_1,\hdots,\xi_d \rangle \ra \Ox_{\frakX_n} \langle \eta_1,\hdots,\eta_d \rangle$ is faithfully flat. This is the morphism of structure rings of $P_n \ra P_n' \times_{\frakX_n'} \frakX_n$ whose underlying map of sets is surjective. For the flatness, since $\frakX_n$ is flat over $W/p^nW$ by assumption, we can work modulo $p.$ We abusively denote by $\xi_i,$ $\eta_i$ and $b_i$ the classes of $\xi_i,$ $\eta_i$ and $b_i$ modulo $p.$ By \eqref{eq4243}, modulo $p,$ we have
$$
\nu(\xi_i)=-\eta_i^{[p]}+\sum_{k=1}^{p-1} \frac{(-1)^{k+1}}{k}\eta_i^k+p_2^{\#}b_i-p_1^{\#}b_i.
$$
Let $\calI_1$ be the ideal of the strict closed immersion $\Delta_1:X \ra Y.$ Since the composition $Y \ra X\times_SX \ra X$ is log smooth and strict, it is smooth and the immersion $\Delta_1$ is hence regular. It follows that the canonical morphism
$$
\bigoplus_{n\in \N} \op{Sym}^n_{\Ox_X} \left ( \calI_1/\calI_1^2 \right ) \ra \bigoplus_{n\in \N} \calI_1^n/\calI_1^{n+1}
$$
is an isomorphism. Using this and the fact that $\eta_i^p=p!\eta_i^{[p]}=0$ for all $1\le i\le d,$ we prove that there exist local sections $a_I$ of $\Ox_{X}$ such that
$p_2^{\#}b_i-p_1^{\#}b_i=\sum_{I\in \llbracket 0,\hdots,p-1 \rrbracket ^d,\ I\neq 0}a_I\eta^I.$ It follows that $\nu$ is of the form
$$
\nu(\xi_i)=-\eta_i^{[p]}+\sum_{\substack{I\in \llbracket 0,p-1\rrbracket^d\\ I\neq 0}}c_I\eta^I,
$$
for local sections $c_I$ of $\Ox_X.$ We conclude that $\Ox_X\langle \eta_1,\hdots,\eta_d\rangle$ is free over $\Ox_X \langle \xi_1,\hdots,\xi_d \rangle$ of basis $\left (\eta^{[I]} \right )_{I\in \llbracket 0,p-1 \rrbracket^d}$ (see the proof of \cite{Shiho} 3.3 for more details) and $\nu$ is hence flat.

For the log flatness, note that $P_n \ra P_n'$ factors as $P_n \ra P_n'\times_{\frakX_n'} \frakX_n\ra P_n'.$ The second morphism is log flat by \ref{logflatlift} and the first is flat, by the computation above, and strict, hence log flat too.
\end{proof}

\subsection{Indexed stratifications}

\begin{parag}
Recall the hypothesis and notation of \ref{indparag91} and the logarithmic formal scheme $R_{\frakX'/\frakS}$ given in \ref{prop124}. Fix a positive integer $n.$ Let $P_{\frakX/\frakS,0}$ be the log PD envelope of $\frakX \ra \frakX\times_{\frakS}\frakX$ and $P_{\frakX'/\frakS,1}$ the PD envelope of $\frakX' \ra R_{\frakX'/\frakS}$ equipped with the logarithmic structure pullback of that of $R_{\frakX'/\frakS}.$ We denote by $P_n$ and $P_n'$ their respective reductions modulo $p^n$ and by $\calP_n$ and $\calP'_n$ the respective structural rings of these reductions. Let $\upmu=\ov{\calM}_X^{\mathrm{gp}}$ and identify the underlying étale sites of $X,$ $X',$ $\frakX_n,$ $\frakX_n',$ $P_n$ and $P_n'$ via the Frobenius $F_{X/S}:X\ra X'$ and the canonical immersions.

Recall the result \ref{thmKato} (resp. \ref{equivstrat}) stating that the data of a $\calP_n$-stratifications (resp. $\calP_n'$-stratification) on an $\Ox_{\frakX_n}$-module (resp. $\Ox_{\frakX_n'}$-module) is equivalent to the data of an integrable quasi-nilpotent connection (resp. integrable quasi-nilpotent $p$-connection) on this same module.

Consider the canonical projections 
$$p_1,p_2:P_n\ra \frakX_n,\quad p_1',p_2':P'_n \ra \frakX_n'.$$
We denote by $\calP'^{\calB}_n,$ $\prescript{\calB}{}\calP'_n,$ $\calP^{\calA}_n$ and $\prescript{\calA}{}\calP_n$ the $\upmu$-indexed algebras
\begin{equation}\label{PB}
\begin{matrix} \calP'^{\calB}_n=\calP'_{n,\upmu}\otimes_{\Ox_{\frakX'_n,\upmu}}\calB_{\frakX_n/\frakS_n},\quad \prescript{\calB}{}\calP'_n=\calB_{\frakX_n/\frakS_n}\otimes_{\Ox_{\frakX'_n,\upmu}}\calP'_{n,\upmu} \\
\calP^{\calA}_n=\calP_{n,\upmu}\otimes_{\Ox_{\frakX_n,\upmu}}\calA_{\frakX_n},\quad \prescript{\calA}{}\calP_n=\calA_{\frakX_n}\otimes_{\Ox_{\frakX_n,\upmu}}\calP_{n,\upmu}.
\end{matrix}
\end{equation}
Note that the projections $p_1$ and $p_2$ are strict and so, by \ref{strictA}, they induce isomorphisms of $\upmu$-indexed algebras
$
\calP^{\calA}_n=p_1^*\calA_{\frakX_n} \isom \calA_{P_n}$ and $\prescript{\calA}{}\calP_n=p_2^*\calA_{\frakX_n} \isom \calA_{P_n}.
$
Next, we introduce and study indexed stratifications and their relation to connections.
\end{parag}

\begin{definition}\label{def622}
Let $\calE$ be an $\Ox_{\frakX_n,\upmu}$-module and $\lambda \in \Gamma \left ( \frakX_n,\Ox_{\frakX_n} \right ).$ A $\lambda$-connection on $\calE$ is a morphism of abelian groups $\nabla : \calE \ra \calE \otimes_{\Ox_{\frakX_n,\upmu}} \omega^1_{\frakX_n / \frakS_n,\upmu}$ satisfying the Leibniz rule
$$
\nabla (ax) = a\nabla(x) + \lambda x\otimes da,
$$
for all local sections $x$ and $a$ of $\calE$ and $\Ox_{\frakX_n,\upmu}$ respectively. We say that $\nabla$ is integrable (resp. integrable and quasi-nilpotent) if the $\lambda$-connection $\nabla_s$ on $\calE_s$ is integrable (resp. integrable and quasi-nilpotent) for all local sections $s$ of $\upmu$ (\ref{parKhaminei} and \ref{def227}).
\end{definition}

\begin{definition}
Let $\calE$ be an $\Ox_{\frakX_n,\upmu}$-module (resp. $\Ox_{\frakX_n',\upmu}$-module). A \emph{$\calP_{n,\upmu}$-stratification on $\calE$} (resp. \emph{$\calP_{n,\upmu}'$-stratification on $\calE$}) is a $\calP_{n,\upmu}$-linear (resp. $\calP'_{n,\upmu}$-linear) isomorphism
$$\epsilon:\calP_{n,\upmu}\otimes_{\Ox_{\frakX_n,\upmu}}\calE \isom \calE\otimes_{\Ox_{\frakX_n,\upmu}}\calP_{n,\upmu}\quad \left (\text{resp.}\ \epsilon:\calP_{n,\upmu}'\otimes_{\Ox_{\frakX_n',\upmu}}\calE \isom \calE\otimes_{\Ox_{\frakX_n',\upmu}}\calP_{n,\upmu}' \right )$$
satisfying the following two conditions:
\begin{enumerate}
\item The base change of $\epsilon$ by the canonical projection $\calP_{n,\upmu} \ra \Ox_{\frakX_n,\upmu}$ (resp. $\calP_{n,\upmu}' \ra \Ox_{\frakX_n',\upmu}$) is equal to the identity of $\calE.$ 
\item $p_{13}^*\epsilon=p_{12}^*\epsilon \circ p_{23}^*\epsilon,$ where $p_{ij}$ is the morphism
$P_{n}\times_{\frakX_n} P_n \ra P_n$ (resp. $P_n'\times_{\frakX_n'}P_n'\ra P_n'$),
induced by the projection $\frakX_n \times_{\frakS_n} \frakX_n \times_{\frakS_n} \frakX_n \ra \frakX_n \times_{\frakS_n} \frakX_n$ (resp. $\frakX'_n \times_{\frakS_n} \frakX'_n \times_{\frakS_n} \frakX'_n \ra \frakX'_n \times_{\frakS_n} \frakX'_n$) on the factors $i$ and $j.$
\end{enumerate}
\end{definition}

\begin{proposition}\label{prop12.3}
Let $\calE$ be an $\Ox_{\frakX_n,\upmu}$-module (resp. $\Ox_{\frakX_n',\upmu}$-module). The following data are equivalent:
\begin{enumerate}
\item A $\calP_{n,\upmu}$-stratification (resp. $\calP_{n,\upmu}'$-stratification) on $\calE.$
\item An integrable quasi-nilpotent connection (resp. $p$-connection) on $\calE.$
\end{enumerate}
\end{proposition}

\begin{proof}
We prove the equivalence for $p$-connections. The case of connections is similar.
The data of an integrable quasi-nilpotent $p$-connection
$\nabla:\calE \ra \calE\otimes_{\Ox_{\frakX'_n,\upmu}}\omega^1_{\frakX'_n/\frakS_n,\upmu}$
is equivalent to the data of integrable quasi-nilpotent connections
$$\nabla_s:\calE_s \ra \calE_s\otimes_{\Ox_{U'}}\omega^1_{U'/\frakS_n}$$
for every section $s\in\Gamma(U',\upmu)$ over an étale $\frakX_n'$-scheme $U',$ which are functorial in $s.$
Similarly, the data of an indexed $\calP_{n,\upmu}'$-stratification
$\epsilon:\calP_{n,\upmu}'\otimes_{\Ox_{\frakX'_n,\upmu}}\calE \isom \calE \otimes_{\Ox_{\frakX_n',\upmu}}\calP_{n,\upmu}'$
is equivalent to the data of $\calP_n'$-stratifications
$$\epsilon_s:\calP_{n}'\otimes_{\Ox_{U'}}\calE_s \isom \calE_s \otimes_{\Ox_{U'}}\calP_{n}'$$
for every section $s\in\Gamma(U',\upmu)$ over an étale $\frakX_n'$-scheme $U',$ which are functorial in $s.$ The result then follows from \ref{equivstrat}.
\end{proof}

\begin{definition}
We denote by
\begin{equation}\label{epsilonA}
\epsilon_{\calA_{\frakX_n}}:\calP_{n,\upmu}\otimes_{\Ox_{\frakX_n,\upmu}}\calA_{\frakX_n} \isom \calA_{\frakX_n}\otimes_{\Ox_{\frakX_n,\upmu}}\calP_{n,\upmu},\quad \epsilon_{\calB_{\frakX_n/\frakS_n}}':\calP_{n,\upmu}'\otimes_{\Ox_{\frakX_n',\upmu}}\calB_{\frakX_n/\frakS_n} \isom \calB_{\frakX_n/\frakS_n} \otimes_{\Ox_{\frakX_n',\upmu}}\calP_{n,\upmu}',
\end{equation}
the stratifications corresponding respectively to the canonical connection $d_{\calA_{\frakX_n}}$ and the $p$-connection $\theta_{F_n}$ given in \ref{defthetaBF}.
\end{definition}

\begin{remark}\label{remdefeA}
Let $p_1,p_2:P_n \ra \frakX_n$ the two projections. Both of them are strict so, by \ref{strictA}, we have isomorphisms $\nu_i:p_i^*\calA_{\frakX_n} \isom \calA_{P_n}.$ The stratification $\epsilon_{\calA_{\frakX_n}}$ is then given by the composition
$$
\epsilon_{\calA_{\frakX_n}} : \calP_{n,\upmu} \otimes_{\Ox_{\frakX_n,\upmu}} \calA_{\frakX_n} = p_2^*\calA_{\frakX_n} \xrightarrow{\nu_2} \calA_{P_n} \xrightarrow{\nu_1^{-1}} p_1^*\calA_{\frakX_n}=\calA_{\frakX_n} \otimes_{\Ox_{\frakX_n,\upmu}} \calP_{n,\upmu}.
$$
Indeed, if $s$ is a local section of $\calM_{\frakX_n}^{\mathrm{gp}},$ then $\nu_i(1\otimes e_s)=e_{p_i^{\flat}s}.$ Let $\iota:\frakX_n \ra P_n$ be the canonical immersion and $\ov{\calI}$ its ideal. By (\cite{SF3} (4.3.1)), we have the following commutative diagram, where the upper row is an exact sequence of monoids:
$$
\begin{tikzcd}
0 \ar{r} & 1+\ov{\calI} \ar{r} \ar[equal]{d} & \calM_{P_n} \ar{r} \ar[hook]{d} & \calM_{\frakX_n} \ar[hook]{d} \ar{r} & 0 \\
0 \ar{r} & 1+\ov{\calI} \ar{r}{\alpha_{P_n}^{-1}} & \calM_{P_n}^{\mathrm{gp}} \ar{r}{\iota^{\flat}}  & \calM_{\frakX_n}^{\mathrm{gp}}  \ar{r} & 0. 
\end{tikzcd}
$$
It follows that the lower sequence of abelian groups is exact.
There exists a unique local section $\eta$ of $\ov{\calI}$ such that $\alpha_{P_n}^{-1}(1+\eta)+p_1^{\flat}s=p_2^{\flat}s.$ By (\cite{Kat89} 5.8.1 and 6.5), the isomorphism $\ov{\calI}/\ov{\calI}^{[2]} \isom \omega^1_{\frakX_n/\frakS_n}$ identifies the class of $\eta$ modulo $\ov{\calI}^{[2]}$ with $\op{dlog}s.$ Let $\pi:\calP_n \ra \calP_n/\left ( \ov{\calI}^{[2]} \right )$ be the projection. We then have
\begin{alignat*}{2}
& (\op{Id} \otimes \pi) \left ( \nu_1^{-1} \circ \nu_2(1\otimes e_s)-e_s\otimes 1 \right ) \\ =& (\op{Id} \otimes \pi ) \left ( \nu_1^{-1}(e_{p_2^{\flat}s})-e_s\otimes 1\right ) 
= (\op{Id} \otimes \pi )\left ( \nu_1^{-1} \left ( (1+\eta)e_{p_1^{\flat}s} \right ) - e_s\otimes 1\right )=e_s \otimes \op{dlog}s.
\end{alignat*}
We deduce, by \eqref{dAX} and \ref{thmKato} that the composition $\nu_1^{-1} \circ \nu_2$ corresponds to $d_{\calA_{\frakX_n}}.$
\end{remark}

\begin{lemma}\label{lemAP}
Equip $P_n$ and $P_n'$ with the logarithmic structure $p_1^*\calM_{\frakX_n}$ and $p_1'^*\calM_{\frakX_n'}$ (resp. $p_2^*\calM_{\frakX_n}$ and $p_2'^*\calM_{\frakX_n'}$) respectively. Consider the morphism of logarithmic schemes $\varphi_{1,n}$ (resp. $\varphi_{2,n}$) given in \ref{philog}. The projections $p_1$ and $p_1'$ (resp. $p_2$ and $p_2'$) induce isomorphisms of $\upmu$-indexed algebras
\begin{equation}
\calP_n^{\calA} \circledast_{\calP_n'^{\calB}} \calP_n^{\calA} \isom \calA_{P_n \times_{P_n'} P_n},\quad \left (\text{resp.}\ \prescript{\calA}{} \calP_n \circledast_{\prescript{\calB}{}\calP_n'} \prescript{\calA}{}\calP_n \isom \calA_{P_n \times_{P_n'} P_n} \right ).
\end{equation}
\end{lemma}

\begin{proof}
Both assertions are proved similarly so we only prove the first one. Let $g:P_n \times_{P_n'} P_n \ra \frakX_n \times_{\frakX_n'} \frakX_n$ be the morphism induced by two copies of $p_1$ and $p_1'.$ Since $p_1$ and $p_1'$ are strict, $g$ is also strict and so, by \ref{strictA}, we have an isomorphism of $\upmu$-indexed algebras $g^*\calA_{\frakX_n \times_{\frakX_n'} \frakX_n} \isom \calA_{P_n \times_{P_n'} P_n}.$ By \ref{thmkeyiso}, we obtain canonical isomorphisms
$$
g^*\calA_{\frakX_n \times_{\frakX_n'} \frakX_n} \isom g^* \left ( \calA_{\frakX_n} \circledast_{\calB_{F_n}} \calA_{\frakX_n} \right ) \isom p_1^*\calA_{\frakX_n} \circledast_{p_1'^*\calB_{F_n}} p_1^*\calA_{\frakX_n} \isom \calP_n^{\calA} \circledast_{\calP_n'^{\calB}} \calP_n^{\calA}.
$$ 
\end{proof}

\begin{proposition}\label{prop535}
Consider the morphism $\varphi:\calP_{n,\upmu}' \ra \calP_{n,\upmu}$ induced by \ref{prop433}. Then $\epsilon_{\calA_{\frakX_n}}$ is $\varphi^*\epsilon'_{\calB_{\frakX_n/\frakS_n}}$-linear.
\end{proposition}

\begin{proof}
Let $j:F_n^*\calB_{\frakX_n/\frakS_n} \ra \calA_{\frakX_n}$ be the canonical morphism. We have to prove that the commutativity of the diagram
$$
\begin{tikzcd}
\calP_{n,\upmu} \otimes_{\Ox_{\frakX_n,\upmu}} F_n^*\calB_{\frakX_n/\frakS_n} \ar{r}{\op{Id}\otimes j} \ar[swap]{d}{\varphi^*\epsilon_{\calB_{\frakX_n/\frakS_n}}} & \calP_{n,\upmu} \otimes_{\Ox_{\frakX_n,\upmu}} \calA_{\frakX_n} \ar{d}{\epsilon_{\calA_{\frakX_n}}} \\
F_n^*\calB_{\frakX_n/\frakS_n} \otimes_{\Ox_{\frakX_n,\upmu}} \calP_{n,\upmu} \ar{r}{j\otimes \op{Id}} & \calA_{\frakX_n} \otimes_{\Ox_{\frakX_n,\upmu}} \calP_{n,\upmu}.
\end{tikzcd}
$$
For that, it is sufficient to prove that $j$ is horizontal with respect to $d_{\calA_{\frakX_n}}$ and the connection $\nabla$ corresponding to $\varphi^*\epsilon'_{\calB_{\frakX_n/\frakS_n}}.$ In other words, we have to check the commutativity of
$$
\begin{tikzcd}
F_n^*\calB_{\frakX_n/\frakS_n} \ar{r}{j} \ar[swap]{d}{\nabla} & \calA_{\frakX_n} \ar{d}{d_{\calA_{\frakX_n}}} \\
(F_n^*\calB_{\frakX_n/\frakS_n}) \otimes_{\Ox_{\frakX_n,\upmu}} \omega^1_{\frakX_n/\frakS_n,\upmu} \ar{r}{j \otimes \op{Id}} & \calA_{\frakX_n} \otimes_{\Ox_{\frakX_n,\upmu}} \omega^1_{\frakX_n/\frakS_n,\upmu}.
\end{tikzcd}
$$
We can work locally and assume we have the same hypothesis and notation of \ref{loccoord}. We denote by $\xi_i$ and $\eta_i$ the reductions of $\frac{\widetilde{\eta}_i}{p}$ and $\widetilde{\eta}_i$ modulo $p^n$ i.e. their images in $\calP_n'$ and $\calP_n$ respectively, so that we have
$$
\calP_n'=\Ox_{\frakX_n'} \langle \xi_1,\hdots,\xi_d \rangle,\quad \calP_n = \Ox_{\frakX_n} \langle \eta_1,\hdots,\eta_d \rangle.
$$
Then, for all local sections $x$ of $\calB_{\frakX_n/\frakS_n},$ we have $\epsilon'_{\calB_{\frakX_n/\frakS_n}}(1\otimes x)=\sum_{I\in \N^d} x_I \otimes \xi^{[I]},$ where $x_I=\theta_{F_n}\left ( \partial_I' \right )(x)$ for all $I\in \N^d$ (where $\theta_{F_n}\left ( \partial_I' \right )$ is defined using the dual basis $\left (\partial_{\epsilon_1}',\hdots,\partial_{\epsilon_d}' \right )$ of $\left (\op{dlog}m_1',\hdots,\op{dlog}m_d' \right )$ as in \ref{def227}). Then
$$
\varphi^*\epsilon'_{\calB_{\frakX_n/\frakS_n}}(1\otimes x)=\sum_{I\in \N^d} F_n(x_I) \otimes \varphi \left ( \xi^{[I]} \right )=\sum_{J\in \N^d}y_J\otimes \eta^{[J]}.
$$
By \ref{equivstrat}, $\nabla(x)=\sum_{i=1}^dy_{\epsilon_i}\otimes \op{dlog}m_i.$
By \eqref{eq4243**}, we have
$$
y_{\epsilon_i} = F_n(x_{\epsilon_i})+\sum_{j=1}^d \frac{\partial_{\epsilon_i}(db_j)}{\alpha_{\frakX_n}(u_j)}F_n(x_{\epsilon_j}).
$$
Let $m'$ be a local section of $\calM_{\frakX_n'}$ and the unique local sections $a_1',\hdots,a_d'$ of $\Ox_{\frakX_n'}$ such that
$$\op{dlog}m'=\sum_{i=1}^da_i'\op{dlog}m_i'.$$
Then $\theta_{F_n} \left (\partial_{\epsilon_i}' \right )\left (e_{F_n^{\flat}m'} \right )=pa_i'e_{F_n^{\flat}m'}.$ It follows that, for $x=e_{F_n^{\flat}m'},$
\begin{alignat*}{2}
\frac{\partial_{\epsilon_i}(db_j)}{\alpha_{\frakX_n}(u_j)} F_n(x_{\epsilon_j}) &= F_n(a_j')\frac{\partial_{\epsilon_i} \left (d(1+pb_j) \right )}{\alpha_{\frakX_n}(u_j)}e_{F_n^{\flat}m'} = F_n(a_j')\frac{\partial_{\epsilon_i}\left (d\alpha_{\frakX_n}(u_j) \right )}{\alpha_{\frakX_n}(u_j)}e_{F_n^{\flat}m'}\\
&= F_n(a_j')\partial_{\epsilon_i}\left (\op{dlog}u_j \right )e_{F_n^{\flat}m'}.
\end{alignat*}
Then,
\begin{alignat*}{2}
\nabla( e_{F_n^{\flat}m'} ) &= \sum_{i=1}^dy_{\epsilon_i}\otimes \op{dlog}m_i = e_{F_n^{\flat}m'}\otimes \sum_{i=1}^d \left ( pF_n(a_i')+ \sum_{j=1}^d F_n(a_j')\partial_{\epsilon_i}(\op{dlog}u_j ) \right )  \op{dlog}m_i \\
&= e_{F_n^{\flat}m'} \otimes \sum_{i=1}^d F_n(a_i') \left ( p\op{dlog}m_i+\op{dlog}u_i \right ) =  e_{F_n^{\flat}m'} \otimes \op{dlog}F_n^{\flat}m' = d_{\calA_{\frakX_n}}(e_{F_n^{\flat}m'}).
\end{alignat*}
The equality for $\alpha_{\frakX_n}(m)e_{-m}$ is clear and the equality over $\Ox_{\frakX_n'}$ can be checked similarly.
\end{proof}

\begin{definition}\label{defadmstrat}
Denote by $\sigma:\upmu^2\ra \upmu$ the addition map.
\begin{enumerate}
\item Let $(\calE,\epsilon)$ be a $\upmu$-indexed $\calA_{\frakX_n}$-module equipped with a $\calP_{n,\upmu}$-stratification $\epsilon$ and denote by $\pi:\calA_{\frakX_n}\boxtimes \calE \ra \sigma^*\calE$ the multiplication map. We say that $\epsilon$ is \emph{admissible} if the following diagram is commutative
\begin{equation}\label{admindstratdiag}
\begin{tikzcd}
\calP_{n,\upmu^2}\otimes_{\Ox_{\frakX_n,\upmu^2}}(\calA_{\frakX_n}\boxtimes \calE) \ar{d}{\epsilon_{\calA_{\frakX_n}}\boxtimes \op{Id}_{\calE}} \ar{rr}{\op{Id}_{\calP_{n,\upmu^2}}\otimes \pi} & & \calP_{n,\upmu^2}\otimes_{\Ox_{\frakX_n,\upmu^2}} \sigma^*\calE \ar{dd}{\sigma^*\epsilon} \\
\calA_{\frakX_n} \boxtimes (\calP_{n,\upmu}\otimes_{\Ox_{\frakX_n,\upmu}} \calE) \ar{d}{\op{Id}_{\calA_{\frakX_n}}\boxtimes \epsilon} & & \\
(\calA_{\frakX_n}\boxtimes \calE)\otimes_{\Ox_{\frakX_n,\upmu^2}}\calP_{n,\upmu^2}\ar{rr}{\pi\otimes \op{Id}_{\calP_{n,\upmu^2}}} & & \sigma^*\calE\otimes_{\Ox_{\frakX_n,\upmu^2}}\calP_{n,\upmu^2}
\end{tikzcd}
\end{equation}
\item Let $(\calE,\epsilon)$ be a $\upmu$-indexed $\calB_{\frakX_n/\frakS_n}$-module equipped with a $\calP'_{n\upmu}$-stratification. We say that $\epsilon$ is \emph{admissible} if the diagram \eqref{admindstratdiag} is commutative when we replace $\calP_n,$ $\frakX_n,$ $\calA_{\frakX_n}$ and $\epsilon_{\calA_{\frakX_n}}$ by $\calP'_n,$ $\frakX_n',$ $\calB_{\frakX_n/\frakS_n}$ and $\epsilon'_{\calB_{\frakX_n/\frakS_n}}$ respectively.
\end{enumerate}
\end{definition}

\begin{remark}\label{rem126}
\begin{enumerate}
\item Let $(\calE,\epsilon)$ be a $\upmu$-indexed $\calA_{\frakX_n}$-module equipped with a $\calP_{n,\upmu}$-stratification. Then, by (\ref{PB}) and (\cite{SF3} 3.26), the isomorphism $\epsilon$ corresponds to a $\calP_{n,\upmu}$-linear isomorphism
$$\epsilon:\calP^{\calA}_n\circledast_{\calA_{\frakX_n}}\calE \isom\calE\circledast_{\calA_{\frakX_n}}\prescript{\calA}{}\calP_n.$$
Then, by definition, $\epsilon$ is admissible if and only if $\epsilon$ is $\epsilon_{\calA_{\frakX_n}}$-linear. Indeed, if $s,t\in \Gamma(U,\calM_{\frakX_n}^{\mathrm{gp}})$ are sections over an étale $\frakX_n$-scheme $U$ and $\ov{s},\ov{t}$ are their images in $\upmu,$ then the fiber of (\ref{admindstratdiag}) at $(\ov{s},\ov{t})$ is equal to
\begin{equation}\label{rem126diag}
\begin{tikzcd}
\calP_{n|U}\otimes_{\Ox_U}\calA_{\frakX_n,\ov{s}}\otimes_{\Ox_U} \calE_{\ov{t}} \ar{d}{\epsilon_{\calA_{\frakX_n},\ov{s}}\otimes\op{Id}_{\calE_{\ov{t}}}} \ar{rr}{\op{Id}_{\calP_{n|U}}\otimes \pi} & & \calP_{n|U}\otimes_{\Ox_{U}} \calE_{\ov{s}+\ov{t}} \ar{dd}{\epsilon_{\ov{s}+\ov{t}}} \\
\calA_{\frakX_n,\ov{s}} \otimes_{\Ox_U} \calP_{n|U}\otimes_{\Ox_{U}} \calE_{\ov{t}} \ar{d}{\op{Id}_{\calA_{\frakX_n,\ov{s}}}\otimes \epsilon_{\ov{t}}} & & \\
\calA_{\frakX_n,\ov{s}}\otimes_{\Ox_U} \calE_{\ov{t}}\otimes_{\Ox_{U}}\calP_{n|U}\ar{rr}{\pi\otimes \op{Id}_{\calP_{n|U}}} & & \calE_{\ov{s}+\ov{t}}\otimes_{\Ox_{U}}\calP_{n|U}
\end{tikzcd}
\end{equation}
and so the admissibility condition is equivalent to
$$\epsilon_{\ov{s}+\ov{t}}(1\otimes (e_s x))=\epsilon_{\calA_{\frakX_n},\ov{s}}(1\otimes e_s)\epsilon_{\ov{t}}(1\otimes x)$$
for all $x\in \Gamma(U,\calE_{\ov{t}}).$
\item Similarly, let $(\calE',\epsilon')$ be a $\upmu$-indexed $\calB_{\frakX_n/\frakS_n}$-module equipped with a $\calP_{n\upmu}'$-stratification
$$\epsilon':\calP'_{n,\upmu}\otimes_{\Ox_{\frakX_n',\upmu}}\calE' \isom \calE'\otimes_{\Ox_{\frakX_n',\upmu}}\calP'_{n,\upmu}.$$
Then, by (\ref{PB}) and (\cite{SF3} 3.26), the isomorphism $\epsilon'$ corresponds to a $\calP'_{n,\upmu}$-linear isomorphism
$$\epsilon':\calP'^{\calB}_n\circledast_{\calB_{\frakX_n/\frakS_n}}\calE' \isom\calE'\circledast_{\calB_{\frakX_n/\frakS_n}} \prescript{\calB}{}\calP'_n.$$
By definition, $\epsilon'$ is admissible if and only if it is $\epsilon'_{\calB_{\frakX_n/\frakS_n}}$-linear.
\end{enumerate}
\end{remark}

\begin{parag}
We denote by $\boldsymbol{\op{MHS}}(\calP_{n,\upmu}')$ the category of $\upmu$-indexed $\Ox_{\frakX_n',\upmu}$-modules $\calE'$ equipped with a $\calP_{n,\upmu}'$-stratification and by $\boldsymbol{\op{MHS}}(\calB_{\frakX_n/\frakS_n})$ the subcategory of $\boldsymbol{\op{MHS}}(\calP_{n,\upmu}')$ consisting of $\upmu$-indexed $\calB_{\frakX_n/\frakS_n}$-modules $\calE'$ equipped with an admissible $\calP_{n,\upmu}'$-stratification (\ref{defadmstrat}) and $\calB_{\frakX_n/\frakS_n}$-linear morphisms.
We also denote by $\boldsymbol{\op{MHS}}(\calP_{n,\upmu})$ the category $\Ox_{\frakX_n,\upmu}$-modules $\calE$ equipped with a $\calP_{n,\upmu}$-stratification and by $\boldsymbol{\op{MHS}}(\calA_{\frakX_n})$ the subcategory consisting of $\upmu$-indexed $\calA_{\frakX_n}$-modules $\calE$ equipped with an admissible $\calP_{n,\upmu}$-stratification (\ref{defadmstrat}) and $\calA_{\frakX_n}$-linear morphisms. For what follows, we recall the following notation: let $\nabla$ be an integrable connection on an $\Ox_{\frakX_n,\upmu}$-module $\calE$ and suppose we have local coordinates $m_1,\hdots,m_d$ for $\frakX_n/\frakS_n.$ Let $\left (\partial_{\epsilon_1},\hdots,\partial_{\epsilon_d} \right )$ be the dual basis of $\left (\op{dlog}m_1,\hdots,\op{dlog}m_d \right ).$ Recall the endomorphisms \eqref{def227}
\begin{equation}\label{eq4371}
\nabla \left ( \partial_{\epsilon_i} \right ) : \calE \xrightarrow{\nabla} \calE \otimes_{\Ox_{\frakX_n,\upmu}} \omega^1_{\frakX_n/\frakS_n,\upmu} \xrightarrow{\op{Id}_{\calE} \otimes \partial_{\epsilon_i}} \calE\quad \forall 1\le i\le d,\qquad 
\nabla \left ( \partial_I \right ) = \prod_{i=1}^d \prod_{j=0}^{I_i-1} \left ( \nabla \left (\partial_{\epsilon_i} \right )-j \right )\quad \forall I\in \N^d.
\end{equation}
\end{parag}

\begin{proposition}\label{prop12.9}
Let $\calE$ be a $\upmu$-indexed $\calA_{\frakX_n}$-module (resp. $\upmu$-indexed $\calB_{\frakX_n/\frakS_n}$-module) equipped with a $\calP_{n,\upmu}$-stratification (resp. $\calP'_{n,\upmu}$-stratification)
$$\epsilon:\calP_{n,\upmu}\otimes_{\Ox_{\frakX_n,\upmu}}\calE \isom\calE\otimes_{\Ox_{\frakX_n,\upmu}}\calP_{n,\upmu}\quad \left (\text{resp.}\ \calP'_{n,\upmu}\otimes_{\Ox_{\frakX'_n,\upmu}}\calE \isom\calE\otimes_{\Ox_{\frakX'_n,\upmu}}\calP'_{n,\upmu}\right ).$$
Denote by $\nabla:\calE \ra \calE \otimes_{\Ox_{\frakX_n,\upmu}}\omega^1_{\frakX_n/\frakS_n,\upmu}$ (resp. $\nabla:\calE \ra \calE \otimes_{\Ox_{\frakX'_n,\upmu}}\omega^1_{\frakX'_n/\frakS_n,\upmu}$) the connection (resp. $p$-connection) corresponding to $\epsilon.$ If $\nabla$ is admissible then so is $\epsilon.$
\end{proposition}

\begin{proof}
We just prove the assertion for connection as the other case is similar. Let $s,t\in \Gamma(U,\calM_{\frakX_n}^{\mathrm{gp}})$ and $x\in \Gamma(U,\calE_{\ov{s}})$ be sections over an étale $\frakX_n$-scheme $U.$ We want to prove that
$$\epsilon_{\ov{s}+\ov{t}}(1\otimes(e_tx))=\epsilon_{\calA_{\frakX_n},\ov{t}}(1\otimes e_t)\epsilon_{\ov{s}}(1\otimes x).$$
By \ref{thmKato}, we have
\begin{alignat}{2}
\epsilon_{\calA_{\frakX_n},\ov{t}}(1\otimes e_t) &= \sum_{I \in \mathbb{N}^d}d_{\calA_{\frakX_n},\ov{t}}(\partial_I)(e_t)\otimes \eta^{[I]},\quad 
\epsilon_{\ov{s}}(1\otimes x) = \sum_{J\in \mathbb{N}^d} \nabla_{\ov{s}}(\partial_J)(x)\otimes \eta^{[J]}.
\end{alignat}
So
$$\epsilon_{\calA_{\frakX_n},\ov{t}}(1\otimes e_t)\epsilon_{\ov{s}}(1\otimes x)=\sum_{K\in \mathbb{N}^d}\left (\sum_{I+J=K}\begin{pmatrix} K\\ I\end{pmatrix} d_{\calA_{\frakX_n},\ov{t}}(\partial_I)(e_t)\nabla_{\ov{s}}(\partial_J)(x) \right )\otimes \eta^{[K]}.$$
It is thus sufficient to prove that, for any multi-index $K\in \mathbb{N}^d,$
\begin{equation}\label{diffeq0}
\nabla_{\ov{s}+\ov{t}}(\partial_K)(e_tx)=\sum_{I+J=K} \begin{pmatrix} K\\ I\end{pmatrix} d_{\calA_{\frakX_n},\ov{t}}(\partial_I)(e_t)\nabla_{\ov{s}}(\partial_J)(x).
\end{equation}
We prove this in the next lemma.
\end{proof}

\begin{lemma}
Keep the same hypothesis of the proposition \ref{prop12.9} and its proof. Then \eqref{diffeq0} holds.
\end{lemma}

\begin{proof}
To simplify notations, we drop the subscripts $\ov{s}$ and $\ov{t}.$ 
We prove the lemma by induction on $|K|:$ for $|K|=0,$ the equality is obvious. Fix an integer $n\ge 0$ and suppose \eqref{diffeq0} holds for all $K\in\mathbb{N}^d$ such that $|K|\le n.$ Let $1\le i\le d,$ $K\in\mathbb{N}^d$ such that $|K|=n$ and let us prove that \eqref{diffeq0} holds for $K+\epsilon_i.$
By the definition \eqref{eq4371}, we have $\nabla \left ( \partial_{K+\epsilon_i} \right ) = \nabla \left ( \partial_{\epsilon_i} \right )\circ \nabla \left ( \partial_K \right )  - K_i \nabla \left ( \partial_K \right )$
and, for local sections $a$ of $\calA_{\frakX_n},$ the admissibility of $\nabla$ yields
$$
\nabla(\partial_{\epsilon_i})(ax)=a\nabla(\partial_{\epsilon_i})(x)+d_{\calA_{\frakX_n}}\left ( \partial_{\epsilon_i} \right )(a)x.
$$
It follows that
\begin{alignat*}{2}
 & \nabla \left (\partial_{\epsilon_i} \right )(d_{\calA_{\frakX_n}} \left (\partial_{I} \right )(e_t) \nabla \left (\partial_{J} \right ) (x)) \\
=& d_{\calA_{\frakX_n}} \left (\partial_{I} \right )(e_t) \nabla \left (\partial_{\epsilon_i} \right ) \circ \nabla \left (\partial_{J} \right ) (x) + \left (d_{\calA_{\frakX_n}} \left (\partial_{\epsilon_i} \right ) \circ d_{\calA_{\frakX_n}} \left (\partial_{I} \right )(e_t) \right ) \nabla \left (\partial_{J} \right ) (x) \\
=& d_{\calA_{\frakX_n}} \left (\partial_{I} \right )(e_t)  \left ( \nabla \left (\partial_{J+\epsilon_i} \right ) + J_i \nabla \left (\partial_{J} \right ) \right )(x) + \left (d_{\calA_{\frakX_n}} \left (\partial_{I+\epsilon_i} \right ) +I_i d_{\calA_{\frakX_n}} \left (\partial_{I} \right )(e_t) \right ) \nabla \left (\partial_{J} \right ) (x).
\end{alignat*}
If $I+J=K,$ we get
\begin{alignat*}{2}
&\nabla \left (\partial_{\epsilon_i} \right )(d_{\calA_{\frakX_n}} \left (\partial_{I} \right )(e_t) \nabla \left (\partial_{J} \right ) (x)) \\
=& d_{\calA_{\frakX_n}} \left (\partial_{I} \right )(e_t)  \nabla \left (\partial_{J+\epsilon_i} \right )(x) + d_{\calA_{\frakX_n}} \left (\partial_{I+ \epsilon_i} \right )(e_t)  \nabla \left (\partial_{J} \right )(x) + K_i d_{\calA_{\frakX_n}} \left (\partial_{I} \right )(e_t)  \nabla \left (\partial_{J} \right )(x).
\end{alignat*}
This equality and the induction hypothesis yield
\begin{alignat*}{2}
 & \nabla \left (\partial_{K+\epsilon_i} \right )(e_tx) = \left ( \nabla \left (\partial_{\epsilon_i} \right ) -K_i \right ) \circ \nabla \left (\partial_{K} \right )(e_tx) \\
=& \sum_{I+J=K} \begin{pmatrix} K\\ I\end{pmatrix} \left (\nabla \left (\partial_{\epsilon_i} \right )(d_{\calA_{\frakX_n}} \left (\partial_{I} \right )(e_t) \nabla \left (\partial_{J} \right ) (x)) - K_i d_{\calA_{\frakX_n}} \left (\partial_{I} \right )(e_t)  \nabla \left (\partial_{J} \right )(x) \right ) \\
=& \sum_{I+J=K} \begin{pmatrix} K\\ I\end{pmatrix} \left ( d_{\calA_{\frakX_n}} \left (\partial_{I} \right )(e_t)  \nabla \left (\partial_{J+\epsilon_i} \right )(x) + d_{\calA_{\frakX_n}} \left (\partial_{I+ \epsilon_i} \right )(e_t)  \nabla \left (\partial_{J} \right )(x)  \right ) \\
=&  \sum_{I'+J'=K+\epsilon_i} \left ( \begin{pmatrix} K\\ I'-\epsilon_i \end{pmatrix} + \begin{pmatrix} K\\ I'\end{pmatrix} \right ) \left ( d_{\calA_{\frakX_n}} \left (\partial_{I'} \right )(e_t)  \nabla \left (\partial_{J'} \right )(x) \right ),
\end{alignat*}
with the convention $\begin{pmatrix} K\\ -\epsilon_i\end{pmatrix}=\begin{pmatrix} K\\ K+\epsilon_i\end{pmatrix}=0$ in the last line. We conclude by $$\begin{pmatrix} K\\ I'-\epsilon_i \end{pmatrix} + \begin{pmatrix} K\\ I'\end{pmatrix} = \begin{pmatrix} K+\epsilon_i\\ I'\end{pmatrix}.$$
\end{proof}

\begin{proposition}\label{indpropstrat2}
The equivalences of categories given in \ref{prop12.3}, induce equivalences
$$\boldsymbol{\op{MHS}}(\calA_{\frakX_n}) \isom \boldsymbol{\op{MIC}}^{\mathrm{qn}}(\calA_{\frakX_n}),\quad \boldsymbol{\op{MHS}}(\calB_{\frakX_n/\frakS_n}) \isom p\text{-}\boldsymbol{\op{MIC}}^{\mathrm{qn}}(\calB_{\frakX_n/\frakS_n}).$$
\end{proposition}

\begin{proof}
Both equivalences are deduced in a similar way, so we just prove the one about connections.
For any $\Ox_{\frakX_n,\upmu}$-module $\calE$ and any local section $m\in\Gamma(U,\upmu)$ over an étale $\frakX_n$-scheme $U,$ let
$$\pi_{m}:\calE_{m}\otimes_{\Ox_{U}}\calP_{n|U} \ra \calE_m\otimes_{\Ox_U} \left ( \calP_{n|U}/\ov{\calI}^{[2]} \right )$$
be the canonical projection, where $\ov{\calI}$ is the PD ideal of $P_n.$ 
Let $(\calE,\epsilon)$ be an object of $\boldsymbol{\op{MHS} \left (\calP_{n,\upmu} \right )}$ and $(\calE,\nabla)$ its image in $\boldsymbol{\op{MIC}}^{\mathrm{qn}}(\frakX_n/\frakS_n).$

Suppose that $(\calE,\epsilon)$ is an object of the subcategory $\boldsymbol{\op{MHS} \left ( \calA_{\frakX_n} \right )}.$ For all sections $s,t\in \Gamma(U,\calM^{\mathrm{gp}}_{\frakX_n})$ over an étale $\frakX_n$-scheme, the morphism
$\nabla_{\ov{s}+\ov{t}}:\calE_{\ov{s}+\ov{t}}\ra \calE_{\ov{s}+\ov{t}}\otimes_{\Ox_U}\omega^1_{U/\frakS_n}$ is given by
$$\nabla_{\ov{s}+\ov{t}}(e_sx)=\pi_{\ov{s}+\ov{t}}(\epsilon_{\ov{s}+\ov{t}}(1\otimes(e_sx))-(e_sx)\otimes 1),$$
for all local sections $x\in \Gamma(U,\calE_{\ov{t}}).$
Since $\epsilon$ is admissible,
\begin{equation}\label{diffeq10}
\epsilon_{\ov{s}+\ov{t}}(1\otimes(e_sx))=\epsilon_{\calA_{\frakX_n},\ov{s}}(1\otimes e_s)\epsilon_{\ov{t}}(1\otimes x).
\end{equation}
By the usual equivalence between integrable quasi-nilpotent connections and $\calP_n$-stratifications (\ref{thmKato}), we have
\begin{alignat*}{2}
\pi_{\ov{s}}(\epsilon_{\calA_{\frakX_n},\ov{s}}(1\otimes e_s)) &= d_{\calA_{\frakX_n},\ov{s}}(e_s)+e_s\otimes 1,\quad
\pi_{\ov{t}}(\epsilon_{\ov{t}}(1\otimes x)) = \nabla_{\ov{t}}(x)+x\otimes 1. 
\end{alignat*}
so, in $\calE_{\ov{s}+\ov{t}} \otimes_{\Ox_U} \left ( \calP_{n|U}/\ov{\calI}^{[2]} \right ),$ we get
\begin{alignat*}{2}
\pi_{\ov{s}+\ov{t}}(\epsilon_{\ov{s}+\ov{t}}(1\otimes e_sx)) &= (d_{\calA_{\frakX_n},\ov{s}}(e_s)+e_s\otimes 1)(\nabla_{\ov{t}}(x)+x\otimes 1)\\
&= (e_sx)\otimes 1+e_s\nabla_{\ov{t}}(x)+x\otimes d_{\calA_{\frakX_n},\ov{s}}(e_s).
\end{alignat*}
It follows that
$$\nabla_{\ov{s}+\ov{t}}(e_sx)=e_s\nabla_{\ov{t}}(x)+x\otimes d_{\calA_{\frakX_n},\ov{s}}(e_s)$$
and $\nabla$ is hence admissible.

Conversely, if $(\calE,\nabla)$ is an object of the subcategory $\boldsymbol{\op{MIC}}^{\mathrm{qn}}(\calA_{\frakX_n}),$ then the fact that $\epsilon$ is admissible follows from \ref{prop12.9}.
\end{proof}

\subsection{Stratified Shiho functor and the Cartier transform}
\begin{parag}
Let $(\calE',\epsilon')$ be an object of $\boldsymbol{\op{MHS}}(\calB_{\frakX_n/\frakS_n})$ and consider the morphism $\varphi:P_n \ra P'_n$
deduced from the morphism $P_{\frakX/\frakS,0} \ra P_{\frakX'/\frakS,1},$ given in \ref{prop433}, by reduction modulo $p^n.$ We abusively denote the morphism of structure rings of $\varphi$ by
$\varphi:\calP'_n \ra \calP_n.$
By \ref{rem126}, the admissible $\calP_{n,\upmu}'$-stratification $\epsilon'$ is an $\epsilon_{\calB_{\frakX_n/\frakS_n}}'$-linear isomorphism
\begin{equation}\label{eq4411}
\epsilon':\calP_n'^{\calB}\circledast_{\calB_{\frakX_n/\frakS_n}}\calE' \ra \calE' \circledast_{\calB_{\frakX_n/\frakS_n}}\prescript{\calB}{}{\calP_n'}.
\end{equation}
Similarly, by \ref{prop535}, the canonical stratification $\epsilon_{\calA_{\frakX_n/\frakS_n}}$ corresponds to an $\epsilon_{\calB_{\frakX_n/\frakS_n}}'$-linear isomorphism
$\epsilon_{\calA_{\frakX_n}}:\calP_n^{\calA}\ra \prescript{\calA}{}{\calP_n}.$
The morphism $\varphi$ and the canonical embedding $\calB_{\frakX_n/\frakS_n}\hookrightarrow \calA_{\frakX_n}$ induce a morphism
$
\calP_n'^{\calB} \ra \calP_n^{\calA}
$
which turns $\calP_n^{\calA}$ into a $\upmu$-indexed $\calP_n'^{\calB}$-algebra.
Let $\calE=\calA_{\frakX_n} \circledast_{\calB_{\frakX_n/\frakS_n}} \calE'$ and
$$\epsilon=\epsilon' \circledast_{\calP_n'^{\calB}} \epsilon_{\calA_{\frakX_n}}:\calP_{n,\upmu}\otimes_{\Ox_{\frakX_n,\upmu}}\calE \ra \calE \otimes_{\Ox_{\frakX_n,\upmu}}\calP_{n,\upmu}.$$
Then $\epsilon$ is an admissible $\calP_{n,\upmu}$-stratification making $(\calE,\epsilon)$ an object of $\boldsymbol{\op{MHS}}(\calA_{\frakX_n})$ and we obtain a functor
\begin{equation}\label{indpsi}
\Psi:\begin{array}[t]{clc}
\boldsymbol{\op{MHS}}(\calB_{\frakX_n/\frakS_n}) & \ra & \boldsymbol{\op{MHS}}(\calA_{\frakX_n}) \\
(\calE',\epsilon') & \mapsto & (\calE,\epsilon),
\end{array}
\end{equation}
where $\calE$ and $\epsilon$ are defined above.
\end{parag}

\begin{lemma}\label{lemPD}
Consider the morphism $\varphi_i$ given in \ref{philog}. The ideals of the closed immersions $\frakX \ra \frakX \times_{\frakX'} \frakX$ and $\frakX \ra P_{\frakX/\frakS,0}\times_{P_{\frakX'/\frakS,1}} P_{\frakX/\frakS,0}$ have unique PD structures.
\end{lemma}

\begin{proof}
The morphisms $F:\frakX \ra \frakX'$ and $\varphi_i:P_{\frakX/\frakS,0} \ra P_{\frakX'/\frakS,1}$ are log flat by \ref{logflatlift} and \ref{faithfulflatphi}. It follows that the projections $\frakX \times_{\frakX'} \frakX \ra \frakX$ and $P_{\frakX/\frakS,0} \times_{P_{\frakX'/\frakS,1}} P_{\frakX/\frakS,0} \ra P_{\frakX/\frakS,0}$ are also log flat. Since they are strict, they are flat and $\frakX \times_{\frakX'} \frakX$ and $P_{\frakX/\frakS,0}\times_{P_{\frakX'/\frakS,1}} P_{\frakX/\frakS,0}$ are hence flat over $W.$ It is thus sufficient to prove that for all local sections $x$ of the ideal of $\frakX \ra \frakX \times_{\frakX'} \frakX$ (resp. $\frakX \ra P_{\frakX/\frakS,0}\times_{P_{\frakX'/\frakS,1}} P_{\frakX/\frakS,0}$) and all integers $k\ge 1,$ $x^k\in k!\Ox_{\frakX \times_{\frakX'} \frakX}$ (resp. $x^k\in k!\Ox_{P_{\frakX/\frakS,0}\times_{P_{\frakX'/\frakS,1}} P_{\frakX/\frakS,0}}$). The rest of the proof is similar to the proofs of (\cite{SF1} 12.6) and (\cite{Shiho} 3.7).
\end{proof}

\begin{theorem}\label{Psiequiv}
The functor $\Psi$ \eqref{indpsi} is an equivalence of categories.
\end{theorem}

\begin{proof}
Consider the PD structure on the ideal of the diagonal immersion $\frakX_n \ra \frakX_n\times_{\frakX_n'}\frakX_n,$ induced by \ref{lemPD}. This immersion is strict by \ref{prop114}. It follows that, by the universal property of the logarithmic PD envelope (\cite{Kat89} 5.3), the canonical morphism $\frakX_n \times_{\frakX_n'} \frakX_n \ra \frakX_n \times_{\frakS_n} \frakX_n$ factors uniquely through $P_n.$ Denote by $\psi:\frakX_n \times_{\frakX_n'} \frakX_n \ra P_n$ the resulting morphism. By a similar argument to the proof of (\cite{SF1} 12.9), the composition $\varphi \circ \psi$ factors through the canonical immersion $\iota':\frakX_n' \ra P_n'.$ We hence have the commutative diagram
\begin{equation}
\begin{tikzcd}
 & & & \frakX_n'\ar{d}{\iota'} & \\
 & & P_n \ar{r}{\varphi} \ar{d} & P_n'\ar{d} \ar[shift right=0.5ex]{dr}{p_1'} \ar[shift left=0.5ex,swap]{dr}{p_2'} & \\
\frakX_n & \frakX_n\times_{\frakX_n'}\frakX_n \ar{ur}{\psi} \ar[bend right=-30]{uurr}{\nu} \ar{r} \ar[shift right=0.5ex,swap]{l}{q_1} \ar[shift left=0.5ex]{l}{q_2} & \frakX_n\times_{\frakS_n}\frakX_n \ar{r}{F_n\times_{\frakS_n}F_n} & \frakX_n'\times_{\frakS_n}\frakX_n' \ar[shift left=0.5ex]{r} \ar[shift right=0.5ex]{r} & \frakX_n',
\end{tikzcd}
\end{equation}
where $p_1',$ $p_2',$ $q_1$ and $q_2$ are the projections. Let us prove the fullness of $\Psi:$ let $(\calE_1',\epsilon_1')$ and $(\calE_2',\epsilon_2')$ be two objects of $\boldsymbol{\op{MHS}}\left ( \calB_{\frakX_n/\frakS_n} \right ).$ Let $(\calE_i,\epsilon_i)=\Psi(\calE_i',\epsilon_i')$ for $i=1,2.$ Let $g:\calE_1\ra \calE_2$ be an $\calA_{\frakX_n}$-linear morphism compatible with the stratifications. Note that
$$
 \psi^* \epsilon'_i =  \left (\varphi \circ \psi \right )^*\left ( \epsilon'_i \right )=\nu^*\iota'^*\epsilon'_i= \nu^*\op{Id}_{\calE_i'},
$$
$$
\psi^*\epsilon_i=\psi^*\left (\epsilon_{\calA_{\frakX_n}} \circledast \epsilon'_i \right ) = \psi^*\epsilon_{\calA_{\frakX_n}} \circledast \psi^* \epsilon'_i = \psi^*\epsilon_{\calA_{\frakX_n}} \circledast\nu^*\op{Id}_{\calE_i'}.
$$
By \ref{remdefeA}, $\psi^*\epsilon_{\calA_{\frakX_n}}$ is the isomorphism that exchanges factors.
Pulling back the following left diagram along $\psi,$ we get the commutative diagram on the right:
$$
\begin{tikzcd}
\calP_{n,\upmu} \otimes_{\Ox_{\frakX_n,\upmu}} \calE_1 \ar{r}{\op{Id} \otimes g} \ar[swap]{d}{\epsilon_1} & \calP_{n,\upmu} \otimes_{\Ox_{\frakX_n,\upmu}} \calE_2 \ar{d}{\epsilon_2} \\
\calE_1  \otimes_{\Ox_{\frakX_n,\upmu}} \calP_{n,\upmu} \ar{r}{g\otimes \op{Id}} & \calE_2  \otimes_{\Ox_{\frakX_n,\upmu}} \calP_{n,\upmu},
\end{tikzcd}\quad
\begin{tikzcd}
\Ox_{\frakX_n \times_{\frakX_n'} \frakX_n,\upmu} \otimes_{\Ox_{\frakX_n,\upmu}} \calE_1 \ar{r}{\op{Id} \otimes g} \ar[swap,sloped]{d}{\sim} & \Ox_{\frakX_n \times_{\frakX_n'} \frakX_n,\upmu} \otimes_{\Ox_{\frakX_n,\upmu}} \calE_2 \ar[sloped]{d}{\sim} \\
\calE_1  \otimes_{\Ox_{\frakX_n,\upmu}} \Ox_{\frakX_n \times_{\frakX_n'} \frakX_n,\upmu} \ar{r}{g\otimes \op{Id}} & \calE_2  \otimes_{\Ox_{\frakX_n,\upmu}} \Ox_{\frakX_n \times_{\frakX_n'} \frakX_n,\upmu},
\end{tikzcd}
$$
where the vertical isomorphisms of the right diagram are induced by the isomorphism that exchanges factors on $\frakX_n \times_{\frakX_n'} \frakX_n.$ By \ref{prop114}, the projections $q_1$ and $q_2$ are strict. Hence they induce isomorphisms of $\upmu$-indexed algebras $\Ox_{\frakX_n \times_{\frakX_n'} \frakX_n,\upmu} \otimes_{\Ox_{\frakX_n,\upmu}} \calA_{\frakX_n} \isom \calA_{\frakX_n\times_{\frakX_n'} \frakX_n}.$ By \ref{thmkeyiso}, we deduce an isomorphism of $\upmu$-indexed algebras
$$\Ox_{\frakX_n \times_{\frakX_n'} \frakX_n,\upmu} \otimes_{\Ox_{\frakX_n,\upmu}} \calA_{\frakX_n} \isom \calA_{\frakX_n} \circledast_{\calB_{\frakX_n/\frakS_n}} \calA_{\frakX_n}.$$
Since $g$ is $\calA_{\frakX_n}$-linear, the diagram to the right becomes equal to
$$
\begin{tikzcd}
\left ( \calA_{\frakX_n} \circledast_{\calB_{\frakX_n/\frakS_n}} \calA_{\frakX_n} \right ) \circledast_{\calA_{\frakX_n}} \calE_1 \ar{r}{\op{Id} \circledast g} \ar[swap,sloped]{d}{\sim} & \left ( \calA_{\frakX_n} \circledast_{\calB_{\frakX_n/\frakS_n}} \calA_{\frakX_n} \right ) \circledast_{\calA_{\frakX_n}} \calE_2 \ar[sloped]{d}{\sim} \\
\calE_1  \circledast_{\calA_{\frakX_n}} \left ( \calA_{\frakX_n} \circledast_{\calB_{\frakX_n/\frakS_n}} \calA_{\frakX_n} \right ) \ar{r}{g\circledast \op{Id}} & \calE_2  \circledast_{\calA_{\frakX_n}} \left ( \calA_{\frakX_n} \circledast_{\calB_{\frakX_n/\frakS_n}} \calA_{\frakX_n} \right ).
\end{tikzcd}
$$
By \ref{inddescentmor}, there exists a unique $\calB_{\frakX_n/\frakS_n}$-linear morphism $f:\calE_1' \ra \calE_2'$ whose base change by $\calB_{\frakX_n/\frakS_n} \ra \calA_{\frakX_n}$ is equal to $g.$ It remains to prove that $f$ is compatible with the stratifications $\epsilon_1'$ and $\epsilon_2'.$ Note that the base change of the diagram
\begin{equation}\label{diag4421}
\begin{tikzcd}
\calP_n'^{\calB} \circledast_{\calB_{\frakX_n/\frakS_n}} \calE_1' \ar{r}{\op{Id} \circledast f} \ar[swap]{d}{\epsilon_1'} & \calP_n'^{\calB} \circledast_{\calB_{\frakX_n/\frakS_n}} \calE_2' \ar{d}{\epsilon_2'} \\
\calE_1' \circledast_{\calB_{\frakX_n/\frakS_n}} \prescript{\calB}{}\calP_n' \ar{r}{f\circledast \op{Id}} &  \calE_2' \circledast_{\calB_{\frakX_n/\frakS_n}} \prescript{\calB}{}\calP_n'
\end{tikzcd}
\end{equation}
by $\calP_n'^{\calB} \ra \calP_n^{\calA}$ is commutative by compatibility of $g$ with the stratifications $\epsilon_1$ and $\epsilon_2.$ Since $\calP_n'^{\calB} \ra \calP_n^{\calA}$ is faithfully flat by \ref{faithfulflatphi}, we deduce that \eqref{diag4421} is commutative by \ref{inddescentmor}. We have hence proved the fullness of $\Psi.$ The faithfulness follows from \ref{inddescentmor} by a similar argument.

Let us now prove the essential surjectivity: let $(\calE,\epsilon)$ be an object of $\boldsymbol{\op{MHS}}\left ( \calA_{\frakX_n} \right ).$ Just as before, we deduce from $\epsilon$ an isomorphism $\alpha:\left ( \calA_{\frakX_n} \circledast_{\calB_{\frakX_n/\frakS_n}} \calA_{\frakX_n} \right ) \circledast_{\calA_{\frakX_n}} \calE \isom \calE \circledast_{\calA_{\frakX_n}} \left ( \calA_{\frakX_n} \circledast_{\calB_{\frakX_n/\frakS_n}} \calA_{\frakX_n} \right ).$
The cocyclicity property of $\epsilon$ implies that $\alpha$ satisfies the cocyclicity property given in \ref{inddescentmod}. It follows that there exists a $\upmu$-indexed $\calB_{\frakX_n/\frakS_n}$-module $\calE'$ such that $\calA_{\frakX_n} \circledast_{\calB_{\frakX_n/\frakS_n}} \calE'=\calE.$ All that remains is to descend the stratification $\epsilon.$ We identify $\calP_n^{\calA}$ with $\prescript{\calA}{}\calP_n$ via $\epsilon_{\calA_{\frakX_n}}$ and $\calP_n'^{\calB}$ with $\prescript{\calB}{}\calP'_n$ via $\epsilon'_{\calB_{\frakX_n/\frakS_n}}.$ Consider the PD structure on the ideal of $\frakX_n \ra P_n \times_{P_n'} P_n$ induced by \ref{lemPD}. For $1\le i<j\le 4,$ we have a unique $\pi_{ij}:P_n\times_{P_n'}P_n \ra P_n$ induced by the projection $\frakX_n\times_{\frakS_n} \frakX_n \times_{\frakS_n} \frakX_n \times_{\frakS_n} \frakX_n\ra \frakX_n \times_{\frakS_n}\frakX_n$ on the factors $i$ and $j.$ Similarly, for $1\le i\le 4,$ let $\pi_i:P_n \times_{P_n'} P_n \ra \frakX_n$ be the morphism induced by the projection on the $i$th factor. By \ref{lemAP}, $\pi_{ij}$ and $\pi_i$ induce morphisms of indexed algebras, abusively denoted
$$
\pi_{ij}:\calP_n^{\calA} \ra \calP_n^{\calA} \circledast_{\calP_n'^{\calB}} \calP_n^{\calA},\quad \pi_i:\calA_{\frakX_n} \ra \calP_n^{\calA} \circledast_{\calP_n'^{\calB}} \calP_n^{\calA}.
$$
Consider the descent datum on $\calP_n^{\calA} \circledast_{\calA} \calE$ (resp. $\calE \circledast_{\calA} \prescript{\calA}{}\calP_n$) given by the pullback of $\epsilon$ along the composition
$$
\begin{tikzcd}
P_n \times_{P_n'} P_n \ar{r}{p_2\times p_2} \ar{dr}{\pi_{24}} & \frakX_n \times_{\frakX_n'} \frakX_n \ar{d}{\psi} \\
 & P_n
\end{tikzcd}
\quad \left (\text{resp.}\
\begin{tikzcd}
P_n \times_{P_n'} P_n \ar{r}{p_1\times p_1} \ar{dr}{\pi_{13}} & \frakX_n \times_{\frakX_n'} \frakX_n \ar{d}{\psi} \\
 & P_n
\end{tikzcd}
\right ).
$$
To prove that $\epsilon$ descends to a $\calP_n'^{\calB}$-linear morphism $\epsilon':\calP_n'^{\calB} \circledast_{\calB_{\frakX_n / \frakS_n}} \calE' \ra \calE'  \circledast_{\calB_{\frakX_n / \frakS_n}} \prescript{\calB}{} \calP_n',$ using the faithful flatness of $\calP_n'^{\calB} \ra \calP_n^{\calA}$ and \ref{inddescentmor}, it is sufficient to prove the commutativity of
$$
\begin{tikzcd}
\pi_4^*\calE \ar{r}{\pi_{34}^*\epsilon} \ar[swap]{d}{\pi_{24}^*\epsilon} & \pi_3^*\calE \ar{d}{\pi_{13}^*\epsilon} \\
\pi_2^*\calE \ar{r}{\pi_{12}^*\epsilon} & \pi_1^*\calE.
\end{tikzcd}
$$
The cocyclicity property of the stratification $\epsilon$ implies that $\pi_{13}^*\epsilon \circ \pi_{34}^*\epsilon = \pi_{14}^*\epsilon = \pi_{12}^*\epsilon \circ \pi_{24}^*\epsilon.$ Hence the existence of $\epsilon'.$ The fact that $\epsilon'$ is an isomorphism is inherited from $\epsilon,$ by \ref{inddescentmor}. The cocyclicity property of $\epsilon'$ is inherited from the cocyclicity of $\epsilon,$ again by \ref{inddescentmor}, and by the faithful flatness of
$$
\calP_n'^{\calB} \circledast_{\calB_{\frakX_n/\frakS_n}} \calP_n'^{\calB} \ra \calP_n^{\calA} \circledast_{\calA_{\frakX_n}} \calP_n^{\calA},
$$
which results from the faithful flatness of \eqref{eq5260}.
\end{proof}

\begin{proposition}\label{prop443}
Consider the diagram
\begin{equation}\label{inddiagann}
\begin{tikzcd}
\boldsymbol{\op{MHS}(\calB_{\frakX_n/\frakS_n})} \ar{r}{\Psi} \ar{d} & \boldsymbol{\op{MHS}(\calA_{\frakX_n})}\ar{d} \\
p\text{-}\boldsymbol{\op{MIC}}^{\mathrm{qn}}(\calB_{\frakX_n/\frakS_n}) \ar[hook]{d} & \boldsymbol{\op{MIC}}^{\mathrm{qn}}(\calA_{\frakX_n})\ar[hook]{d} \\
\boldsymbol{p\op{\text{-}MIC}(\calB_{\frakX_n/\frakS_n})} \ar{r}{\Phi} & \boldsymbol{\op{MIC}(\calA_{\frakX_n})}
\end{tikzcd}
\end{equation}
where $\Psi$ is defined in \eqref{indpsi}, $\Phi$ is defined in \eqref{Phiind}, the upper vertical arrows are the equivalences of categories in \ref{indpropstrat2} and the lower vertical arrows are the canonical embeddings. Then \eqref{inddiagann} is commutative.
\end{proposition}

\begin{proof}
We can work locally on $\frakX$ and $\frakS$ and hence suppose that we have local coordinates $m_1,\hdots,m_d\in \Gamma(X,\calM_{X})$ of $X/S.$ We deduce local coordinates $m_1'=\pi^{\flat}m_1,\hdots,m_d'=\pi^{\flat}m_d$ of $X'/S,$ where $\pi:X'\ra X$ is given in \eqref{diagfrob}. We may assume that $m_i$ and $m_i'$ lift to sections $\tilde{m}_i$ and $\tilde{m}_i'$ of $\calM_{\frakX}$ and $\calM_{\frakX'}$ respectively. Define $\tilde{\eta}_i$ and $\tilde{\eta}_i'$ as in \ref{loccoord}. We denote by $\xi_i$ and $\eta_i$ respectively the reductions of $\frac{\tilde{\eta}_i'}{p}$ and $\widetilde{\eta}_i$ modulo $p^n.$ We also abusively denote by $m_i$ and $m_i'$ the images of $\tilde{m}$ and $\tilde{m}'$ in $\calM_{\frakX_n}$ and $\calM_{\frakX_n'}$ respectively. This way, we have
$$
\calP'_n=\Ox_{\frakX_n'} \langle \xi_1,\hdots,\xi_d \rangle,\quad \calP_n=\Ox_{\frakX_n}\langle \eta_1,\hdots,\eta_d\rangle,\quad \omega^1_{\frakX_n/\frakS_n}=\bigoplus_{i=1}^d\Ox_{\frakX_n}\op{dlog}m_i,\quad \omega^1_{\frakX_n'/\frakS_n}=\bigoplus_{i=1}^d\Ox_{\frakX_n'}\op{dlog}m_i'.
$$
Let $\ov{\calI}$ and $\ov{\calI}'$ be the PD ideals of $P_n$ and $P_n'$ respectively. For an $\Ox_{\frakX_n,\upmu}$-module $\calE,$ an $\Ox_{\frakX_n',\upmu}$-module $\calE'$ and a local section $s\in\Gamma(U,\upmu)$ over an étale $\frakX_n$-scheme $U,$ we consider the projections
$$
\pi_{s}:\calE_{s}\otimes_{\Ox_{U}}\calP_{n|U} \ra \calE_s\otimes_{\Ox_U} \left (\calP_{n|U} / \ov{\calI}^{[2]} \right ),\quad \pi'_{s}:\calE'_{s}\otimes_{\Ox_{U'}}\calP'_{n|U'} \ra \calE'_s\otimes_{\Ox_{U'}} \left ( \calP'_{n|U'}/\ov{\calI}'^{[2]} \right ).
$$
Let $(\calE',\epsilon')$ be an object of $\boldsymbol{\op{MHS}(\calB_{\frakX_n/\frakS_n})},$ $(\calE',\theta')$ its image in $\boldsymbol{p\op{\text{-}MIC}(\calB_{\frakX_n/\frakS_n})},$ $(\calE,\epsilon)=\Psi(\calE',\epsilon')$ and $(\calE,\nabla')$ (resp. $(\calE,\nabla)$) the image of $(\calE',\theta')$ (resp. $(\calE,\epsilon)$) in $\boldsymbol{\op{MIC}(\calA_{\frakX_n})}.$ We have to prove that $\nabla=\nabla'.$
Let $s\in \Gamma(U,\calM_{\frakX_n}^{\mathrm{gp}})$ for some étale $\frakX_n$-scheme $U$ and $\ov{s}$ its image in $\upmu.$ The $p$-connection
$$\theta'_{\ov{s}}:\calE'_{\ov{s}} \ra \calE'_{\ov{s}}\otimes_{\Ox_{U'}}\omega^1_{U'/\frakS_n}$$
is given, for all local sections $x$ of $\calE'_{\ov{s}},$ by
\begin{alignat*}{2}
\theta'_{\ov{s}}(x) &= \sum_{i=1}^d\theta'_{\ov{s}}(\partial_{\epsilon_i})(x)\otimes \op{dlog}m'_i,
\end{alignat*}
where $\theta'_{\ov{s}}(\partial_{\epsilon_i})$ is given in \ref{def227}.
By \eqref{explicitconnimage}, the connection
$$\nabla'_{\ov{s}}:\calE_{\ov{s}}=(\calA_{\frakX_n}\circledast_{\calB_{\frakX_n/\frakS_n}}\calE')_{\ov{s}} \ra \calE_{\ov{s}}\otimes_{\Ox_U}\omega^1_{U/\frakS_n}$$
is induced by
\begin{equation}\label{eq4421}
\begin{array}[t]{clc}
\bigoplus\limits_{\ov{t}+\ov{r}=\ov{s}}\calA_{\frakX_n,\ov{t}}\otimes_{\Ox_U'}\calE'_{\ov{r}} & \ra & \bigoplus\limits_{\ov{t}+\ov{r}=\ov{s}}\left ( \calA_{\frakX_n,\ov{t}}\otimes_{\Ox_U'}\calE'_{\ov{r}} \right ) \otimes_{\Ox_U}\omega^1_{U/\frakS_n} \\
(ae_t)\otimes x & \mapsto & \sum\limits_{i=1}^d (ae_t)\otimes \theta'_{\ov{r}}(\partial_{\epsilon_i})(x)\otimes (p^{-1}dF)(\op{dlog}m'_i)+d_{\calA_{\frakX_n},\ov{t}}(ae_t)\otimes x,
\end{array}
\end{equation}
where $a$ and $x$ are local sections of $\Ox_U$ and $\calE'_{\ov{r}}$ respectively.
Now let us look at $\epsilon.$ By definition of the functor $\Psi,$
$$\epsilon=\epsilon_{\calA_{\frakX_n}}\circledast_{\calP'^{\calB}_n}\epsilon'.$$
The morphism $\epsilon_{\ov{s}}$ is then induced by
$$
\begin{array}[t]{clc}
\calP_{n|U}\otimes_{\Ox_U}\bigoplus\limits_{\ov{t}+\ov{r}=\ov{s}}\calA_{\frakX_n,\ov{t}}\otimes_{\Ox_{U'}}\calE'_{\ov{r}} & \ra & \bigoplus\limits_{\ov{t}+\ov{r}=\ov{s}} \left ( \calA_{\frakX_n,\ov{t}}\otimes_{\Ox_{U'}}\calE'_{\ov{r}} \right ) \otimes_{\Ox_U}\calP_{n|U} \\
1\otimes (ae_t)\otimes x & \mapsto & \epsilon_{\calA_{\frakX_n},\ov{t}}(1\otimes ae_t)\otimes \epsilon'_{\ov{r}}(1\otimes x)
\end{array}
$$
The connection
$$\nabla_{\ov{s}}:\calE_{\ov{s}}=(\calA_{\frakX_n}\circledast_{\calB_{\frakX_n/\frakS_n}}\calE')_{\ov{s}} \ra \calE_{\ov{s}}\otimes_{\Ox_U}\omega^1_{U/\frakS_n}$$
is then induced by
\begin{equation}\label{eq4422}
\begin{array}[t]{clc}
\bigoplus\limits_{\ov{t}+\ov{r}=\ov{s}}\calA_{\frakX_n,\ov{t}}\otimes_{\Ox_{U'}}\calE'_{\ov{r}} & \ra & \bigoplus\limits_{\ov{t}+\ov{r}=\ov{s}}\left ( \calA_{\frakX_n,\ov{t}}\otimes_{\Ox_{U'}}\calE'_{\ov{r}} \right ) \otimes_{\Ox_U}\omega^1_{U/\frakS_n} \\
(ae_t)\otimes x & \mapsto & \pi_{\ov{s}}(\epsilon_{\calA_{\frakX_n},\ov{t}}(1\otimes ae_t)\otimes \epsilon'_{\ov{r}}(1\otimes x))-(ae_t)\otimes x\otimes 1
\end{array}
\end{equation}
We conclude by the following computations
\begin{alignat*}{2}
\pi_{\ov{t}}(\epsilon_{\calA_{\frakX_n},\ov{t}}(1\otimes ae_t)) &= (ae_t)\otimes 1+\sum_{i=1}^dd_{\calA_{\frakX_n},\ov{t}}(\partial_{\epsilon_i})(ae_t)\otimes \eta_i, \\
\pi'_{\ov{r}}(\epsilon'_{\ov{r}}(1\otimes x)) &= x\otimes 1+\sum_{i=1}^d\theta'_{\ov{r}}(\partial_{\epsilon_i})(x)\otimes \xi_i.
\end{alignat*}
It follows that, modulo $\ov{\calI}^{[2]},$ we have
$$
\epsilon_{\calA_{\frakX_n},\ov{t}}(1\otimes ae_t)\otimes \epsilon'_{\ov{r}}(1\otimes x)-(ae_t)\otimes x\otimes 1=
\sum_{i=1}^dd_{\calA_{\frakX_n},\ov{t}}(\partial_{\epsilon_i})(ae_t)\otimes x\otimes \eta_i
+\sum_{i=1}^d (ae_t) \otimes \theta'_{\ov{r}}(\partial_{\epsilon_i})(x)\otimes \varphi \left (\xi_i \right ).
$$
The equality of \eqref{eq4421} and \eqref{eq4422} follows then from \eqref{phidF}.
\end{proof}

\begin{theorem}\label{Phiequiv}
The functor $\Phi,$ defined in \eqref{Phiind}, induces an equivalence of categories
$$
p\text{-}\boldsymbol{\op{MIC}}^{\mathrm{qn}}(\calB_{\frakX_n/\frakS_n}) \isom  \boldsymbol{\op{MIC}}^{\mathrm{qn}}(\calA_{\frakX_n}).
$$
\end{theorem}

\begin{proof}
The fact that $\Phi$ preserves quasi-nilpotence follows from \ref{prop443}. The fact that it is an equivalence follows then from \ref{Psiequiv} and the commutative diagram
\begin{equation*}
\begin{tikzcd}
\boldsymbol{\op{MHS}(\calB_{\frakX_n/\frakS_n})} \ar{r}{\Psi} \ar[swap,sloped]{d}{\sim} & \boldsymbol{\op{MHS}(\calA_{\frakX_n})}\ar[sloped]{d}{\sim} \\
p\text{-}\boldsymbol{\op{MIC}}^{\mathrm{qn}}(\calB_{\frakX_n/\frakS_n}) \ar{r}{\Phi} & \boldsymbol{\op{MIC}}^{\mathrm{qn}}(\calA_{\frakX_n}).
\end{tikzcd}
\end{equation*}
\end{proof}

\bibliographystyle{plain}
\bibliography{Indexed_local_Cartier_transform_modulo_p_n}

\end{document}